\documentclass[12pt,reqno]{amsart}
\usepackage[foot]{amsaddr}    
\usepackage{comment,listings,framed,graphicx,url,xcolor,wrapfig}
\usepackage{fullpage,caption} 
\usepackage{url}
\usepackage{cases,float,bbm}
\newcommand*\diff{\mathop{}\!\mathrm{d}}

\newcommand{\kn}[3]{#1_{\mathrm{KIN#2}}^{#3}}
\newcommand{\okn}[1]{\overline{k_{#1}}}
\newcommand{\Em}{E_*^{-}}
\newcommand{\Ep}{E_*^{+}}
\newcommand{\Dm}{D^0_{-}}
\newcommand{\Dp}{D^0_{+}}
\newcommand{\qbar}{\overline{Q}}
\newcommand{\Ubar}{\overline{U}}
\newcommand{\sbar}{\overline{S}}
\newcommand{\chibar}{\overline{X}}
\newcommand{\psibar}{\overline{\Psi}}
\newcommand{\abs}[1]{\left\vert\,#1\,\right\vert}
\newcommand{\wt}{\widetilde}

\newcommand{\R}{\mathbb{R}}

\newtheorem{example}{Example}
\newtheorem{remark}{Remark}
\newtheorem{lemma}{Lemma}
\newtheorem{proposition}{Proposition}

\usepackage{soul}

\newcommand{\sign}{\mathrm{sgn}}
\newcommand{\norm}[2]{\|#1\|_{#2}}
\newcommand{\gnorm}[1]{\|#1\|_{\Delta,1}}
\newcommand{\ba}{\begin{eqnarray}}
\newcommand{\ea}{\end{eqnarray}}
\newcommand{\bsub}{\begin{subequations}}
\newcommand{\esub}{\end{subequations}}
\newcommand{\bas}{\begin{eqnarray*}}
\newcommand{\eas}{\end{eqnarray*}}
\newcommand{\Meta}{{X}}
\newcommand{\dt}[1]{\frac{d#1}{dt}}
\newcommand{\ddt}{\frac{d}{dt}}
\newcommand{\lchistar}{L_{\chi^*}}
\newcommand{\myre}{{\mathcal{R}}}

\newcommand{\mpunit}[1]{\mathrm{[#1]}}

\newcommand\myskip[1]{}

\newcommand{\mpm}[1]{}

\newcommand{\mps}{\mpm{\bf Shortened. }}
\newcommand\mpskip[1]{\mpm{Removed}}

\definecolor{csgreen}{rgb}{0.0, 0.8, 0.6}

\usepackage{etoolbox}
\makeatletter
\patchcmd\maketitle
  {\uppercasenonmath\shorttitle}
  {}
  {}{}
\patchcmd\maketitle
  {\@nx\MakeUppercase{\the\toks@}}
  {\the\toks@}
  {}
  {}{}
\patchcmd\@settitle
  {\uppercasenonmath\@title}
  {}
  {}{}
\patchcmd\@setauthors
  {\MakeUppercase{\authors}}
  {\authors}
  {}{}
\renewcommand{\email}[2][]{%
\ifx\emails\@empty\relax\else{\g@addto@macro\emails{,\space}}\fi%
\@ifnotempty{#1}{\g@addto@macro\emails{\textrm{(#1)}\space}}%
\g@addto@macro\emails{#2}%
}
\makeatother

\title[Stability for Transport in Hydrate Zone]{Stability of a numerical scheme for methane transport in hydrate zone under equilibrium and non-equilibrium conditions}
\author{Malgorzata Peszynska and Choah Shin}
\email[M. Peszynska]{mpesz@oregonstate.edu}
\email[C. Shin]{shinc@oregonstate.edu}
\address{Department of Mathematics, Oregon State University, Corvallis, OR 97330}
\thanks{This research was partially supported by NSF DMS-1522734 ``Phase transitions in porous media across multiple scales'' and DMS-1912938 ``Modeling with Constraints and Phase Transitions in Porous Media'', NSF IRD plan 2019-20, and Martin-O'Neill COS Fellowship.}
\date{Article accepted for publication by the Journal of Computational Geosciences on March 16, 2021}

\begin{document}
\maketitle
\begin{abstract}
In this paper we carry out numerical analysis for a family of {simplified} gas transport models with hydrate formation and dissociation  in subsurface, in equilibrium and non-equilibrium conditions. {These models are adequate for simulation of hydrate phase change} at basin and at shorter time scales, {but the analysis does not account directly for the related effects of evolving hydraulic properties}. 

To our knowledge this is the first analysis of such a model. It is carried out for the transport steps while keeping the pressure solution fixed. We frame the transport model as conservation law
with a non-smooth space-dependent flux function; the kinetic model {approximates} this equilibrium. We prove weak stability of the upwind scheme applied to the regularized conservation law. We illustrate the model, confirm convergence with numerical simulations, {and illustrate its use for some relevant equilibrium and non-equilibrium scenarios.}

\end{abstract}
\section{Introduction}

In this paper we analyze a computational model for transport of methane in hydrate zone at the equilibrium and kinetic time scales. Our interest in methane hydrate comes from collaborations with geophysicists who aim to explain and predict hydrate deposits found in nature. 
 The simulations of hydrate evolution have been carried out by many researchers
 including in \cite{DB2001,DB2003,LF05,LF2007,LF11,MKTW08,NRuppel03,PHTK,PMHT,XuRuppel99,Xu04}; however, our paper appears to be the first one to analyze  the numerical schemes. 

\medskip
Methane hydrate, also known as ``Ice That Burns'' is an ice-like crystalline substance made of methane molecules enclosed in a cage made by water molecules. Methane hydrate is abundant in deep sub-sea sediments whenever favorable conditions of high pressure, low temperature, and large supply of methane hold. Methane hydrate is also found in the Arctic below permafrost. 

To explain the presence and shape of hydrate deposits found in nature, as well as to understand the methane fluxes as a response to the climate change, various simulations were carried out; see, e.g.,  \cite{Garg08,HP18,LF2007,PHTK,PMHT,TT04,WH-svalbard-2018,Xu04,XuRuppel99}. These simulations are typically carried out at the {\em basin time scales} of several kilo-years or at least years or months. The presence of hydrate is explained with a postulate of supply of gas from deep Earth sources or by existence of biogenic sources of methane such as microbial species; see \cite{LF2007,PHTK,TT04}. Recent studies focus also on the dissociation of hydrate deposits in response to the  environmental conditions such as an increase in average 
9
 temperatures and address the impact of hydrate on the balance of greenhouse gases \cite{berndt2014,DB2001,Dickens03,Hunter13,FL14} at the time scale of years or decades. 

Hydrate has also been  evaluated as a potential energy source \cite{Intercomparison,Moridis2002,ruppel2007tapping}. In particular, in the pilot projects in Japan and Alaska
\cite{NETLmh,Moridis2002,MHNETL}, the recovery of methane is enabled by lowering the pressure in the wells which triggers hydrate dissociation and release of large amounts of gas.   A similar mechanism contributes to hazard while drilling \cite{HenCle99-stefan,HenCle99,Sloan}, with the characteristic time scale of days. 

\medskip
{\bf Overview:} In this paper we analyze a discrete model of hydrate formation and dissociation which describes methane transport {by advection and diffusion}  coupled to phase behavior {in equilibrium and non-equilibrium conditions and two-phase liquid-hydrate conditions, and which treats phase behavior in a sequential way with macro-time steps}. 
The equilibrium model is {a simplified version} of the comprehensive model \cite{LF2007} presented earlier in \cite{PHTK,PMHT}. {Kinetic models from the literature \cite{EKDB1987,GHW2015,GWH2016,KB1987,WH-svalbard-2018} are formulated for the general context of three-phase equilibria; our model resembles these but covers both unsaturated as well as saturated conditions in liquid-hydrate conditions.}
The scheme combines finite volume spatial discretization with implicit-explicit time discretization, {and uses the formal mathematical framework of multivalued graphs. This framework for the equilibrium model is equivalent to variable switching as we demonstrated in \cite{GMPS}. Our analysis of the kinetic model with this framework supports the understanding of the equilibrium as the limit of kinetic model under fast reaction rates.} 

Our main contribution is the analysis of {numerical stability of the advective model in equilibrium and non-equilibrium, as well as demonstration of the convergence of the scheme at the rate $O(\sqrt{h})$ common for generic scalar conservation laws.   The analysis applies to the transport in liquid-hydrate zone under various simplifying assumptions including long time range, close to geothermal temperature distribution, close to hydrostatic pressure distribution, modest gas supply, and constant salinity. We make a-priori assumptions on the data and illustrate the sensitivity of the model to data as guided by the analysis}. In particular, our analyses explain  the feature of discontinuous hydrate lenses observed in nature and exacerbated in heterogeneous sediments. The scheme for the kinetic model is shown to be robust across the saturated and unsaturated conditions. 

To our knowledge, our analysis is the first of this kind for advective transport in either the equilibrium and kinetic setting for hydrate models.  {Our analysis applies only to the 
simplified model with which we simulate hydrate phase change, but does not directly apply to the possibly strongly coupled effects like evolving hydraulic properties which are critical for simulation of hydrate evolution and recycling.} Our results are {therefore} the first step towards  {future work on the } analysis of schemes for more comprehensive models. {The paper also includes auxiliary supporting results which can be used in a more general context.} 

\medskip
 
{\bf Plan of the paper:} We start with {auxiliary notation on evolution with multivalued monotone graphs} in Sec.~\ref{sec:monotone}. 
{In Sec.~\ref{sec:model} we provide details of the transport model, starting with the equilibria, and in Sec.~\ref{sec:kinmodel} we describe the kinetic models. In Sec.~\ref{sec:LAG} we define the numerical schemes. In Sec.~\ref{sec:EQ} we analyze the numerical scheme for the equilibrium model and in Sec.~\ref{sec:KIN} we analyze the scheme for the kinetic model. Finally, in Sec.~\ref{sec:results} we present numerical results: we present convergence studies in the case covered by the theory as well as simulation results in  realistic examples including the comparison of equilibrium and kinetic models. We conclude in Sec.~\ref{sec:conclusions}.} The Appendix in Sec.~\ref{sec:appendix} provides extensive auxiliary results. 

\section{Notation and background for ODEs with monotone graphs.}
\label{sec:monotone}
We recall here the notation and a few elements of the mathematical framework of evolution equations with monotone multivalued graphs on $\R$ to extend what is known for the initial value problem $\dt a + G(a)=f, a(0)=a^0$ when $G:\R \to \R$ is a monotone increasing function. {These extensions are useful for modeling phase equilibria and kinetic schemes}. We refer to the comprehensive details on the general abstract Hilbert space setting and monotone multivalued operators provided, e.g., in \cite{Brezis73,Show-monotone}. {We need the notation and basic properties in our estimates and analyses}.

We recall that for a relation (graph) $G\subset \R \times \R$, its $\mathrm{domain}(G)=\{a: {\exists b:} (a,b) \in G\}\subset \R$,  and the inverse $G^{-1}=\{(b,a): (a,b)\in G\}$. The graph $G$ is monotone if $\forall (a_1,b_1), (a_2,b_2) \in G$ we have $(b_2-b_1)(a_2-a_1)\geq 0$.  If $(a,b) \in G$, and $G$ is monotone multivalued, we will write $b \in G(a)$ to denote some selection $b$ out of the set $G(a)$. {This selection is not unique, hence the symbol ``$\in$''. }
Further, $G$ is maximal monotone if $I+G$ is onto $\R$. For a maximal monotone $G$ and $\lambda>0$, the resolvent  
\ba
\label{eq:re}
\myre^G_{\lambda}=(I+\lambda G)^{-1}
\ea
is a contractive function, and the solution to $a+\lambda G(a) \ni f$ is unique and given by $a= \myre^G_{\lambda}(f)$. 

{\bf Evolution ODE with graph:} {The resolvent $\myre^G_{\lambda}$ helps to define the solution to an evolution problem}  
\ba
\label{eq:odegraph} \dt a + G(a) \ni f;\;\; a(0)=a^0,
\ea
where $f\in L^1(0,T)$ is some given input and where $a^0\in \mathrm{domain}(G)\subset \R$ is some initial data. 
The $C^0$ solution $a(t)$ {to} \eqref{eq:odegraph} is defined as the limit as $\tau \to 0$ of the fully implicit finite difference step function solutions $a^n \approx a(t^n)$, with $t^n=n\tau$,  to the inclusion
\ba
\label{eq:fd}
\frac{a^n-a^{n-1}}{\tau}+G(a^n) \ni f^n, \; n\geq 1.
\ea
{In spite of the symbol $\in$}, the step solution $a^n\in \mathrm{domain}(G)$ to \eqref{eq:fd} is uniquely defined 
$a^n=\myre^G_{\tau} (a^{n-1}+\tau f^n)$.
Once we know $a^n$, the actual selection $G(a^n)=f^n-\tfrac{a^n-a^{n-1}}{\tau}$ is given uniquely from \eqref{eq:fd}.   

{In this paper we use various single-valued approximations $G_{\lambda} \approx G$} which are maximal monotone when $G$ is. One is the Yosida approximation $G_{\lambda}=\frac{1}{\lambda}(I-\myre^G_{\lambda})$ which provides another way to define the solution $a(t)$ to \eqref{eq:odegraph} as the limit as $\lambda \to 0$ of $a_{\lambda}(t)$, the family of solutions to the ODE $\dt {a_{\lambda}} + G_{\lambda}(a_{\lambda})= f$. 

\medskip

{\bf Evolution system with graphs:} 
{In addition to \eqref{eq:odegraph}, we consider the following system on $\R \times \R$}
\ba
\label{eq:odesys}
\ddt{a} =Q, \ddt{b}=-Q;\;\; a(0)=a^0,b(0)=b^0.
\ea
{We are interested in the case when $Q(a,b)\in b-G(a)$ is multivalued with $G(\cdot)$ is monotone. Here the first and the second equations have similar properties to \eqref{eq:odegraph} but are coupled.  Adding the two equations leads to  $a(t)+b(t)=\mathrm{const}=a^0+b^0$. With the abstract theory from \cite{Show-monotone}, it is easy to show that the system \eqref{eq:odesys} is well-posed in $\R \times \R$. We also see that the solutions $(a(t),b(t))$ evolve towards some $(a^{\infty},b^{\infty})$ which is at the intersection of $G$ with the manifold $a+b=a^0+b^0$}.

\medskip
{\bf Special graphs used in this paper:} 
The graph $\sign(x)$ assigns $-1$ to $x<0$, $1$ to $x>0$, and the set $[-1,1]$ to $x=0$, and we write $\sign(x)=(-\infty,0)\times\{-1\} \cup \{0\} \times [-1,1] \cup (0,\infty) \times \{1\}$. This graph $\sign(x)$ is distinct from the single valued discontinuous function $\sign_0(x)$ which agrees with $\sign(x)$ for $x \neq 0$ but assigns $0$ to $x=0$. 
The Heaviside graph $H(x)=\frac{1}{2}(1+\sign(x))$ assigns $0$ to $x<0$, 1 to $x>0$, and the set $[0,1]$ to $x=0$. We also use $x_+=\mathrm{max}(0,x)$. 

\section{Transport Model under Equilibrium  and Kinetic Phase Constraints}
\label{sec:model}
Methane is present in sub-ocean sediments and Arctic regions due to biogenic sources and from upward fluxes from the deeper Earth layers \cite{LF2007}. It is transported by diffusion and advective fluxes, and can be present in liquid, gas or solid phases.  The partition of methane component between phases depends on the pressure and temperature and on the amount of  methane component. 
With small amounts of methane and at large depths (i.e., large pressures), methane is dissolved in the aqueous (liquid) brine phase denoted by $l$. With larger amounts of methane and at low temperatures, the solid phase made of methane and water in fixed proportions precipitates; this solid phase denoted by $h$ is called methane hydrate (clathrate or methane ice). At higher temperatures, the hydrate phase is not stable and free gas phase forms. A typical distribution of phases in sub-ocean sediments is that the solid hydrate phase is stable and present at low temperatures; specifically, this occurs above the so-called Bottom Hydrate Stability Zone (BHSZ). In turn, below BHSZ, only the gas phase is stable. Phase equilibria represent the tendency of a system to maintain low energy, and correspond to the most stable distribution of components between phases. 

The phase distributions may not always follow equilibria; this is common at short time scales, e.g., after seismic events which alter the distribution of gases and sediments, or during production of gas from subsurface. Our work in this paper addresses equilibrium models as well as certain selected scenario of non-equilibria. In this paper we do not account for the presence of free gas such as  ex-solved gas or from buoyant gas travelling upwards. The analysis of a model involving gas is the subject of current work. 

In this section we describe a model accounting for methane transport above BHSZ in a porous reservoir $\Omega \subset \R^d$, over time $t$, either under the assumption of phase equilibria or allowing  non-equilibria. For phase equilibria we follow closely the comprehensive model in \cite{LF2007} in the hydrate zone which we simplified in \cite{PHTK,PMHT} for the purposes of efficient simulations; we used real reservoir data and experimental observations  to motivated the reduced model; see also our well-posedness analysis of the reduced model in  \cite{GMPS,PSW}. 

Our objectives in this paper are to (i) provide an analysis of  the finite volume discretization of the reduced equilibrium model from \cite{PHTK,PMHT}, which we accomplish for a regularization of the original model. Our analysis also  (ii) explains  the shape of hydrate deposits observed in nature as well as the hydrate saturation profiles found in computational simulations. Our main result is  (iii) the study of a  new non-equilibrium model for hydrate in under-saturated and over-saturated conditions above BHSZ which extends literature, and for which we are able to prove rigorous numerical stability.

Below we briefly recall the equilibrium model from \cite{PSW}; it is the same as that in \cite{PHTK,PMHT} under the assumption of constant salinity.  We also motivate and explain the non-equilibrium models which we compare to other kinetic models in the literature. 

\subsection{Assumptions.} We make the following assumptions in the model development.

(A1) The reservoir $\Omega$ is in the hydrate stability zone, i.e.,  only the liquid and hydrate phases {are stable}. 

(A2) Free gas is not present in $\Omega$, {i.e., there is abundant water present for hydrate formation}; see \cite{LF2007}

\noindent
{In addition, the following assumptions are made for the sake of presentation and analysis, but are not needed for the computational model or for simulations.} 

(A3) Liquid and hydrate phases are incompressible.

(A4) The sediment is rigid, and the porosity $\phi(x,t) \approx \phi(x)$ is fixed. 

(A5) Salinity $\chi_{lS}(x,t)=\chi_{lS}^{sw}=\mathrm{const}$ equals seawater salinity $\chi_{lS}^{sw}$. 

\subsection{Mass conservation equations}
Consider $(x,t)$ at a point $x\in \Omega$ and time $t>0$. We denote the mass fraction of methane in the liquid phase by $\chi(x,t)$ and the volume fraction of hydrate by $S(x,t)$.  With (A1) and (A2), 
the total mass density of methane at $(x,t)$ is $\rho_l u=(1-S)\rho_l\chi+ S\rho_h\chi_{Mh}$ with $\rho_l,\rho_h$ denoting mass density of brine and of hydrate, respectively, and $\chi_{Mh}$ denoting the mass fraction of methane in the hydrate phase, which is a known fixed constant.  The mass conservation equation for methane component in porous sediment of porosity $\phi$ is thus 
\bsub
\label{eq:massMC}
\ba
\label{eq:massM}
\partial_t\left(\phi \rho_l u\right)+\nabla \cdot( q\rho_l \chi) -\nabla \cdot (\rho_l {d_m} \nabla \chi)=F_M.
\ea
Here {$d_m$} is diffusivity and $q$ is the Darcy flux defined below. 
Also, $F_M$ accounts for methane sources, e.g., biogenic production of methane by microbes, or sink terms relevant for the production scenarios. We also rewrite 
\ba
\label{eq:udef}
u&=&(1-S)\chi+RS=\chi+S(R-\chi).
\ea
\esub
Here $R = \frac{\rho_h \chi_{Mh}}{\rho_l} \approx 
\chi_{Mh}$ since $\rho_h \approx \rho_l$. In practice, with  the values reported in  \cite{PHTK}; $R = 0.1203$ kg/kg,  we have $\rho_l=1030$ kg/$m^3$, $\rho_h=925$ kg/$m^3$, while ${\max_x \chi^*(x)} \approx 2.4\times 10^{-3}$ kg/kg for the case of UBGH2-7 as given in \cite{PHTK}.
From now on we will assume 
\ba
\label{eq:Rassum}
0 < \chi^*<R.
\ea

With the single equation \eqref{eq:massMC} involving two variables $u$ and $S$ or $\chi$ and $S$, we close the system by either assuming equilibrium conditions binding $\chi$ and $S$, or setting up a non-equilibrium model which evolves towards the equilibrium. We describe the equilibrium model in Sec.~\ref{sec:EQp} and the non-equilibrium models in Sec.~\ref{sec:kinmodel}.

\medskip
The model \eqref{eq:massMC} is complemented with the mass conservation equation for water component whose concentration is $1-\chi$ in liquid phase and $1-\chi_{Mh}$ in the hydrate phase. \mps The model for water mass conservation is
\ba
\label{eq:massW}
\partial_t\left(\phi\left[(1-S)\rho_l(1-\chi)+ S\rho_h(1-\chi_{Mh})\right]\right)+\nabla \cdot(q \rho_l (1-\chi)) =0.
\ea

\subsection{Pressure equation}
The pressure equation follows by adding \eqref{eq:massM} with \eqref{eq:udef} and \eqref{eq:massW}; with Darcy's law we obtain
\begin{subequations}
\label{eq:pressure}
\ba
\label{eq:press}
\partial_t\left(\phi\left[(1-S)\rho_l+ S\rho_h\right]\right)+\nabla \cdot(q \rho_l) -\nabla \cdot (\rho_l {d_m} \nabla \chi)=F_M,
\\
\label{eq:Darcy}
q=-\frac{K}{\mu_l}(\nabla P-\rho_l G\nabla d),
\ea
\end{subequations}
with Darcy flux $q$, pressure $P(x,t)$, permeability $K$, liquid phase viscosity $\mu_l$, and depth $d(x)$. Here 
$K=K(x;S)$
depends on the presence of hydrate in the pore-space, with empirical data, e.g., in  \cite{LF2007}.\mpskip{see also the multiscale models based on porescale simulations in \cite{PTISW}.} Typically $K(\cdot;S)$ decreases with $S$, and the porous matrix is plugged up when hydrate saturation is close to $1$.  

At large time scales such as in basin modeling the pressure follows distribution close to hydrostatic with $q=0$. Otherwise there can be gas fluxes with $q \neq 0$, e.g., from deep in the Earth's crust upwards, and we must solve \eqref{eq:pressure} under given boundary conditions. One practical scenario is when $q$ is given at the bottom of the reservoir, and a fixed pressure is known at the top, e.g., from the known height of water column.
Rewriting \eqref{eq:pressure} as
\ba
\label{eq:spress}
\nabla \cdot q = \frac{F_M}{\rho_l} + \nabla \cdot ({d_m} \nabla \chi) + \partial_t \left(\phi S \frac{\Delta \rho}{\rho_l}\right), \;\; \Delta \rho = \rho_l -\rho_h, \ea
allows to study contributions to local variations of $q$. We see that the magnitude of the first and second terms on the right hand side is modest in realistic settings \cite{PHTK}. However, the third term may contributes to the local increase of velocity due to the density difference $\Delta \rho$ whenever $S\uparrow$ {increases} rapidly.

\subsection{Phase equilibria for hydrate crystal formation} 
\label{sec:EQp}
\mps
The formation of a hydrate crystal out of liquid phase  usually involves the processes of nucleation, diffusion of molecules towards the existing cages, and the adsorption of new crystals, see, e.g.,  molecular dynamics simulations in \cite{Sloan,walsh09} at time scales of $10^{-6}\mpunit{s}$. At the reservoir time scales of transport the hydrate formation or dissociation is modeled by an aggregate of the microscopic processes; {one assumes either} an equilibrium presented here or a kinetic model discussed in Sec.~\ref{sec:kinmodel}.

In equilibrium, the hydrate  forms only if the methane concentration $\chi$ in water has  reached its maximum solubility denoted by $\chi^*$. 
When the hydrate crystals form, we have that $S>0$ and $\chi=\chi^*$, the saturated case. When $\chi <\chi^*$, no hydrate exists, and $S=0$.  This is expressed by the constraint 
\ba
\label{eq:phase}
\left\{
\begin{array}{cc}
\chi\leq \chi^*,& S=0,\\
\chi=\chi^*,&S \geq 0.
\end{array}
\right.
\ea
Next we write an explicit formula for the dependence of the total amount of methane $u$ on $(\chi,S)$. At a given $(x,t)$ with $\chi^*(x,t)$ known, we have
\ba
\label{eq:es}
u= \left\{
\begin{array}{cc}
\chi,& u \leq \chi^*, S=0,\\
(1-S)\chi^*+SR,&u\ge \chi^*,S \geq 0.
\end{array}
\right.
\ea
Conversely, given $u$, the equilibrium values $\chi$ and $S$ on $E_*$ are given uniquely
\begin{equation}\label{eq:invert}
    \chi = \min\{\chi^*,u\} \quad \mathrm{ and } \quad S = \frac{(u-\chi^*)_+}{R-\chi^*}.
\end{equation}

The relationship \eqref{eq:invert} can be used at any $(x,t)$ to get the unique values $\chi(x,t)$ and $S(x,t)$ from $u(x,t)$. 
{The quantity $\chi^*=\chi^*(P,T,\chi_{lS})$ depends on the pressure $P$ and temperature $T$, and salinity $\chi_{lS}$. This dependence is resolved sequentially in our computational model: {over some macro time step $P,T$ are kept constant and $\chi^*$ depends on $x$ only}. This is discussed in detail later in Sec.~\ref{sec:macro}.
}

\subsubsection{Equilibrium model with multivalued graphs}

\begin{figure*}
\centering
\includegraphics[width=.4\textwidth]{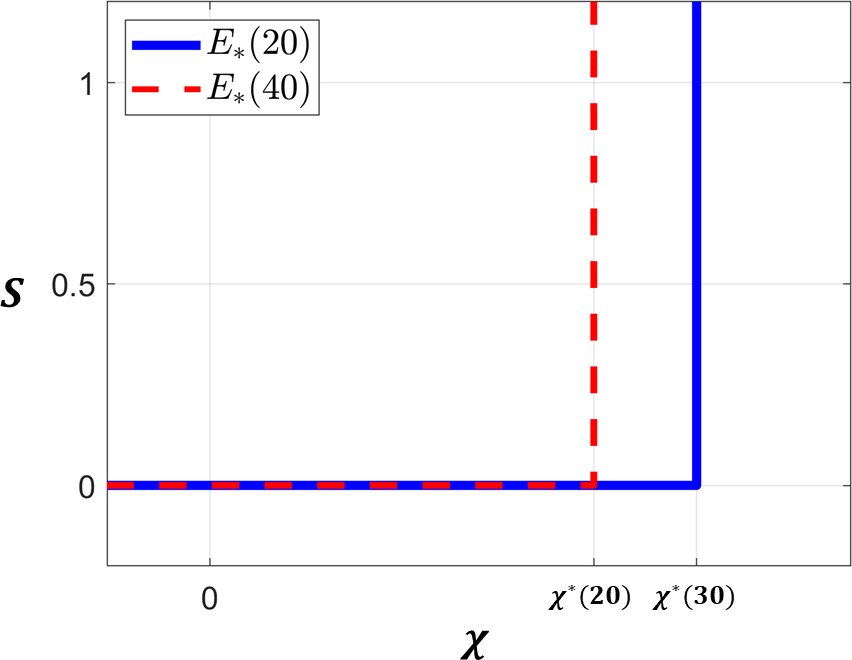}
\includegraphics[width=.4\textwidth]{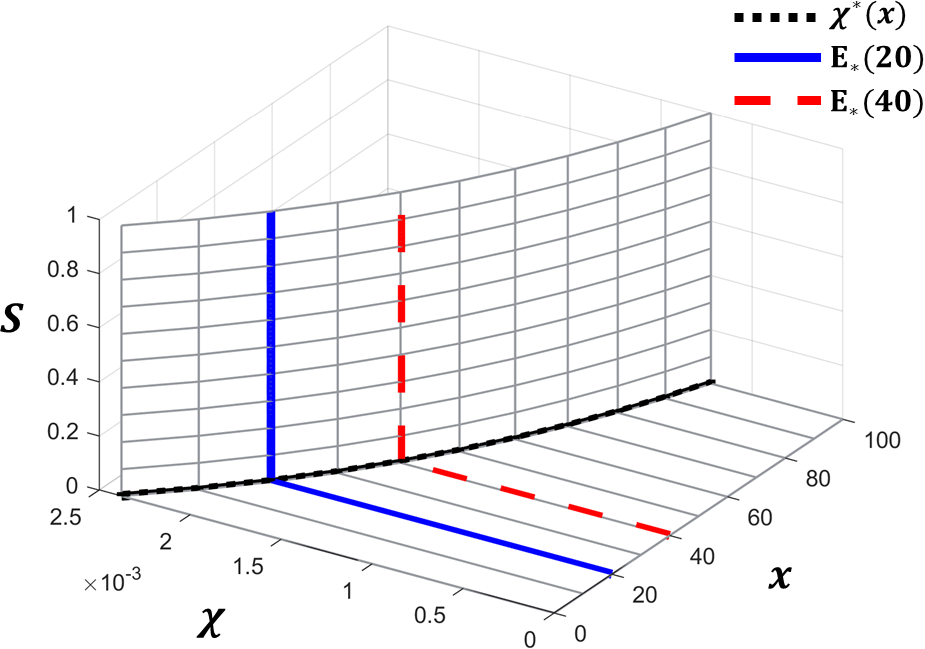}
\caption{Illustration of graph $E_*$ with data from Ulleung Basin case UBGH2-7 \cite{PHTK}, where $\chi^*(x) \approx 0.0024e^{-0.012x}$. Left: the portion of $E_*=E_*(x)$ for a fixed $x$. Right: multivariate view of $(\chi,S) \in E_*(x)$.
\label{fig:Estar_pm}}
\end{figure*}

The formulas \eqref{eq:phase}--\eqref{eq:invert} are simple and explicit. For the needs of the kinetic model to be defined,  we note that the variables $(\chi,S)$ ``live'' on the graph $E_*=\Em\cup \Ep$ defined as 
\ba
\label{eq:phasee}
(\chi,S) \in E_* = (-\infty,\chi^*] \times \{0\} \cup \{ \chi^*\} \times [0,\infty),
\ea
with $\Em=(-\infty,\chi^*] \times \{0\},\;\; \Ep=
\{ \chi^*\} \times [0,\infty)$, as illustrated in Fig.~\ref{fig:Estar_pm}.

The inverse graph $W_*=E_*^{-1}$ is
\bas
 W_*(S) \ni \chi \equiv (S,\chi) \in W_* = \{0\} \times (-\infty,\chi^*]   \cup [0,\infty) \times \{ \chi^*\}.
\eas
It is easy to see that both $E_*$ and $W_*$ are maximal monotone. In what follows we write $S \in E_*(\chi)$ or $\chi \in W_*(S)$. These graphs are set-valued, but in evolution models, the particular selection out of these graphs is actually unique,  as discussed in Sec.~\ref{sec:monotone}.

\begin{remark}
\label{rem:uineq}
Not all the points on the graph $E_*$ \eqref{eq:phasee} are physically meaningful. In particular, any reasonable calculated values of concentrations and saturations should satisfy $\chi\geq 0$ and $0 \leq S <1$. We denote by $E_*^0=\{(\chi,S) \in E_*: \chi\geq 0;  S < 1\}$ the physically meaningful portion of $E_*$. In addition, from \eqref{eq:Rassum} and \eqref{eq:es} we see that if $(\chi,S) \in E_*^0$, then $u$ satisfies
\ba
\label{eq:uphys}
0 \leq u(\chi,S)< R.
\ea
Conversely, for any $u$ which satisfies \eqref{eq:uphys}, we have from \eqref{eq:Rassum} that $(\chi,S)$ given by \eqref{eq:invert} satisfies $0\leq \chi \leq \chi^*$ and $0 \leq S < 1$.
\end{remark}

\section{Transport model with Kinetic Phase Constraints}
\label{sec:kinmodel}
Kinetic models are common in geochemistry and chemical engineering  \cite{lasaga,zhang-book} and describe the evolution of a system towards thermodynamic equilibrium from some initial conditions  out of equilibrium, e.g., in the processes of adsorption, phase transitions, and crystal precipitation and dissolution. 

The time scale of hydrate formation or dissociation is on the order of hours or days $O([\mathrm{h}])-O([\mathrm{days}])$ \cite{Sloan,FleBry}. In production scenarios \cite{Islam1994,Ruan2012,Moridis2002,MS2007,YHW2016} this time scale is comparable to that  of the transport processes. Comprehensive subsurface transport simulators including STOMP, TOUGH, PFLOTRAN, GEOS, Geo-COUS implement the complex kinetic exchange model in the applications using depressurization or thermal stimulation to aid methane recovery from hydrate; see e.g., the recent international code comparison studies \cite{IGHCCS2,IGHCCS1} led by DOE/NETL.

For modeling methane in the environment at large spatial scales, e.g., methane flux response to environmental temperature variations or abrupt geological events, some authors use kinetic models \cite{CCC2013,GHW2015,GWH2016,TT04,WH-svalbard-2018}. Finally, some computational models use kinetics rather than equilibria to implement or to approximate phase behavior regardless of the time scale considered \cite{DB2001,RB1997,TT04}. 

A general kinetic model must predict the evolution of all relevant variables  towards an equilibrium from some out-of equilibrium state, and is complemented by other equations which describe the evolution of all of $(T,P,\chi,S)$ towards some equilibrium $(T^{\infty},P^{\infty},\chi^{\infty},S^{\infty})$, starting from some initial $(T^0,P^0,\chi^0,S^0)$. The kinetics is coupled to the transport  and constitutive equations, and would account for the presence of gas phase and capillary effects. 

{In the framework of our reduced model for liquid-hydrate zone we assume $(P,T)$ are fixed over some time interval $(t^{old},t^{new})$ with $t^{new}=t^{old}+\Delta t$. In equilibrium the variables $(\chi(x,t),S(x,t)) \in E_*(x)$ at every $t$. If the $(P,T)$ conditions change at $t^{new}$, and a new $\chi_{new}^*=\chi^*(x,t^{new})$ is given, the variables 
$(\chi,S)$ are out of the equilibrium with respect to the new graph $E_*^{new}$. If $\Delta t$ is really large we can assume they immediately adjust to the new equilibrium. Otherwise we need a kinetic model to describe the evolution of $(\chi,S)$ towards $E_*^{new}$. }

Below we discuss kinetic models for hydrate, starting
with literature review 
and a homogeneous ``batch reactor'' model, which we couple later with a transport model.

\subsection{Kinetic models of hydrate formation: literature background}

{Following \cite{CB2001,EKDB1987,GHW2015,GWH2016,KB1987,YASS1991}, the kinetics of gas-liquid-hydrate phase system involves an exchange term $Q$  proportional to the driving force in the three phase conditions, 
\ba
\label{eq:drivingforce}
\ddt{S}=Q=k (f_g-f_{eq}),
\ea
where $k>0$ is the hydrate formation {or the dissociation rate}, and where $f_g,f_{eq}$ are the local gas fugacity, and the equilibrium fugacity at the given pressure and temperature, respectively. This expression \eqref{eq:drivingforce}  predicts that the hydrate forms when $f_g > f_{eq}$, and dissociates when $f_g < f_{eq}$. In \cite{GHW2015,GWH2016} the authors propose $Q \propto (P-P_{eq})$, with $P_{eq}$ equal the equilibrium pressure for a given fixed $T$, and this approach models a response to the increase or decrease in pressure.  A physically grounded expression for $k$ is complex \cite{Tishchenko2005}.
The rate $k \propto A_s$, the surface area available for the reaction to occur which is proportional to the effective porosity $\phi(1-S)$. In a three phase system the hydrate formation rate $k$ also depends on the availability of water and methane (thus on the gas and aqueous phase saturations $S_g$ and $S_w$); but for hydrate dissociation the rate $k$ depends on availability of hydrate (thus on $S=1-S_g-S_w$). Therefore \eqref{eq:drivingforce} is in general hysteretic; see also \cite{FleBry}. 
From mathematical point of view the presence of $S$ or $(1-S)$ in $k$ keeps the variable $S$ in physically meaningful domain $S\in[0,1]$. Model \eqref{eq:drivingforce} is designed to work in the saturated case when $S^{\infty}>0$, and $\chi^{\infty}=\chi^*$. }

{Our focus in this paper is on liquid-hydrate systems. For these,  according to \cite{CCC2013,DB2001}, the driving force $f_g-f_{eq}$ in \eqref{eq:drivingforce} can be expressed by the difference of methane concentration at the liquid-gas equilibrium and maximum methane solubility $\chi^*$ at three-phase equilibrium state for the given $P$ and $T$. With no free gas, $Q$ becomes 
\begin{eqnarray}
\label{eq:dbkin}
Q=k(\chi-\chi^*),
\end{eqnarray}
similar to that for crystal formation from  saturated or oversaturated mixtures in geochemistry \cite{zhang-book,lasaga}.
In this paper we extend \eqref{eq:dbkin} so it can 
work well across the two-phase saturated as well as in single phase unsaturated conditions when  $S^{\infty}=0$, and $\chi^{\infty} <\chi^*$.  We explain this extension in Sec.~\ref{sec:batch} with further details given in Sec.~\ref{sec:app-KIN}. Our model is robust, also when coupled to the transport model. In the future we hope to extend it to the three phase equilibria extending \eqref{eq:drivingforce}.}

\subsection{Kinetic batch reactor model model for hydrate evolution in liquid-hydrate conditions \mps}
\label{sec:batch}

\begin{figure}
\begin{center}
\includegraphics[width=5cm]{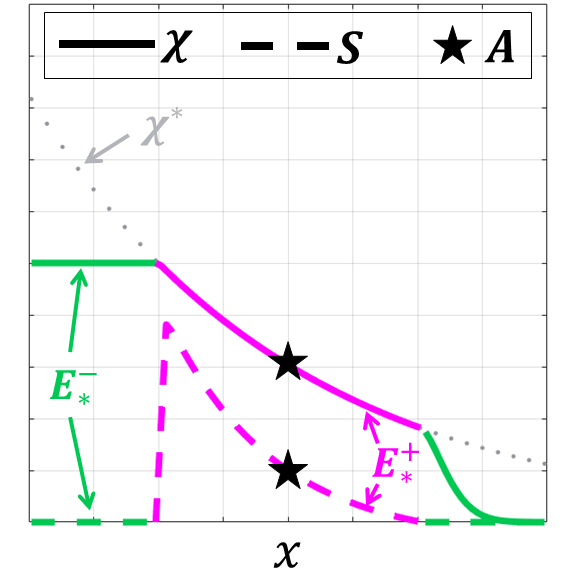}
\\
\includegraphics[width=8cm]{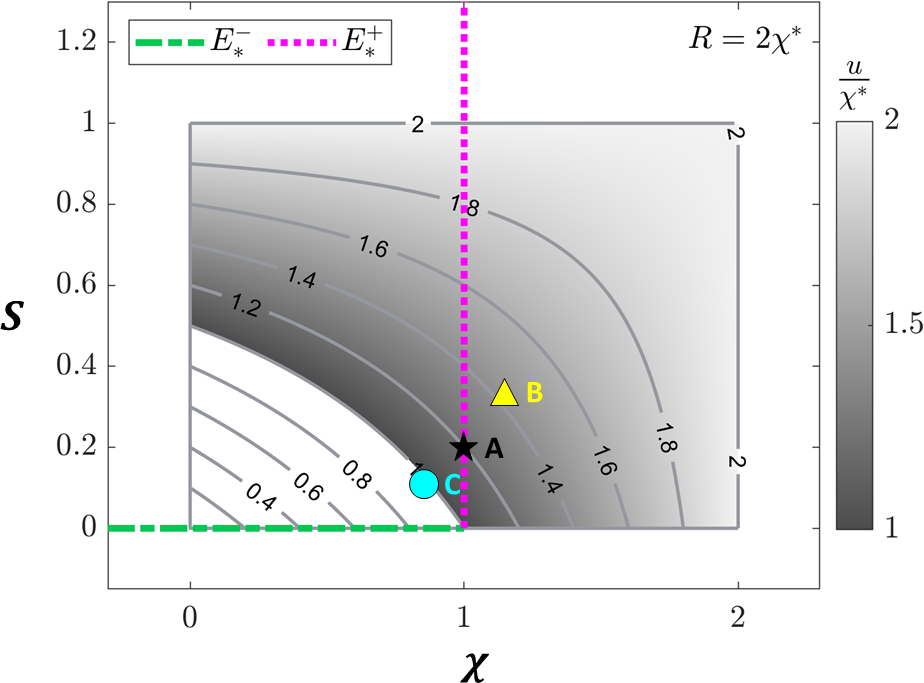}
\\
\includegraphics[width=5cm]{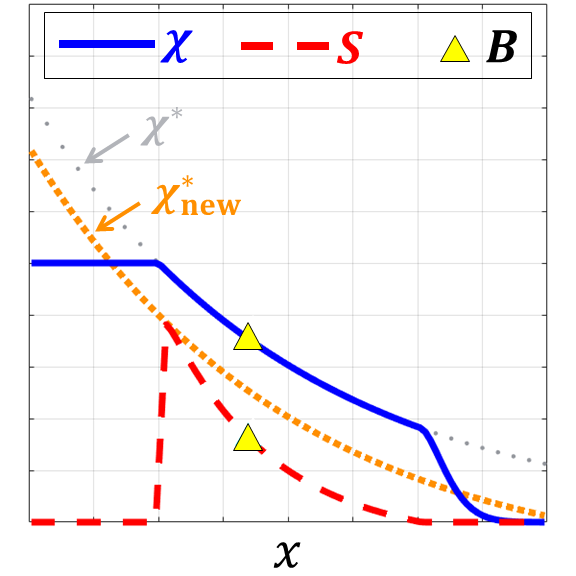}
\includegraphics[width=5cm]{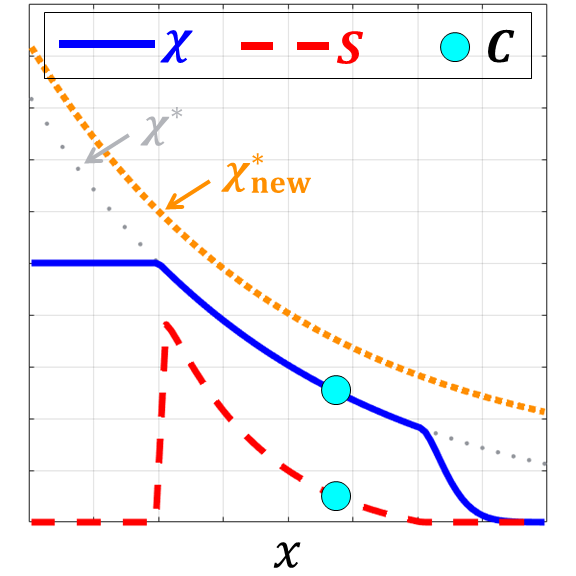}
\end{center}
\caption{Illustration of typical $\chi(x,t),S(x,t))$ at some $t$ in hydrate reservoir (top) $x\in \Omega$. For this illustration we choose $R=2\chi^*$. Top: plot of $\chi,S)$ in equilibrium, with $\chi(x,t)\leq \chi^*(x), S(x,t)\geq 0$ and $(\chi(x,t),S(x,t)) \in E_*(x)$ as in \eqref{eq:phase}. Middle: illustration of the graph $(\chi,S) \in E_*$ at the point A (equilibrium), and at the points B and C out of equilibrium (not on $E_*$) but within the physically meaningful region $(\chi,S) \in {D^0}=[0,R) \times [0,1]$.  The
contours $u(\chi,S)=u^0$ of \eqref{eq:curve} in $D^0$ for $\frac{u^0}{\chi^*}$ equal $0.4, 0.6,.., 1.8$, with the curve $u(\chi,S) = u^0=\chi^*$ separating the ``saturated'' region $\Dp$  shaded in gray from the ``un-saturated'' region $\Dm$ which is in white. The parts $\Ep$ and $\Em$ of $E_*$ are in green and magenta. Bottom: an example of $\chi^*(x),\chi(x),S(x)$ in a reservoir in out of equilibrium conditions when (B) $\chi_{new}^*<\chi^*$, and (C) when $\chi_{new}^*>\chi^*$. 
\label{fig:Domain_IC}}
\end{figure}
Consider an isolated system, with $P,T$ fixed and a fixed amount $u(\chi,S)=u^0$ of methane, and ignore any transport contributions or sources to focus on the distribution of methane between liquid and hydrate phases. 
The values $(\chi(t),S(t))$ live on a fixed curve
\ba
\label{eq:curve}
u(\chi,S)=\chi(1-S)+RS=u^0=u(\chi^0,S^0), 
\ea
in the $(\chi,S)$ plane. With some given $\chi^*$, and a corresponding fixed multi-valued graph $E_*$,  for a given $u^0$, the equilibrium point  lies at the intersection of the curve \eqref{eq:curve} with the graph $E_*$, and can be found from \eqref{eq:invert}. {The graph $E_*$ and the curves \eqref{eq:curve} are illustrated in 
Fig.~\ref{fig:Domain_IC}, with point (A) corresponding to an equilibrium case}.

The case out of equilibrium {(points (B) and (C) in Fig.~\ref{fig:Domain_IC})} is when the pair $(\chi,S)$ on \eqref{eq:curve} is away from $E_*$. For example we can have $\chi(t)>\chi^*$ {(B)}, or $S(t)>0$ with $\chi(t) <\chi^*$ {(C)}. As $t \uparrow \infty$, the points $(\chi(t),S(t))$ evolve from some $(\chi^0,S^0)$ towards some $(\chi^{\infty},S^{\infty})$ on  $E_*$ {along the curve \eqref{eq:curve}} according to some kinetic model with exchange rate $Q$.

We postulate now some conditions on $(\chi^0,S^0)$ and $u^0$ to guarantee that the kinetics leads to physically meaningful  $(\chi^{\infty},S^{\infty})$ on $E_*$.  In particular, 
from Remark~\ref{rem:uineq} we see  that $u^{\infty}=u(\chi^{\infty},S^{\infty})$ should satisfy $0\leq u^{\infty}<R$, thus we must have $0 \leq u^0=u(\chi^0,S^0)< R$. Also, non-negativity must be imposed on $(\chi,S)$. In summary, we consider the physically meaningful region $(\chi,S) \in D^0=[0,R) \times [0,1)$. 

Next, we aim to predict whether a given $u^0=u(\chi^0,S^0)$ leads to $(\chi^{\infty},S^{\infty}) \in \Em$ or to $(\chi^{\infty},S^{\infty}) \in \Ep$. In the latter saturated case  we have $\chi^{\infty} = \chi^*$ and $S^{\infty}\in[0,1)$, and $u^0=u^{\infty} = \chi^{\infty} + (R-\chi^{\infty})S^{\infty} \geq \chi^*$ by \eqref{eq:Rassum}. In the former case we have $u^0\leq \chi^*$. 
It is thus convenient to decompose $D^0=\Dm \cup \Dp$ as follows

\bas
\Dp &=& \{(\chi,S)\in D^0: u(\chi,S)\geq \chi^*\}; \;\; 
\\
\Dm&=& \{(\chi,S)\in D^0: u(\chi,S)\leq \chi^*\}.
\eas
Fig.~\ref{fig:Domain_IC} provides illustration of these definitions, and motivates our subsequent analyses.
%

\subsubsection{Three batch kinetic models} 
{Our objective is} to construct a model which works well in all of $D^0$. We start with \eqref{eq:dbkin} dubbed (KIN1) which works in $\Dp$. We include $S$ in $k_2$ in a simpler model (KIN2) which works well also in $\Dp$ only. Finally to allow the evolution towards a possible equilibrium on $\Em$ or on $\Ep$ we combine these two possible equilibria in (KIN3) using an abstract setting with the graph $E_*$.  Each (KINj) has some rate $k_j$.  In $\Dp$ and under some assumptions all three models are equivalent to one another.  {Only (KIN3) is coupled later with the transport model.} 

\paragraph{(KIN1)} The model \eqref{eq:dbkin} from \cite{lasaga,zhang-book} splits $u(t)$ as a sum of the methane amount in the $h$ phase and of the amount in the $l$ phase, and prescribes the evolution
\ba
\label{eq:KE_model1}
\mathrm{(KIN1)  } \;\;\; \dt{} ( (1-S)\chi) =-Q;\;\;  R\dt S =Q; \;\; Q = k_1(\chi-\chi^*); \;\;
(\chi(0),S(0))=(\chi^0,S^0).
\ea
This model is very intuitive: in particular, we see that $S\uparrow$ when $\chi >\chi^*$. However, (KIN1) works well only in $\Dp$ when $u^0\geq \chi^*$, i.e., when the equilibrium point  $S^{\infty}\geq 0$. When $u^0<\chi^*$ since $k_1$ does not involve $S$, the model leads to an equilibrium outside $D^0$ with $S^{\infty}< 0$. Moreover, the corresponding numerical scheme requires solution of a nonlinear algebraic equation which must be done with some care; see Sec.~\ref{sec:app-KIN}.

\paragraph{(KIN2)} Next we aim to improve (KIN1). We split $u=\chi+S(R-\chi)=\chi+\psi$, with $S=\frac{\psi}{R-\chi}$. The variable $\psi$ interpreted as the ``amount of methane stored in the hydrate phase over the saturated amount in liquid''. Given initial data $(\chi^0,S^0)$, we calculate $\psi^0=S^0(R-\chi^0)$, and postulate the evolution 
\begin{eqnarray}
\label{eq:KE_model2}
\mathrm{(KIN2)  }\;\;\;
\dt \chi =-Q; \;\;
\ddt (\psi)  =Q; \;\;
Q = k_2(\chi-\chi^*); \;\;
(\chi(0),\psi(0))=(\chi^0,\psi^0).
\end{eqnarray}
Now (KIN2) model is linear in $\chi$ and $\psi$, and $Q$ is monotone in $\chi$: the curves $\chi+\psi=u^0$ are simply the lines in the $(\chi,\psi)$ plane. These properties simplify the implementation and analysis.
However, similarly as in (KIN1) $Q$ involves properly only the equilibria on $\Ep$, and thus (KIN2) works well only in $\Dp$.

\paragraph{(KIN3)}
We modify (KIN2) so that when $u^0<\chi^*$, $Q$ leads to  some equilibrium on $\Em$, but when $u^0\geq \chi^*$, the model works identically to (KIN2) and leads correctly to some equilibrium on $\Ep$. An elegant way to do it is  to replace $\chi^*$ in the definition of $Q$ in (KIN2) by a selection $w \in w_*(\psi)$ which defaults to $\chi^*$ on $\Ep$. Here $w_*=e_*^{-1}$, and $e_*=r_*E_*$ is a rescaled version of $E_*$, with a fixed $r_*=R-\chi^*$. When $\psi>0$, we have $w=\chi^*$, but when $\psi=0$, $w\in[0,\chi^*]$. Also, $S \in E_*(\chi)$ is equivalent to  $\psi \in e_*(\chi)$ and $\chi \in w_*(\psi)$. 
The (KIN3) model we {implement and analyze} reads
\ba
\label{eq:KE_model3}
\mathrm{(KIN3)} \;\;\; \dt \chi =-Q;\;\;  \ddt ( \psi) =Q; \;\; Q = 
k_3(\chi-w);\;\; w \in w_*(\psi); \;\; 
(\chi(0),\psi(0))=(\chi^0,\psi^0).
\ea
The solution $(\chi,\psi)$ and the selection $w$ are unique. The exchange term $Q$ is monotone in $\chi$ while $-Q$ is monotone in $\psi$ which make the analysis and implementation easy. As in (KIN2), at any point of time one can calculate $S$ from $\psi$ and $\chi$.  

We provide details on (KIN1), (KIN2) and (KIN3) in Sec.~\ref{sec:app-KIN}. These inform our analysis of methane transport coupled to (KIN3). 

\section{Approximation schemes for methane transport model under equilibrium or kinetic closure}
\label{sec:LAG}
We summarize now the methane transport model in a form amenable to discretization and analyses. First we outline how the thermodynamic conditions on $(P,T)$ and $\chi^*$ are handled.
\bsub
\label{eq:LAG}
\ba
\label{eq:LAGt}
\text{Assume known }T(x,t), \text{ or solve an appropriate energy equation}
\\
\nonumber \text{   under some initial and boundary conditions}.
\\
\label{eq:LAGp}
\text{Assume known }P(x,t), q(x,t),\text{ or find these from \eqref{eq:pressure} under some boundary conditions}.
\\
\label{eq:LAGchi}
\text{Calculate } \chi^*(x,t)=\chi^*(P(x,t),T(x,t),{\chi}_{lS}(x,t)).
\\
\nonumber
\text{For kinetic model, parametrize }w_*(x,t) \text{ with } \chi^*(x,t).
\ea
\esub
For \eqref{eq:LAGchi} we use the approach described in \cite{PHTK} based on estimates of $\chi^*$ generated by CSMGem, semi-empirical model from \cite{Tishchenko2005}, and the parametric model from \cite{DZB2004} using algebraic curve fitting model for equilibrium pressure, $P_{eq}$, given in \cite{Maekawa1995}. In examples for this paper we assume $\chi_{lS} = \chi_{lS}^{sw}$.

Next use mass conservation \eqref{eq:massMC} which we divide by $\rho_l$ upon (A3). 
The equilibrium model, with 
$u(x,t)$ given by \eqref{eq:es} is
\bsub
\label{eq:EQ}
\ba
\label{eq:mass}
\partial_t\left(\phi u\right)+\nabla \cdot(q \chi) -\nabla \cdot ({d_m} \nabla \chi)&=&\frac{F_M}{\rho_l}, \; x \in \Omega, t>0,
\\
\label{eq:esxt}
\chi(x,t) &=& \min\{\chi^*(x,t),u(x,t)\}, 
\\
u(x,0)&=& u_{init}(x),\; x \in \Omega.
\\
\text{Assume boundary conditions for } \chi(x,t),&& \; x \in \partial \Omega, t>0. 
\ea
\esub

The kinetic model rewrites \eqref{eq:mass} 
in terms of $\chi$ and $\psi$. To achieve a convenient symmetrized form, we replace $ \partial_t(\phi u)=\phi \partial_t (\chi+\psi)=\phi  \partial_t (\chi) +\phi Q$ with 
$Q=-k_3(w-\chi)$  given as in \eqref{eq:KE_model3}. The model completed with appropriate initial data for $\chi$ and $\psi$ and boundary data for $\chi$ reads 
\bsub
\label{eq:kinetic}
\ba
\label{eq:transportk}
\partial_t(\phi \chi) -\phi k_3(w -  \chi) + \nabla \cdot (q \chi) - {d_m} \nabla^2 \chi&=& \frac{F_M}{\rho_l}, \;\; x \in \Omega, t>0,
\\
\label{eq:phasekw}
\partial_t (\psi(x,t)) +  k_3(w-\chi(x,t)) &= &0, \;\; w \in w_*(x,t;\psi). 
\\
\chi(x,0)=\chi_{init}(x),&&\;\; \psi(x,0)=\psi_{init}(x),\; \; x \in \Omega.
\\
\text{Assume boundary conditions for } \chi(x,t),&& \; x \in \partial \Omega, t>0. 
\ea
\esub 

{
We comment now on the couplings.  The models \eqref{eq:EQ} and \eqref{eq:kinetic} are strongly coupled to the thermodynamic conditions given in \eqref{eq:LAGchi} and to the flux $q$ found by \eqref{eq:LAGp}. For the time scales of interest in this paper, most significant are the parametrizations of \eqref{eq:esxt} and \eqref{eq:phasekw} by the quantity $\chi^*=\chi^*(x,t)$ found in \eqref{eq:LAGchi}.  On the other hand,  $\chi^*(x,t)$ depends primarily on the temperature and much less on $P(x,t)$.  At the same time, the conductivities in the energy equation are less sensitive to $S$ than the quantities in the pressure equation; see, e.g., data in \cite{LF2007}. }

{In turn, the solution to \eqref{eq:LAG}   depends on the solution to the methane transport \eqref{eq:EQ} or \eqref{eq:kinetic}. In particular, as  \eqref{eq:spress} indicates, the
local variations in $q(x,t)$ are due to 
$
\nabla \cdot q \approx \partial_t (\phi[S \frac{\Delta \rho}{\rho_l}])$,
which require re-computing $q$. In addition, the permeability $K$ in \eqref{eq:LAGp} depends on $S$, and the resulting local pressure variation may affect $\chi^*$ by the appearance of micro-cracks; see, e.g., \cite{DD11,NRuppel03}.}

{These inter--dependencies can be resolved by iteration, time-lagging, or variable freezing; we discuss these next.}

\subsection{Approximation schemes and resolving coupled components}
\label{sec:numerical}
The choice of time-stepping and spatial discretization depends on the objectives of simulation and on the competing demands of modeling accuracy, and efficiency and robustness of the solver. In this paper we are interested in modeling hydrate evolution in natural environment. The simulation scenarios we consider may involve response to changing boundary conditions for temperature or pressure such as due to the warming sea waters or sudden change in the sediment depth. We consider that these inputs vary in time on the scale of years or kiloyears but not as strongly as in production scenarios on the scale of days or hours. This assumption on the time scale motivate the choice of time stepping. 


%
The simplest way to resolve the couplings is to consider the variables $P(x)$ and $T(x)$ as time-independent over the simulation time scale, i.e., ``freeze them'' over $[0,T]$, and to solve the equilibrium model \eqref{eq:EQ}. This strategy is adopted in many hydrate models at basin scale where the pressure and temperature assumed known between any large geologic events and where $P(x)$, $T(x)$ follow closely the hydraulic gradient and geothermal gradient, respectively; see, e.g., \cite{PHTK,PMHT,TT04}. For simulation over shorter time scales this approach may require recomputing $P(x)$ and $T(x)$ over shorter time frames whenever the external controls change. With the equilibrium model \eqref{eq:EQ}, the system is immediately brought to equilibrium in the first transport step. To simulate a gradual return to equilibrium, we must use the kinetic model \eqref{eq:kinetic}. 

The most complex and comprehensive way to resolve the couplings is to use fully implicit coupling for the equilibrium model \cite{LF2007} and for the kinetic model \cite{GHW2015,GWH2016}; see also general subsurface simulators described in   \cite{IGHCCS2}. However, a fully implicit solution for several independent variables including phase behavior requires delicate time-stepping with advanced strategies to ensure global convergence and robustness of the Newton solver. 

As an intermediate strategy between the most complex and most simple, the coupling including the evaluation of thermodynamic conditions can be handled in a sequential manner or by time-lagging. {Our analysis applies in this setting.}  Here \eqref{eq:LAG} is solved at (almost) every time step. Once $(P,T,\chi^*)$ are known, \eqref{eq:EQ} or \eqref{eq:kinetic} follow. 
This is similar to a strategy common in reservoir simulation and compositional models called IMPES or IMPEC in which the pressure equation and thermodynamics conditions in \eqref{eq:LAG} and the concentration equations \eqref{eq:EQ} are solved at separate time schedules  with large pressure time steps $\Delta T$, and small transport steps $\tau=\Delta T/K$. See, e.g.,  \cite{Tchelepi2017,PLW99,WP02} where $K>1$ was used. {The sequential and time-lagging strategies carry some modeling error compared to the fully implicit model; the error decreases when small time steps are used. Additional iterations to decrease this error can be carried out over the macro time step; see recent analysis on multi-rate schemes for coupled flow and geomechanics in \cite{Kumar2016,GWH2016}; additionally, stabilization terms may improve convergence in \cite{kou2010,Radu2015robust}. If needed, we can also set $K=1$ and $\tau=\Delta T$.}

In our computational model we follow the time lagging strategy with macro-time steps but without iteration. We tested this strategy for hydrate basin modeling in \cite{PMHT}. For simplicity below we assume uniform time stepping. 

\subsection
{Time-stepping with macro time steps and concentration time steps.}
\label{sec:macro}

The concentration time step $\tau=\frac{T}{N}$ for \eqref{eq:EQ} or \eqref{eq:kinetic} is chosen to satisfy some stability constraints. The macro time step $\Delta T=K\tau$ for \eqref{eq:LAG} is chosen to be small enough so that $\chi^*(x,t),q(x,t)$ respond to the model inputs for pressure and temperature.
Here $K \geq 1$. Now $MK=N$ and $T=M \Delta T=MK\tau=N\tau$.
\bas
\label{eq:ttran}
0=t^0 < t^1 <\ldots < t^n= n \tau <\ldots < t^N = T = N\tau,\;\; \mathrm{ with\ } t^n=n\tau,\;\; n=0,1, \ldots, N.
\\
\label{eq:tpress}
0=T^0 < T^1 < \ldots <T^{m} = m \Delta T < \ldots < T^{M}=T =M\Delta t,
\;\; \mathrm{ with\ } m=0,1, \ldots, M.
\eas
Note that $T^m=m\Delta T=mK\tau=t^{mK}$. 
We outline our algorithm.

\medskip
\hrule
\vspace{.2cm}
\noindent
{\bf Time stepping (macro time steps) $m=1,2,\ldots M $}.
\\
{\bf (Macro-time step $[T^{m-1},T^m]$)}: 
\\
\phantom{ST}
Assume $S\vert_{T^{m-1}}$ is known. {Recalculate hydraulic properties}.
\\
\phantom{ST}
Solve \eqref{eq:LAG} for $(T,P,q,\chi^*)\vert_{T^m}$.
\\
\phantom{ST}
Set the values $(q,\chi^*), t \in [T^{m-1},T^m]$ from 
$(q,\chi^*)\vert_{T^m}$ 
\\
\phantom{STEP}or by interpolating between these and 
$(q,\chi^*)\vert_{T^{m-1}}$.
\\
\phantom{ST}{\bf (Concentration time steps $n=(m-1)K+1\ldots mK$)}: 
\\
\phantom{STEP}
Assume $(q,\chi^*), t \in
[T^{m-1},T^m]=[t^{(m-1)K},t^{mK}]$ known.
\\
\phantom{STEP}
Solve in each 
$[t^{n-1},t^{n}]$ the concentration problem \eqref{eq:EQ} or \eqref{eq:kinetic} for $(S,\chi)\vert_{t^{n}}$. 
\\
\phantom{ST}
With $n=mK$, set $S \vert_{T^{m}}=S \vert_{t^{n}}$. 
Advance to the next macro-time step with $m:=m+1$. 
\medskip
\hrule
\vspace{.2cm}

We devote Ex.~\ref{ex:warming} in Sec.~\ref{sec:newex} to the study of sensitivity of simulations to $\Delta T$.

\subsection
{Spatial discretization.}
We set up hexahedral grid over $\Omega$ and use finite volume type approximations. Our schemes are first-order in time, with explicit in time upwind treatment of advection, and implicit treatment of phase behavior {and diffusion} at every time step.  For simplicity we define the schemes for 1d case
with $x \in \Omega=(0,{D^{max}})$  with $x=0$ at or above BHSZ, and $x$ pointing upwards, with the flux upwards  {$q(x,t)>0$}. We cover $\Omega$ with uniform size grid cells $[x_{j-1/2},x_{j+1/2}]$, each with center at $x_j=(j+1/2)h$ where $h = x_{j+1/2}-x_{j-1/2}$. 
We denote the grid values $V_j \approx v(x_j)$, and $V_j^n \approx v(x_j,{t^n})$. The Darcy flux $q$ are defined at the cell edges $q^n_{j \pm 1/2}$. The fluxes $\chi q$ are approximated as is done for the space-dependent flux in the ``color equation'' \cite{RedLeveque}[{Chapter 9}].

We skip the presentation of schemes for \eqref{eq:LAG} which are standard; see, e.g., \cite{LF2007,PJW02}. 
However, our treatment of phase equilibria and of kinetics requires care. 
In Sec.~\ref{sec:EQ} we define numerical schemes for the concentration steps \eqref{eq:EQ} and in Sec.~\ref{sec:KIN} for \eqref{eq:kinetic}. These share the mass conservation equation discretized as follows.
We approximate $U_j^n  \approx u(x_j,t^n)$ and $\Meta_j^n \approx \chi(x_j,t^n)$  discretizing the mass conservation \eqref{eq:mass} and \eqref{eq:transportk} parts of \eqref{eq:EQ} and \eqref{eq:kinetic} by
\ba
\label{eq:upwindqd}  
\phi_j(U_j^{n} - U_j^{n-1})
+\frac{\tau}{h}(q_{j-1/2}^{n-1}\Meta_j^{n-1}-q_{j-3/2}^{n-1} \Meta_{j-1}^{n-1}) 
+ \frac{d_m\tau}{h^2}\left[2X_j^n  - X_{j-1}^n -X_{j+1}^n\right]
=\frac{\tau F_M(x_j,t^n)}{\rho_l}.
\ea
We approximate the initial data $
    U_j^0 = \frac{1}{h}\int_{x_{j-1/2}}^{x_{j+1/2}} u_{init}(x)\diff x,
$. The initial data for $(\chi,\psi)$ in the approximation to \eqref{eq:kinetic} is defined analogously. 

The equation \eqref{eq:upwindqd} is complemented with the discrete version of \eqref{eq:esxt} for the equilibrium model or with discrete version of \eqref{eq:phasekw} for the kinetic model, and with appropriate statement on the boundary conditions. These are stated in Sec.~\ref{sec:EQ} and \ref{sec:KIN} along with the analysis of their stability.

\section{Stability analysis for equilibrium model} 
\label{sec:EQ}

{We first recall notation.} 
For some grid function $V=(V_j)_j$ with $V_j \approx v(x_j)$, and $V^n=(V_j^n)_j$, we let $V^{\Delta}$ represent the collection of all $(V^n)^n$.
We recall 
$\|V\|_1 = h\sum_{j} \abs{V_j}$, and the total variation $TV(V^n)$ and total variation in time $TV_T(V^{\Delta})$ defined as 
\bas
TV(V^n) = \sum_{j\in \mathbb{Z}}\abs{V^n_j - V^n_{j-1}}, \;\; 
TV_T(V^{\Delta}) = \sum_{n=0}^{T/\tau} \left[ \tau TV(V^n) + \|V^n-V^{n-1}\|_1\right].
\eas

For the kinetic problem we work with
$\gnorm{(\Meta^n,\Psi^n)} = \gnorm{\Meta^n}+\gnorm{\Psi^n}$, and 
$TV(\Meta^n,\Psi^n)$ and
$TV_T(\Meta^{\Delta},\Psi^{\Delta})$ extended similarly to product space. 
In our analysis we study  $TV_T(U^{\Delta})$ for \eqref{eq:EQ} and  $TV_T(\Meta^{\Delta},\Psi^{\Delta})$ for \eqref{eq:kinetic}. These quantities help to predict the variability and challenges to the numerical solution depending on the data. We show that $TV_T(U^{\Delta})$ and $TV_T(\Meta^{\Delta},\Psi^{\Delta})$  increase in time depending on the 
variability and smoothness of initial data  $\chi^*(x,t)$ and $q(x,t)$ in $x$ and $t$. Stability along with consistency of the discrete schemes lead to the convergence of numerical schemes. 
Our analysis is also useful to understand the sensitivities of the models \eqref{eq:EQ} and \eqref{eq:kinetic} on their data. 

\paragraph{{(AA)}Assumptions for analysis.} 
We analyze only the scheme \eqref{eq:upwindqd}  for transport model complemented by an equilibrium or kinetic closure to be stated, under assumptions (A1-A5). We assume that the data $\chi^*(x,t)$ and $q(x,t)$ found by \eqref{eq:LAG} is known over each macro-time step $[T^{m-1},T^m]$ and varies in some predictable fashion. {As usual, to study the accumulation of the discretization error in time, we set $F_M=0$.  We also set set $d_m=0$ to focus on the advection dominated case}. We consider the transport problem on $x\in \R$ (that is, $j\in \mathbb{Z}$)  rather than $x \in \Omega$, which avoids dealing with a mixture of boundary and initial conditions in the analysis. For this we assume that initial data and the solution to the transport problem have compact support in some $\Omega_S \subset \R$ with measure $\omega_S$, and this reduces summing over $j\in \mathbb{Z}$ to $j \in \mathbb{Z}^0$.
Clearly, realistic simulations consider a bounded domain and boundary conditions.

{Finally, we assume the sediment is homogeneous with $\phi(x)=\phi_0=const$, and we drop $\phi_0$ while keeping the notation unchanged, but the analysis could be amended easily as long as $\phi(x)$ is smooth and bounded away from $0$. In particular, in \eqref{eq:EQ} we could change variables and set $\overline{u}=u(x,t)\phi(x)$ with $\chi(x,t) = \min\{\chi^*(x,t),\frac{\overline{u}(x,t)}{\phi(x)}\}$.}

\begin{figure}
\centering
\includegraphics[width=.4\textwidth]{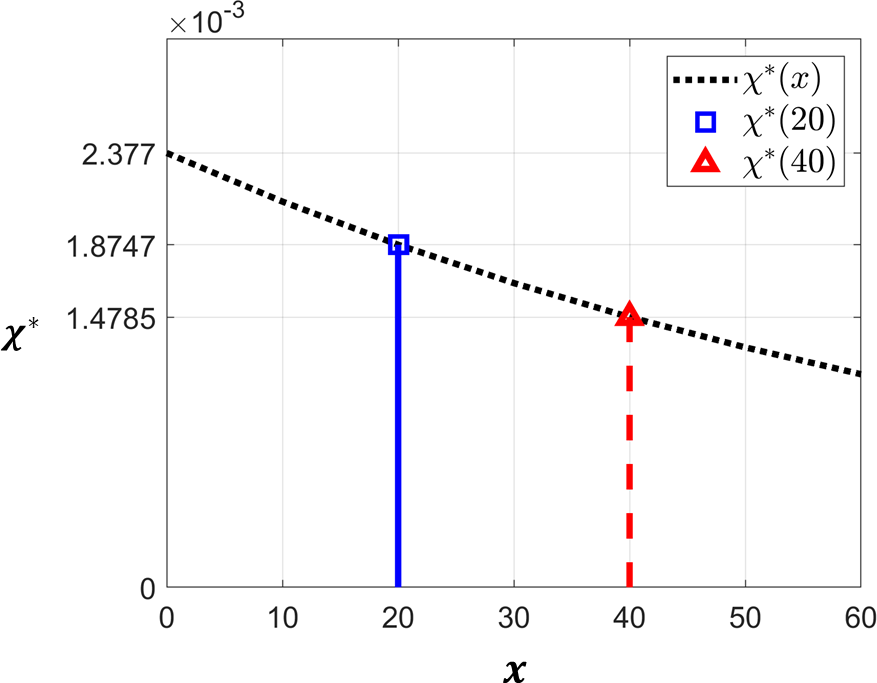}
\includegraphics[width=.4\textwidth]{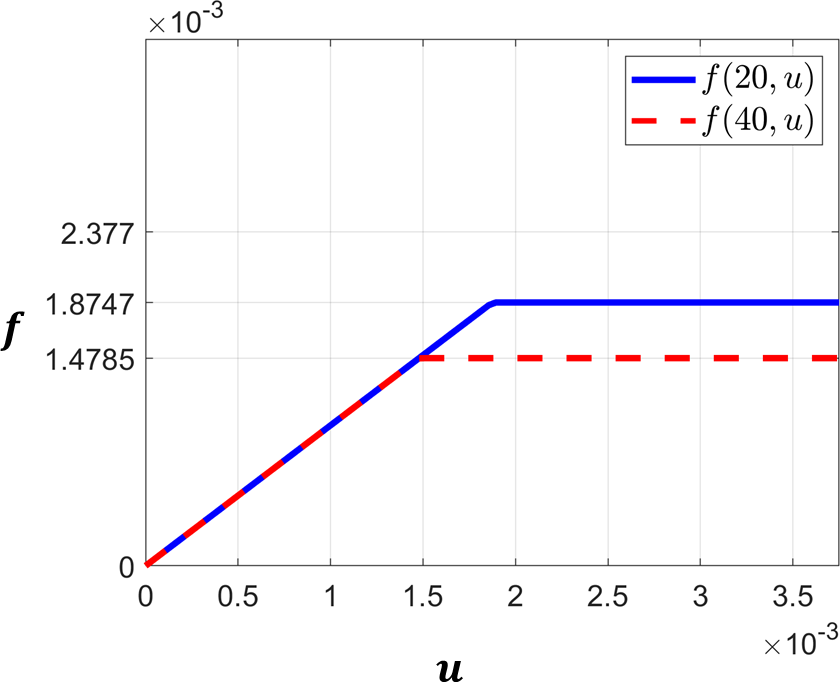}
\caption{\label{fig:chiu}
Illustration of $\chi^*(x)$, and of the flux function $f(x;u)$ with data from Ulleung Basin case UBGH2-7 \cite{PHTK}. Left: typical $\chi^*(x)$ in homogeneous sediment, with values $\chi^*(x)$ highlighted at $x=20$, and $x=40$.  Right: the flux function $f(x;u)$ for $x=20$, and $x=40$. Note that the flux function $f(x,t;u)$ is piecewise linear in $u$ and features a corner at $u=\chi^*(x^*)$.  }
\end{figure}

We now restate \eqref{eq:upwindqd} under the assumptions (AA) amended by the discrete version of \eqref{eq:esxt}
\ba
\label{eq:upwindEQ}  
\frac{1}{\tau}(U_j^{n} - U_j^{n-1})
+\frac{\tau}{h}(q_{j-1/2}^{n-1}\Meta_j^{n-1}-q_{j-3/2}^{n-1} \Meta_{j-1}^{n-1}) 
 =0,\;\;
\Meta_j^n = \min\{\chi^*(x_j,{t^n}),U_j^n\}, \;\; j \in \mathbb{Z}^0.
\ea
We analyze {this scheme recognizing its familiar upwind character}
\ba
\label{eq:upwind}
    U_j^{n} = U_j^{n-1} - \frac{\tau}{h} \left[F_j^{n-1} - F_{j-1}^{n-1}\right],
\;\;
F_j^{n-1} = q_{j-1/2}^{n-1}X_j^{n-1}=f(x_j,t^{n-1};U_j^{n-1})
\ea
for a conservation law with the flux function $f(x,t;u)$ which we set from \eqref{eq:massMC} under (AA)
\begin{subequations}
\label{eq:chi}
\ba
\label{eq:chiu}
\partial_t u+\partial_x f&=&0; \;\; \text{ for } x\in \R, \, t \in [0,T)
\\
\label{eq:chif}
f(x,t;u)&=&q(x,t) \chi(x,t)= q(x,t) \min\{\chi^*(x,t),u(x,t)\},
\\
u(x,0)&=&u_{init}(x).
\ea
\end{subequations}
The function $f$ is illustrated in Fig.~\ref{fig:chiu}
for a typical homogeneous unconsolidated sand reservoir with data from Ulleung Basin \cite{PHTK}, with $x$ pointing upwards and $q(x,t)\approx const>0$.

{For stability of \eqref{eq:upwind} when $f=f(u)$  proving a bound for $TV_T(U^{\Delta})$ is usually the first step in the analysis of convergence \cite{GreenLeveque,RedLeveque}; the first-order  upwind scheme \eqref{eq:upwind} converges at the rate of $O(\sqrt{h})$ in $L^1$ \cite{Kruzhkov1960,Kuznetsov1976,TT1995}, {the best possible} convergence rate in the presence of discontinuities \cite{sabac1997}.}

{However, when $f=f(x;u)$, the solutions to $u_t+f_x=0$ do not in general obey the maximum principle, and its quasilinear form $u_t+f_u(x;u)u_x=-f_x(x;u)$ with $\abs{f_x} \leq L_f$ reveals that the solution may grow pointwise as $O(L_f t)$ while its total variation may grow in time; see the discussion in  \cite{RedLeveque}{[Chapter 9]} and \cite{KT2004,K2003,Towers2000}. Under some assumptions the analysis in \cite{LTJ1995} predicts that $TV(U^n)$ found with Godunov scheme grows linearly in time, but these assumptions are not applicable to \eqref{eq:upwindEQ}. }

We formulate therefore 
our own auxiliary stability result for \eqref{eq:upwind} {similar to those known, e.g., from \cite{RedLeveque}.  Surprisingly, we did not find it stated in literature, thus we provide detailed proof in Sec.~\ref{sec:app-EQ}.}

\begin{proposition}
\label{th:stability}
Consider \eqref{eq:chi} and assume suppose that  $f\in C^2_b(\overline{\Omega}\times\R_+ \times \R)$ 
and is nondecreasing in $u$.  Let the time step size $\tau$ be small enough so that
\ba
\label{eq:tau}
\max_{(x,t,u)} \abs{\frac{\tau}{h}f_u(x,t;u)} \leq 1.
\ea
Let also some initial data $U^0$ be given, with bounded variation, and let also $U^{\Delta}$ be the solution to 
\eqref{eq:upwind} with compact support of measure bounded by $\omega_S$.  Let also
\ba
\label{eq:L12}
L_1 = \max_{(x,t,u)}\{\abs{f_{xu}(x,t;u)},\abs{f_{xx}(x,t;u)}\}; \;\;\;
L_2 = \max_{(x,t,u)}\{\abs{f_u(x,t;u)},\abs{f_x(x,t;u)}\}.
\ea
Then we have for all $n>0$
\begin{subequations}
\ba
\label{th:TVstability_space}
    TV(U^n)&\leq &C_1(T) = TV(U^0)e^{TL_1} + 2\omega_S(e^{TL_1}-1),
\\
\label{eq:gnorm}
    \norm{U^{n+1}-U^n}{1} &\leq &\tau C_2(T), \;\; C_2(T)= L_2(C_1(T)+\omega_S),
\\
\label{th:TVstability_time}
    TV_T(U^{\Delta})&\leq& C_3(T) = T(C_1(T)+C_2(T)).
\ea
\end{subequations}
\end{proposition}

To apply Proposition~\ref{th:stability} to \eqref{eq:upwindEQ}, we consider $f(x,t;u)$ defined in \eqref{eq:chif} depending on the physical data $q,\chi^*$. We see that $f(\cdot;u)$ is continuous nondecreasing and piecewise linear in $u$ and differentiable except where $u(x,t)=\chi^*(x,t)$.   Since $f$ is at best piecewise smooth in $u$, Proposition~\ref{th:stability} applies only to some regularization of \eqref{eq:chi}
\ba
\label{eq:regularization}
u_t^{\epsilon}+f^{\epsilon}(x,t;u^{\epsilon})_x =0, \;\;
u^{\epsilon}(x,0) = u_{init}(x),
\ea
in which $f$ is approximated with some $C^2$ smooth, positive nondecreasing function $f^{\epsilon}$. 
Since we can make $f^{\epsilon}\approx f$ arbitrarily close, we trust that \eqref{eq:regularization} closely resembles \eqref{eq:chi}; we illustrate this regularization in Ex.~\ref{ex:regularization} in Sec.~\ref{sec:results}. 

{Next question is whether the assumptions on $f^{\epsilon}$ itself are reasonable for a real hydrate reservoir simulation.  First, the problem \eqref{eq:chi} is only a simplification of the strongly coupled dynamical problem \eqref{eq:LAG}--\eqref{eq:EQ}, and our stability analysis does not explain or refer to the strength of the couplings. Instead, we make a-priori assumptions on the data which allow to conclude stability and predict the variability of solutions.
In particular, we predict {variability of $u$ quantified by}  $TV_T(U^{\Delta})$ depending on {the constants} $L_1,L_2$ given in \eqref{eq:L12}; these are small only when $q$, $q_x$, the lithology and $\chi^*$ vary smoothly. We discuss this in detail below.}

\subsection{Assumptions required for the stability of \eqref{eq:upwindEQ} in a hydrate reservoir}
\label{sec:eqs}

From the form of \eqref{eq:chi} and properties of $f$, we expect its solution $u(x,t)$ to feature a family of right and left states travelling at different speeds due to the ``corner'' of $f$ at $u=\chi^*$. In particular, the speed of the state for any $u>\chi^*$ is zero; this leads to $S>0$, i.e., the growth of immobile amount of methane trapped as solid hydrate with the appearance of sharp bands of hydrate. We rewrite \eqref{eq:regularization} as
\begin{equation}
    \label{eq:quasilinear}
    u_t^{\epsilon}+f_u^{\epsilon}(x,t;u^{\epsilon})u_x^{\epsilon} = -f_x^{\epsilon}(x,t;u^{\epsilon}),
\end{equation}
which illustrates that the solution $u \approx u^{\epsilon}$ and the corresponding $S \approx S^{\epsilon}$ grow along its characteristics with a rate bounded  by the source $-f_x^{\epsilon}$. 
To quantify, we define
\ba 
\label{eq:l}
{L_q=\max_{(x,t)} \abs{q(x,t)};\; 
L_{q_x} = \max_{(x,t)} \abs{q_x(x,t)};} \; 
L_{\chi^*} = \max_{(x,t)} \abs{\chi_x^*(x,t)}; \; 
L_3 = \max_{(x,t)} \abs{\chi_t^*(x,t)}.
\ea

{
\begin{remark}
\label{rem:stability}
Assume that $\chi^*$ and $q$ vary mildly so that $L_1,L_2,L_{\chi^*}$ are finite and 
that
\ba
\label{eq:cfl}
\frac{\tau}{h} L_q \leq 1,
\ea
so that \eqref{eq:tau} holds.  Then the scheme \eqref{eq:upwindEQ} for the equilibrium model is weakly stable.
\end{remark}
}

We comment now on the constants $L_1,L_2,L_{\chi^*}$ in realistic reservoirs.
Assume first the quasi-static case in a homogeneous reservoir with $P(x)$ and $ T(x)$ fixed in time $t$, and with $0<q=L_q$. In this case $\norm{f_x}{\infty} =L_q L_{\chi^*}$. 
Consider for example $\chi^*(x) = a\exp(-bx)$ from \cite{PHTK} given with some $a>0$ and small $b>0$. Now $L_{\chi^*}=ab$, $f_x<0$, $\norm{f_x}{\infty}=L_q ab$ and $L_1=\norm{f_{xx}}{\infty} =L_q ab^2$ is small. In turn, we can check that $L_2 =L_q \max(1,ab)$. These stability constants correlate well with the predictions of
hydrate band growth in nature which are large when $q$ is large.

{Consider next heterogeneous reservoirs.} Here the maximum solubility $\chi^*(x)$ depends on the type of sediment, e.g., in grain size \cite{DD11}. Consequently, close to some interfaces between different sediment layers, hydrate can accumulate much faster than elsewhere 
\cite{DD11,Rempel2011,VR2018}. The locally high hydrate accumulation can be predicted from \eqref{eq:quasilinear}, since at a discontinuity of $\chi^*$, its weak derivative $\partial_x \chi^*$ is a Dirac term which may cause a dramatic local increase of $U_j^n$ and of the saturation. We illustrate this later in Ex.~\ref{ex:hetero} in Sec.~\ref{sec:results}. 

Finally, we consider the time dependent case {closest to the strongly coupled hydrate systems} when $q=q(x,t)$ and $\chi^*=\chi^*(x,t)$ and when $T=T(x,t)$ and $P=P(x,t)$. Now the magnitude of $f_x$ comes from both $q_x \chi^*$ and $q\chi^*_x$ which may have opposite signs and disparate magnitudes depending, e.g., on the solutions of \eqref{eq:spress}. It is hard to predict these a-priori, and we can only make assumptions that the constants in Remark~\ref{rem:stability} are bounded. Simulation with $\chi^*=\chi^*(x,t)$ which varies in time is considered in Ex.~\ref{ex:warming} in Sec.~\ref{sec:results}.

{This discussion completes our analysis of the equilibrium case. Based on Proposition 1, we expect the rate of convergence $O(\sqrt{h})$ for the solutions to \eqref{eq:upwindEQ}; this is  confirmed by numerical experiments in Sec.~\ref{sec:results}.}

\section{Scheme for kinetic model and its stability}
\label{sec:KIN}

Now we consider {a numerical scheme for} the kinetic model \eqref{eq:kinetic}. {We implement the general case with source terms and diffusion} and approximate $\chi(x_j,{t^n}) \approx \Meta_j^n$ and $\psi(x_j,{t^n}) \approx \Psi_j^n$, with $U_j^n=\Meta_j^n+\Psi_j^n$. Given $(\Meta_j^{n-1},\Psi_j^{n-1})$ we find $(\Meta_j^{n},\Psi_j^{n},W_j^n)$ as solutions to the local nonlinear system at every $j$; in this local problem the kinetic terms $Q_j^n$ are handled implicitly. We set $k=\tau k_3$. 
{The scheme \eqref{eq:upwindqd}  
under assumptions (AA) in the form directly amenable to analysis reads} 
\bsub
\label{eq:schemek0}
\ba
(\Meta_j^n-\Meta_j^{n-1}) -k(W_j^n-\Meta_j^n)
+\frac{\tau}{h} ({q_{j-1/2}^{n-1}}\Meta_j^{n-1}-{q_{j-3/2}^{n-1}} \Meta_{j-1}^{n-1})  &=&0,
\\
(\Psi_j^n-\Psi_j^{n-1})+k(W_j^n-\Meta_j^n)&=&0,\;\;  W_j^n \in w_*(x_j;\Psi_j^n).
\ea
\esub
In practice we solve \eqref{eq:schemek0} as follows, denoting $\wt{k}=\frac{k}{1+k}$. Given previous time step values $(\Meta_j^{n-1},\Psi_j^{n-1})_j$, at every $j$ we solve for $(\Meta_j^n,\Psi_j^n)$ the local nonlinear system 
\bsub
\label{eq:schemek}
\ba
\label{eq:1}
\Meta_j^n -k(W_j^n-\Meta_j^n) &=&  F_j^n, \;\; F_j^n=\Meta_j^{n-1}(1-\frac{\tau}{h}q_{j-1/2}^{n-1})+ \frac{\tau}{h}q_{j-3/2}^{n-1} \Meta_{j-1}^{n-1},
\\
\label{eq:2}
\Psi_j^n+k(W_j^n-\Meta_j^n)&=&G_j^n,\;\; G_j^n = \Psi_j^{n-1},\;\; 
W_j^n=w_*(x_j;\Psi_j^n).
\ea
\esub
This is a $2 \times 2 $ nonlinear stationary system of equations with a maximal monotone graph $w_*$. It is uniquely solvable with the following explicit formulas which follow from Sec.~\ref{sec:monotone}. Since \eqref{eq:1} is linear in $\Meta_j^n$, we can formally calculate $
\Meta_j^n=\frac{1}{1+k}(F_j^n+kW_j^n)$. After we plug this to \eqref{eq:2} 
we get $\Psi_j^n + \wt{k}w_*(\Psi_j^n)=G_j^n + \wt{k}F_j^n$. Applying the resolvent $\myre=\myre^{w_*}_{\wt{k}}=(I + \wt{k} w_*)^{-1}$ of $w_*$ we obtain $\Psi_j^n=\myre(G_j^n+\wt{k}F_j^n)$. We substitute to get $W_j^n$ and $X_j^n$. 
Finally we can calculate the saturations
$S_j^n=\frac{\psi_j^n}{R-\Meta_j^n}$.

For stability of the scheme we need an auxiliary result formulated for \eqref{eq:schemek} with indices dropped and with inputs $F,G$ and outputs $\Meta,\Psi$. 

\begin{lemma}\label{lem:comparison}
\mps Consider \eqref{eq:schemek} with the right hand side $(F,G)$ and  solutions  $(\Meta,\Psi,W)$. Consider also the right hand side $(\overline{F},\overline{G})$ with the corresponding solutions
$(\overline{\Meta},\overline{\Psi},\overline{W})$ to \eqref{eq:schemek}. 
The following comparison principle and stability hold
\bsub
\ba
\label{eq:comparison}
\abs{\Meta-\overline{\Meta}}+\abs{\Psi-\overline{\Psi}} &\leq& {\abs{F-\overline{F}}+\abs{G-\overline{G}}}.
\\
\label{eq:stability}
\abs{\Meta}+\abs{\Psi} &\leq& \abs{F}+\abs{G}.
\ea
\esub
We also have
$
W-\overline{W}=\frac{1}{k}(G-\overline{G})+(\Psi-\overline{\Psi})(1-\frac{1}{k}).
$
\end{lemma}
\begin{proof}
The result \eqref{eq:comparison} is a special case in $\R \times \R$ of the result, in \cite{HS1995} for $a(u)+c(u-v) =f; b(v)-c(u-v)=g$, where $a(\cdot)$ is maximal monotone, $b(\cdot)$ is strongly monotone and continuous, and $c(\cdot)$ is maximal monotone single valued. The stability result \eqref{eq:stability} follows from the comparison principle \eqref{eq:comparison}. In turn, the algebraic formula for $W-\overline{W}$ follows directly from algebra. 
\end{proof}

\subsection{TV-stability for the kinetic scheme \eqref{eq:schemek}}

Now we prove properties of \eqref{eq:schemek}. Throughout we assume that the CFL condition \eqref{eq:cfl} holds and that the constants $L_3,\lchistar$ are finite.

First we apply the stability part of  Lemma~\ref{lem:comparison} directly to \eqref{eq:schemek} to obtain 
\ba
\abs{\Meta_j^n} + \abs{\Psi_j^n} \leq
\abs{\Meta_j^{n-1}}\left(1-\frac{\tau}{h}q_{j-1/2}^{n-1}\right) + \frac{\tau}{h}q_{j-3/2}^{n-1}\abs{\Meta_{j-1}^{n-1}} + \abs{\Psi_j^{n-1}}.
\ea
Multiplying by $h$ and summing both sides over $j\in\mathbb{Z}^0$, and collapsing the first two terms on the right hand side, we obtain the stability result.
We {obtain} that scheme \eqref{eq:schemek} is stable {in the product space}
\ba\label{eq:stablek}
\gnorm{(\Meta^n,\Psi^n)} \leq 
\gnorm{(\Meta^{n-1},\Psi^{n-1})}.
\ea
Next we prove weak TV-stability which reveals the dependence of $w_*=w_*(x_j;\cdot)$ on $x_j$.

\begin{proposition}
\label{prop:kinetic}
Assume $\chi^*(x)$ is smooth so that $L_3$ and $L_{\chi^*}$ given by \eqref{eq:L12} are finite.  Assume also $(X^{\Delta},\Psi^{\Delta})$ have compact support with measure bounded by $\omega_S$. If CFL condition \eqref{eq:cfl} holds, then 
\ba
\label{eq:tvstabilityk}
TV(\Meta^n,\Psi^n) \leq TV(\Meta^0,\Psi^0) +C_4T, ;\;
C_4{=2k_3\omega_S \lchistar}.
\ea
\end{proposition}
\begin{proof} \mps
We write the system \eqref{eq:schemek} at $j$ and at $j$-$1$ and at {$t^n$}. We set  $\Psi=\Psi^n_j$ and $\overline{\Psi}=\Psi^n_{j-1}$, with analogous notation for other variables, and consider 
\bsub
\label{eq:twosys}
\ba
\label{eq:twosys1}
\Meta-kW+k\Meta&=&F,
\\
\label{eq:twosys2}
\Psi+kW-k\Meta&=&G; \; W \in w_*(\Psi),
\\
\label{eq:twosys3}
\overline{\Meta}-k\overline{W}+k\overline{\Meta}&=&\overline{F},
\\
\label{eq:twosys4}
\overline{\Psi}+k\overline{W}-k\overline{\Meta}&=&\overline{G}; \; \overline{W} \in v_*(\overline{\Psi}).
\ea
\esub
Here for shorthand we denoted the graph $w_*(x_j,t^n)$ by $w_*$ and {a different} graph $w_*(x_{j-1},t^n)$ at $x_{j-1}$ by $v_*(\cdot)$. {Since the graphs $w_*$ and $v_*$} are not the same, we cannot directly apply Lemma~\ref{lem:comparison}.  {Instead}, we rewrite the third and fourth equations with $w_*$ instead of $v_*$, move the difference between $w_*$ and $v_*$ to the right hand side, and examine the difference $w_*-v_*$ due to their ``height'', respectively, $\chi^*(x_j,t^n)$ and $\chi^*(x_{j-1},t^n)$. 

For $\psi>0$ we can write
\bas 
v_*(\psi)=w_*(\psi)-A_jh, \;\; \psi>0,\;\;  A_j=\frac{d}{dx}\chi^*(\overline{x}_j,t^n),\;\; \overline{x}_j \in (x_{j-1},x_j).
\eas
When $\psi=0$ both $w_*$ and $v_*$  are set-valued, and we must work with their Yosida approximations $w_{\lambda}$ and $v_{\lambda}$. In fact  for small $\psi$ we have $w_{\lambda}(\psi)= v_{\lambda}(\psi)$, while for any $\psi$ and $\lambda$ we have $v_{\lambda}(\psi)=w_{\lambda}(\psi)-A_j(\lambda;\psi) h$, with $\abs{A_j(\lambda;\psi)} \leq \lchistar$ from \eqref{eq:l}. 

Reconsidering \eqref{eq:twosys} with $v_{\lambda}$ and $w_{\lambda}$ instead of $v_*$ and $w_*$ but keeping the notation unchanged otherwise, we calculate $k\overline{W} = kv_{\lambda}(\overline{\Psi}) = kw_{\lambda}(\overline{\Psi})-khA_j$, and the third and fourth equations read now
\bas
\overline{\Meta}-kw_{\lambda}(\overline{\Psi})
+k\overline{\Meta}&=&\wt{F}=\overline{F}-khA_j(\lambda;\overline{\Psi}),
\\
\overline{\Psi}+kw_{\lambda}(\overline{\Psi})
-k\overline{\Meta}&=&\wt{G}=\overline{G}+khA_j(\lambda,\overline{\Psi}).
\eas
We can now apply the comparison Lemma~\ref{lem:comparison} for the maximal monotone  $w_{\lambda}$ and inputs $F,G,\wt{F},\wt{G}$. We apply the uniform bound on $A_j$ in \eqref{eq:comparison}, notice $\abs{\wt{F}-F}\leq \abs{\overline{F}-F}+kh\lchistar$ and $\abs{\wt{G}-G}\leq \abs{\overline{G}-G}+kh\lchistar$. Taking the limit as $\lambda \to 0$ we obtain, reverting back to the original notation of \eqref{eq:schemek} that
\begin{multline}
\label{eq:twosysstar}
\abs{\Meta_j^n-\Meta_{j-1}^n} + \abs{\Psi_j^n-\Psi_{j-1}^n}
\leq \abs{\Meta_j^{n-1}-\Meta_{j-1}^{n-1}}\left(1-{\frac{\tau}{h}q_{j-1/2}^{n-1}}\right)+ {\abs{\frac{\tau}{h}q_{j-3/2}^{n-1}}}\abs{ \Meta_{j-1}^{n-1}-\Meta_{j-2}^{n-1}}
\\
+ \abs{\Psi_j^{n-1}-\Psi_{j-1}^{n-1}}
+ 2kh\lchistar.
\end{multline}
The term $2kh\lchistar$ will accumulate giving weak rather than strong stability. 
Summing \eqref{eq:twosysstar} over those $j \in \mathbb{Z}^0$  with {$\sum_j h \leq \omega_S$}, we collapse the first two terms on the right hand side, and with  $k=\tau k_3$ we get  
\ba
\label{eq:tvstability-sys}
TV(\Meta^n,\Psi^n)\leq TV(\Meta^{n-1},\Psi^{n-1}) +2\tau k_3 \omega_S\lchistar.
\ea
Applying recursively, we obtain \eqref{eq:tvstabilityk} with $C_4 = 2k_3\omega_S L_{\chi^*}$. 
\end{proof}

\begin{remark}
The weak TV-stability result \eqref{eq:tvstabilityk} in the product space for the kinetic problem \eqref{eq:kinetic} is similar to the weak stability 
\eqref{th:TVstability_space}
we obtained for $U^{\Delta}$ in the equilibrium model \eqref{eq:chi}, with the difference in the constants depending on  $\chi^*(x)$, and the absence of the factor $(1+\tau L_1)$ in \eqref{eq:tvstability-sys} in the product space. 
\end{remark}

\subsection{TV stability in time}
Given the known $(\Meta_j^{n-1},\Psi_j^{n-1})_j$ the next goal is to bound the terms $\Meta_j^n-\Meta_j^{n-1}$ and $\Psi_j^n-\Psi_j^{n-1}$. For this, we need a handle on $Q_j^n \propto W_j^n-\Meta_j^n$ which quantifies the discrepancy from the equilibrium. We estimate $Q_j^n$ in terms of $Q_j^{n-1}$. 

\begin{lemma}
\label{lemma:Q}
Under the assumption of Proposition~\ref{prop:kinetic} we have that 
\bas
\norm{Q^n}{1}\leq C_5(T).
\eas
\end{lemma}

To prove the lemma, we estimate the terms in a regularized version of \eqref{eq:schemek}. {Additional challenge is to allow for possible variability of $w_*$ in time}. We consider some smooth single valued approximations $w_{\lambda}$  of $w_*\vert_{x_j,{t^n}}$ and 
$v_{\lambda}$  of $w_*\vert_{x_j,{t^{n-1}}}$. The difference between these $w_{\lambda}(\psi)-v_{\lambda}(\psi)=B_j^n(\psi) \tau$ can be estimated uniformly in $\psi$ with $\abs{B_j^n} \leq L_3$, where $L_3$ is given in \eqref{eq:l}. 

\begin{proof}
We rearrange \eqref{eq:schemek}, drop $j$, and seek the solution $(\Meta_{\lambda}^{n},\Psi_{\lambda}^{n})$ to the regularized problem
\bsub
\label{eq:kineticbn}
\ba
\label{eq:first}
\Meta_{\lambda}^n - \Meta^{n-1}-kQ_{\lambda}^n    &=& F^{n-1},
\\
\label{eq:second}
\Psi_{\lambda}^n - \Psi^{n-1} + kQ_{\lambda}^n &=& 0,
\ea
\esub
where $Q_{\lambda}^n=w_{\lambda}(\Psi_{\lambda}^n)-\Meta_{\lambda}^n$ and $F^{n-1} = -\frac{\tau}{h}q_{j-1/2}^{n-1}\Meta_j^{n-1} +\frac{\tau}{h}q_{j-3/2}^{n-1}\Meta_{j-1}^{n-1}$. To get the estimates for $Q_{\lambda}^n$ in terms of  $Q_{\lambda}^{n-1}=v_{\lambda}(\Psi^{n-1})-\Meta^{n-1}$, we break the expression 
\bas
w_{\lambda}(\Psi_{\lambda}^n)-v_{\lambda}(\Psi^{n-1})=
w_{\lambda}(\Psi_{\lambda}^n)-w_{\lambda}(\Psi^{n-1})
+w_{\lambda}(\Psi^{n-1})-v_{\lambda}(\Psi^{n-1})=
b(\Psi_{\lambda}^n - \Psi^{n-1}) + B\tau.
\eas
Here $b = w'_{\lambda}(\wt{\Psi_{\lambda}})\geq 0$ with some $\wt{\Psi_{\lambda}}$, and $B=B_j^n$ discussed above, with $\abs{B}\leq L_3$.  Now we multiply \eqref{eq:second} by $b$ and subtract \eqref{eq:first} from \eqref{eq:second}. Rearranging we obtain
\bas
Q_{\lambda}^n(1+k(1+b)) = Q^{n-1}-F^{n-1}+ B\tau.
\eas
We take absolute value, note $b\geq 0$, and pass to the limit with $\lambda$, to obtain, bringing back the index $j$
\bas
(1+k)\abs{Q_j^n} & \leq &\abs{Q_j^{n-1}}+\abs{F_j^{n-1}}+\tau \abs{B_j^n}.
\eas
Here, with $L_q$ and $L_{q_x}$ defined in \eqref{eq:l}, $\abs{F_j^{n-1}}$ is bounded above 
\bas
\abs{F_j^{n-1}} &=& \frac{\tau}{h}\abs{q_{j-1/2}^{n-1}\Meta_j^{n-1} - q_{j-1/2}^{n-1}\Meta_{j-1}^{n-1} + q_{j-1/2}^{n-1} \Meta_{j-1}^{n-1} - q_{j-3/2}^{n-1} \Meta_{j-1}^{n-1}},\\
&\leq &
\frac{\tau}{h}L_q \abs{\Meta_j^{n-1} -\Meta_{j-1}^{n-1}}+ \tau L_{q_x}\abs{\Meta_{j-1}^{n-1}},\\
&\leq&
\frac{\tau}{h}L_q\left(\abs{\Meta_j^{n-1} -\Meta_{j-1}^{n-1}} + \abs{\Psi_j^{n-1} -\Psi_{j-1}^{n-1}}\right) + \tau L_{q_x}\left(\abs{\Meta_{j-1}^{n-1}} + \abs{\Psi_{j-1}^{n-1}}\right).
\eas
Multiply both sides by $h$ and sum over $j\in \mathbb{Z}^0$ with $\sum_{j\in \mathbb{Z}^0} h\leq \omega_S$. Then apply \eqref{eq:stablek} and \eqref{eq:tvstabilityk} to get
\bas
(1+k)\norm{Q^n}{1} &\leq &
\norm{Q^{n-1}}{1} + \tau \left[L_q TV(X^{n-1},\Psi^{n-1})+L_{q_x}\gnorm{(X^{n-1},\Psi^{n-1})} + L_3\omega_S\right],\\
&\leq & \dots \leq
\norm{Q^0}{1} + T\left[C_6(T) + L_3\omega_S\right],
\eas
where $C_6(T) = L_q (TV(X^{0},\Psi^{0})+ C_4 T) + L_{q_x}\gnorm{(X^{0},\Psi^0)}$.
With $k>0$, we have completed the proof.
\end{proof}

\begin{lemma} \label{lem:gridnorm} \mps
Suppose \eqref{eq:tvstabilityk} holds for $(\Meta,\Psi)$. Then we have
\bas
\|(\Meta^n-\Meta^{n-1},\Psi^n-\Psi^{n-1})\|_{\Delta,1} \leq \tau C_7(T).
\eas
\end{lemma}
\begin{proof}
Rewrite \eqref{eq:schemek} estimating in the form
\bas
\abs{\Meta_j^n -\Meta_j^{n-1}} &\leq& k\abs{Q_j^n}+
\frac{\tau}{h}L_q\abs{\Meta_j^{n-1} - \Meta_{j-1}^{n-1}} + \tau L_{q_x}\abs{\Meta_{j-1}^{n-1}}, \\
\abs{\Psi_j^n - \Psi_j^{n-1}} &\leq& k\abs{Q_j^n}.
\eas
Next we multiply by $h$, take the sum over $j\in \mathbb{Z}^0$ and add these to get
\bas
\gnorm{(X^n-X^{n-1},\Psi^n-\Psi^{n-1})} &\leq & \tau \left[2k_3\norm{Q^n}{1} + L_q TV(X^{n-1},\Psi^{n-1}) + L_{q_x} \gnorm{(X^{n-1},\Psi^{n-1})}\right],\\
&\leq &
\tau \left[2k_3\norm{Q^n}{1}+C_6(T)\right].
\eas
With \eqref{eq:stablek}, \eqref{eq:tvstabilityk} and the estimates for $Q$ from Lemma~\ref{lemma:Q}, we get
\bas
\|(\Meta^n-\Meta^{n-1},\Psi^n-\Psi^{n-1})\|_{\Delta,1} &\leq & \tau C_7(T) = 
\tau \left[2k_3\norm{Q^0}{1} + (1+2k_3T)C_6(T) + 2k_3TL_3\omega_S\right].
\eas
\end{proof}

Combining Lemma~\ref{lem:gridnorm} and \eqref{eq:tvstabilityk}, we 
 conclude with this main result. 
\begin{proposition}\label{prop:TVstabilitytime}
Under hypotheses of Proposition~\ref{prop:kinetic} we have
\ba
\label{eq:tvtKIN}
TV_T(\Meta^n,\Psi^n) \leq C_8(T)=T\left[C_4T + C_7(T)\right].
\ea
{Here $C_4 = 2k_3\omega_S L_{\chi^*}$, $C_7(T) = 
\left[2k_3\norm{Q^0}{1} + (1+2k_3T)C_6(T) + 2k_3TL_3\omega_S\right]$, and  
$C_6(T) = L_q (TV(X^{0},\Psi^{0})+ C_4 T) + L_{q_x}\gnorm{(X^{0},\Psi^0)}$.}
\end{proposition}

{As in equilibrium case discussed in Sec.~\ref{sec:eqs}, this stability result depends on the variability of $q$ and $\chi^*$ and on the initial discrepancy from the equilibrium through the constants in \eqref{eq:tvtKIN}.} 

\section{Numerical examples}
\label{sec:results}
In this section we provide examples for equilibrium and kinetic models. Our goal is {to confirm the theory and in particular demonstrate convergence of the schemes for reasonably realistic cases as well as}  to demonstrate the practical limitations. We set 
$
\nu = L_q \frac{\tau}{h}<1
$, and consider only 1d simulations.

For the equilibrium model we compare the numerical solution obtained by our scheme \eqref{eq:upwindEQ} with an analytical solution and we study effects of regularization; we also confirm the rate of convergence of $O(\sqrt{h})$. {We study also the kinetic model and scheme \eqref{eq:kinetic}, illustrate its convergence and show} the dependence of hydrate formation on the properties of the flux function $f$ in \eqref{eq:chi}, and in particular on its variability across heterogeneous sediments.
We also compare the equilibrium and kinetic models: as expected, kinetic solution are close to the equilibrium solution as the kinetic exchange rate increases. 

In the last examples in Sec.~\ref{sec:newex} we illustrate the sensitivity of the model to the choice of macro time steps $\Delta T$ from Sec.~\ref{sec:macro}. 

\subsection{Examples for equilibrium model}
\begin{example}
[Model case {for equilibrium model} with analytical solution] \label{ex:masscons}
Let $\Omega = (-1,3), R = 2, \chi_L = 1, q = 1$, and the initial condition $u_{init}(x) = \chi_L {H(x+1)}H(-x)$ for \eqref{eq:chi} features a ``box''-like profile.  We consider $\chi^*(x) = e^{-0.5x}$ independent of time. For additional interest, we also consider $\chi^*(x) = 1-0.26x$. 
\end{example}

The analytical solution to \eqref{eq:chi} with $u_{init}(x) = \chi_L H(-x)$ can be found in \cite{PSW}. We modify it for the present case of ``box'' shaped $u_{init}(x)$
\bas
\chi(x,t) &=& \min\left(1,\frac{\chi^*(x)}{\chi_L}\right)u_{init}(x-qt),\\
S(x,t) &=& -\frac{\max(0,t-\frac{x}{q})q\chi^*_x(x)\mathbbm{1}_{G_0(t)}(x)}{R-\chi^*(x)},\\
u(x,t) &=& \chi(x,t) + (R-\chi^*(x))S(x,t),
\eas
where $G_0(t) = \{x:x_L<x \leq qt\}$ with $x_L$ satisfying $\chi^*(x_L) = \chi_L$, the position where first hydrate formation is observed. 

\begin{figure}
\centering
\includegraphics[width=.4\textwidth]{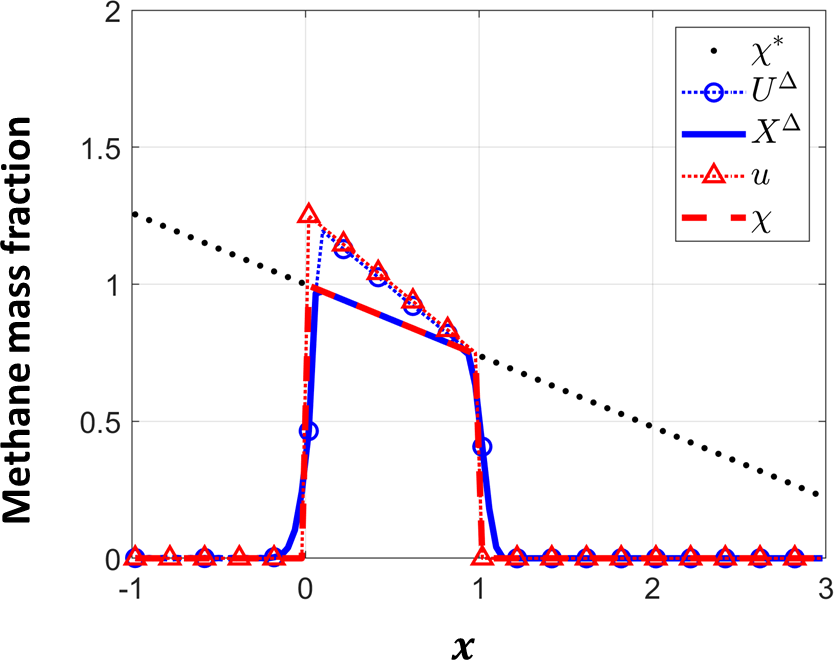}
\includegraphics[width=.4\textwidth]{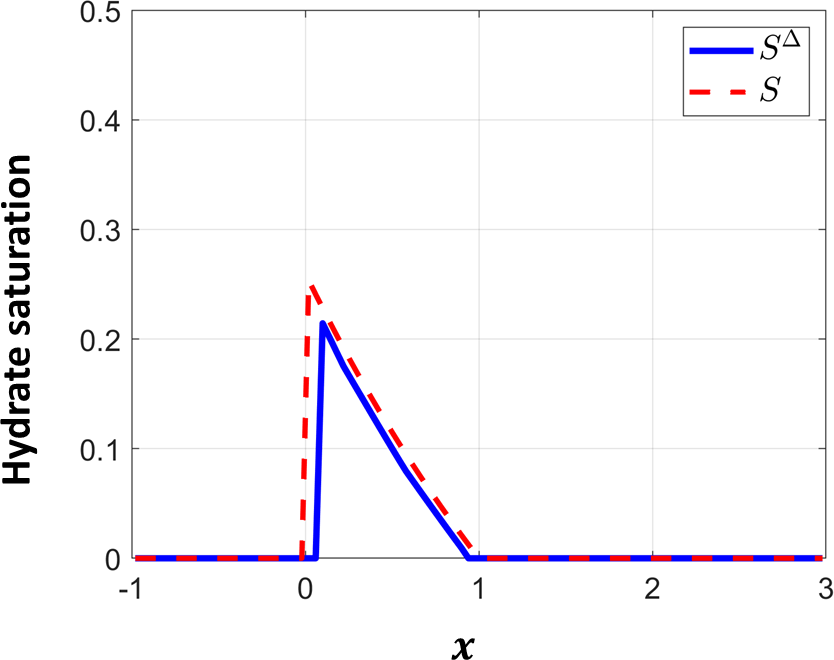}
\\
\includegraphics[width=.4\textwidth]{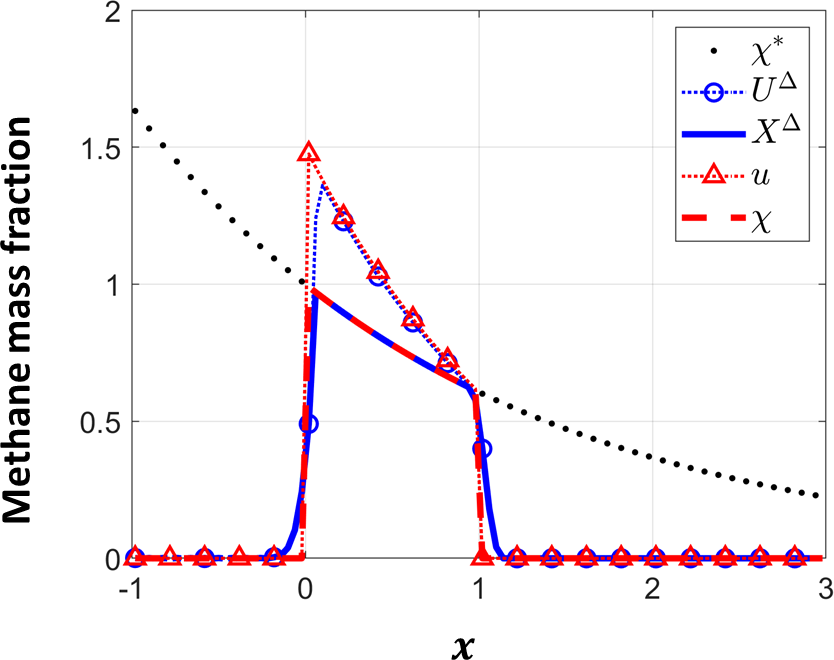}
\includegraphics[width=.4\textwidth]{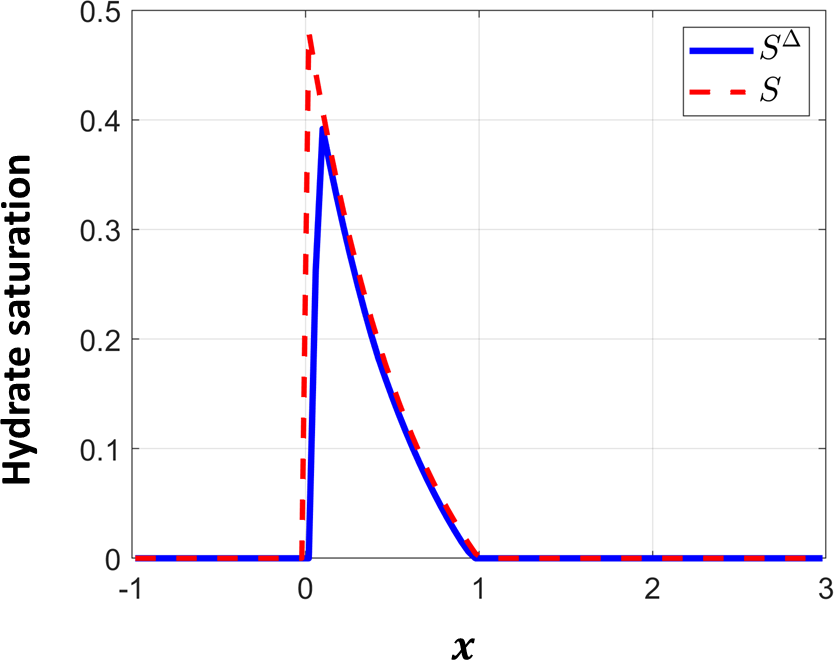}
\caption{Comparison of the numerical solution  $(U^{\Delta},X^{\Delta},S^{\Delta})$  with the analytical solution $(u,\chi,S)$  at $T=1$ with $M=100$ and $\nu = 0.9$ for Ex.~\ref{ex:masscons}.
Top: case with for $\chi^*(x) = 1-0.26x$. Bottom: case with  $\chi^*(x) = e^{-0.5x}$.
} 
\label{fig:masscons}
\end{figure}

We apply scheme \eqref{eq:upwindEQ} to obtain $(U^{\Delta},X^{\Delta},S^{\Delta})$ at $T=1$ with $M = 100$ and $\nu = 0.9$. Illustrations are provided in Fig.~\ref{fig:masscons}.
We see that $(U^{\Delta},X^{\Delta},S^{\Delta})$ are close to the analytical solution $(u,\chi,S)$. As $U^{\Delta}$ propagates to the right, $X^{\Delta}$ satisfies the constraint $X_j^n\leq \chi^*(x_j)$, and the undissolved methane produces $S_j^n>0$, i.e., we see the ``blow-up" behavior of $U^{\Delta}$ with $S_j^n$ as expected.

Comparing the two cases of $\chi^*(x) = 1-0.26x$ and $\chi^*(x) = e^{-0.5x}$, we see that the magnitude of $\chi^*_x$
is more pronounced for the latter case. 
In $U^{\Delta}$ and $S^{\Delta}$ we see small rarefactions at the back of the traveling wave, typical for an increasing concave flux function such as $f$ given by \eqref{eq:chif}.

\medskip
In our next example we evaluate effects of regularization in order to understand the closeness of $u$ and $u^{\epsilon}$, the solutions to \eqref{eq:chi} and \eqref{eq:regularization}. 
With $f^{\epsilon}$ chosen to be really close to $f$, we can make the difference between $U^{\Delta}$ and the solution to the regularized model $U^{\epsilon,\Delta}$ arbitrarily small. For comparison we use the case with the analytical solution from \cite{PSW}, which we adapt to the use of realistic data from Ulleung Basin. 

\begin{example}
\label{ex:regularization}
[Convergence rate and regularization; homogeneous domain, basin time scale] 
\label{ex:pureadv} We consider \eqref{eq:chi} with $u_{init}(x) = \chi_L H(-x)$ on $\Omega = (0,D^{\text{max}})$ where $D^{\text{max}}$ and $\chi^*$ are computed using the reference data measured from the Ulleung basin site UBGH2-7 of \cite{PHTK} with constant salinity of $\chi_{lS}^{sw}=3.5\%$. Let $q = 5\times 10^{-3}\mpunit{m/y}$, $\chi_L = 2\times 10^{-3}$, and $R = 0.1203$. We examine the result at $T = 10\mpunit{ky}$.
\end{example}

\begin{figure}
\centering
\includegraphics[width=.4\textwidth]{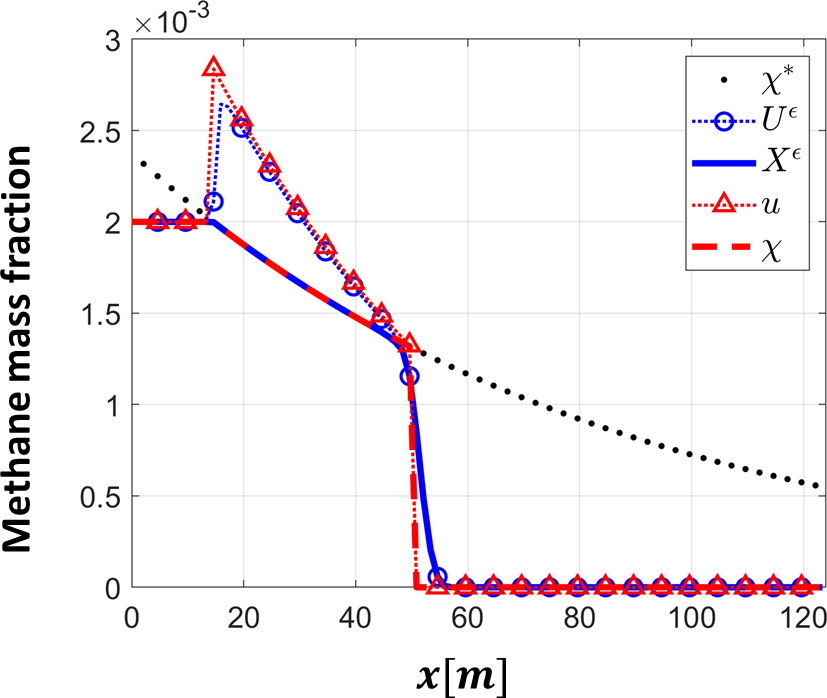}
\includegraphics[width=.4\textwidth]{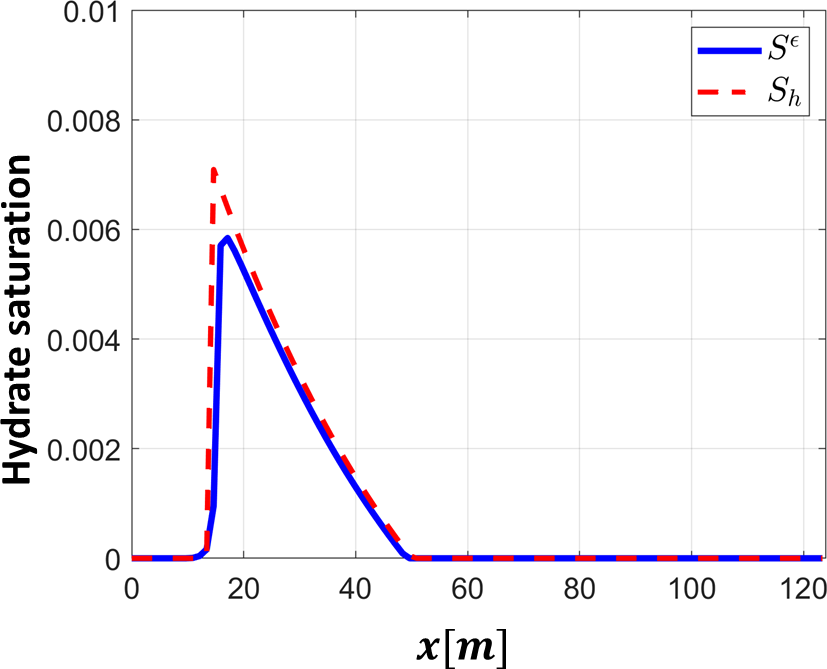}
\caption{Numerical solution $(U^{\epsilon},X^{\epsilon},S^{\epsilon})$ of Ex.~\ref{ex:pureadv} at $t = 10\mpunit{ky}$ with $M = 100$ compared with the analytical solution $(u,\chi,S_h)$.
\label{fig:pureadv}}
\end{figure}

The flux function $f(x;u)$ given by \eqref{eq:chif} has a  corner at $u = \chi^*(x)$ at every $x$. We regularize with $f^{\epsilon}$ which replaces $f$ on $(x,u)\in \Omega\times [\chi^*(x)-\epsilon,\chi^*(x) + \epsilon]$ by a smooth polynomial. Here $\epsilon$ is a regularization parameter; we choose $\epsilon=\alpha h$ with $\alpha = 10^{-4}$. 
In Fig.~\ref{fig:pureadv} we illustrate the analytical solution as well as the numerical solution $U^{\epsilon,\Delta}$ to the regularized problem, at $t = 10$ kyrs with $M=100$ and $\nu = 0.9$. We do not show $U^{\Delta}$ separately because it is virtually indistinguishable from $U^{\epsilon,\Delta}$.

We first examine the qualitative behavior. As predicted by the analytical solution, we observe the rapid growth of total methane content $U^{\epsilon}$ and the hydrate accumulation $S^{\epsilon}$ because $R < 1$ while $S^{\epsilon}$ is inversely related to $R$. 
At $t = 10\,\mpunit{ky}$, the hydrate saturation reaches about $10\%$.

Next we compare $U^{\Delta}$ and $U^{\epsilon,\Delta}$. With $\epsilon=O(h)$, their difference is small. In particular, when $M = 100$, $\|U^{\Delta}-U^{\epsilon,\Delta}\|_1 = 3.52\times 10^{-4}$, $\|X^{\Delta}-X^{\epsilon,\Delta}\|_1 = 1.32\times 10^{-4}$, and $\|S^{\Delta}-S^{\epsilon,\Delta}\|_1 = 2.20\times 10^{-3}$. When $M = 1000$, $\|U^{\Delta}-U^{\epsilon,\Delta}\|_1 = 9.41\times 10^{-6}$, $\|X^{\Delta}-X^{\epsilon,\Delta}\|_1 = 2.99\times 10^{-6}$, and $\|S^{\Delta}-S^{\epsilon,\Delta}\|_1 = 5.46\times 10^{-5}$.

\begin{figure}
    \centering
    \includegraphics[width=.49\textwidth]{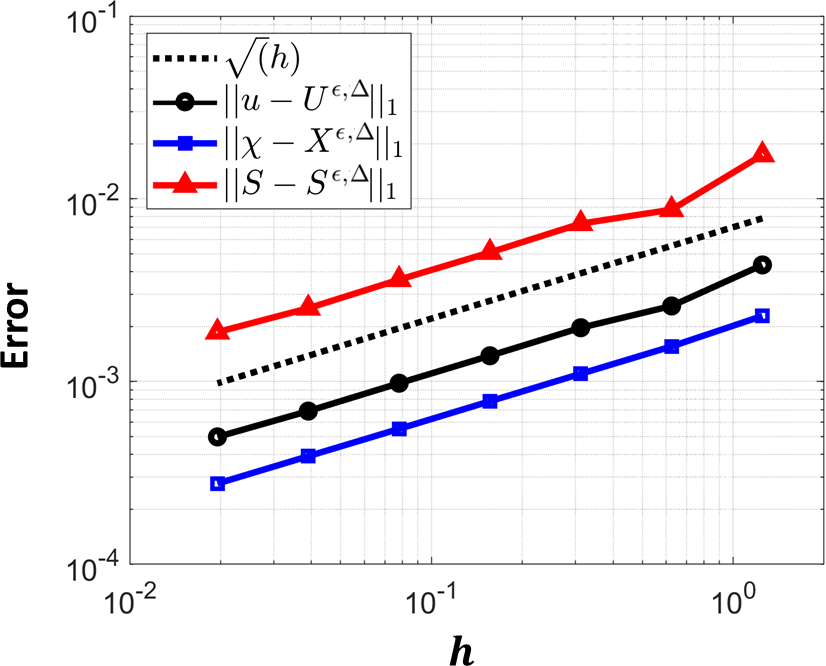}
    \includegraphics[width=.5\textwidth]{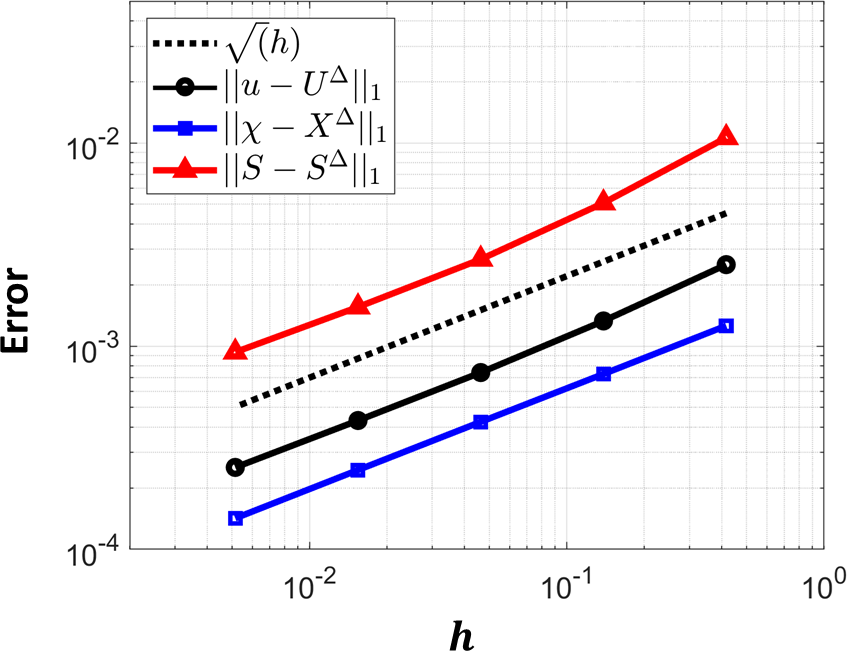}
    \caption{Left: $L_1$-error between the numerical solution $(U^{\epsilon,\Delta},X^{\epsilon,\Delta},S^{\epsilon,\Delta})$ and the analytical solution $(u,\chi,S)$ at $T=10\mpunit{ky}$; from  Ex.~\ref{ex:pureadv}.  Right: Convergence of the numerical solution $(U^{\Delta},X^{\Delta},S^{\Delta})$ to the analytical solution $(u,\chi,S)$ at $T = 10\mpunit{ky}$.} 
    \label{fig:pureadv_error}
\end{figure}

We also check the rate of convergence using a fine grid solution with $100\leq M \leq 6400$, plotted in Fig.~\ref{fig:pureadv_error}.
\bas
\|u-U^{\Delta}\|_1 = O(h^{0.52}),\;\; 
\|\chi-X^{\Delta}\|_1 = O(h^{0.5}),\;\;  
\|S-S^{\Delta}\|_1 = O(h^{0.55}).
\eas
The order is similar for the solutions to the regularized model, with the error slightly bigger due to the modeling error. 
We have 
\bas
\norm{u-U^{\epsilon,\Delta}}{1} = O(h^{0.51}),\;\; 
\norm{\chi-X^{\epsilon,\Delta}}{1} = O(h^{0.50}),\;\;  \norm{S-S^{\epsilon,\Delta}}{1} =O(h^{0.51}).
\eas

\medskip
Our next example challenges the theory since it is set for heterogeneous sediment. This example is inspired by \cite{DD11}; see our 2D simulation in \cite{P18} which accounts also for the flow and fracturing.

\begin{figure}
\centering
\includegraphics[width=.4\textwidth]{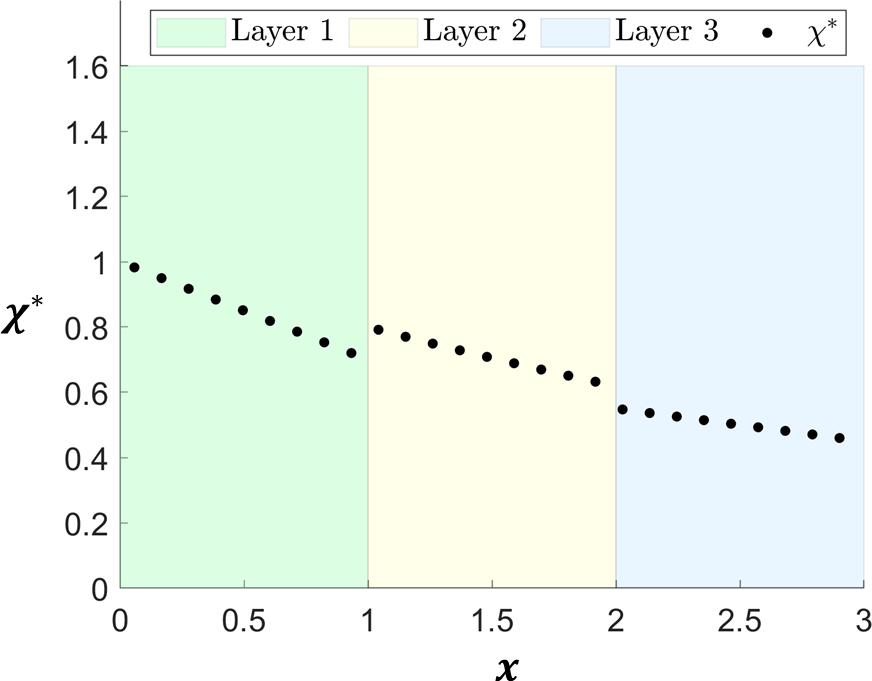}
\\
\includegraphics[width=.4\textwidth]{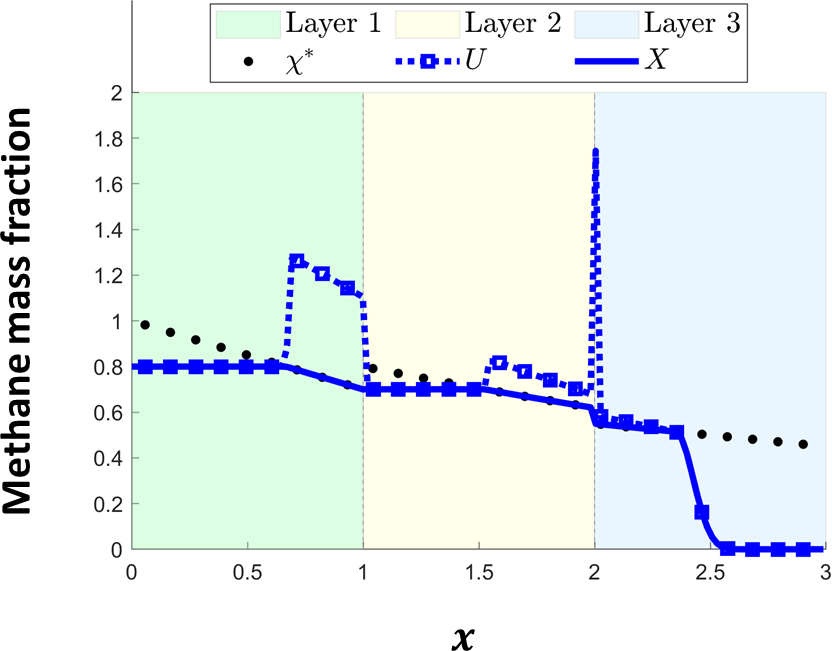}
\includegraphics[width=.4\textwidth]{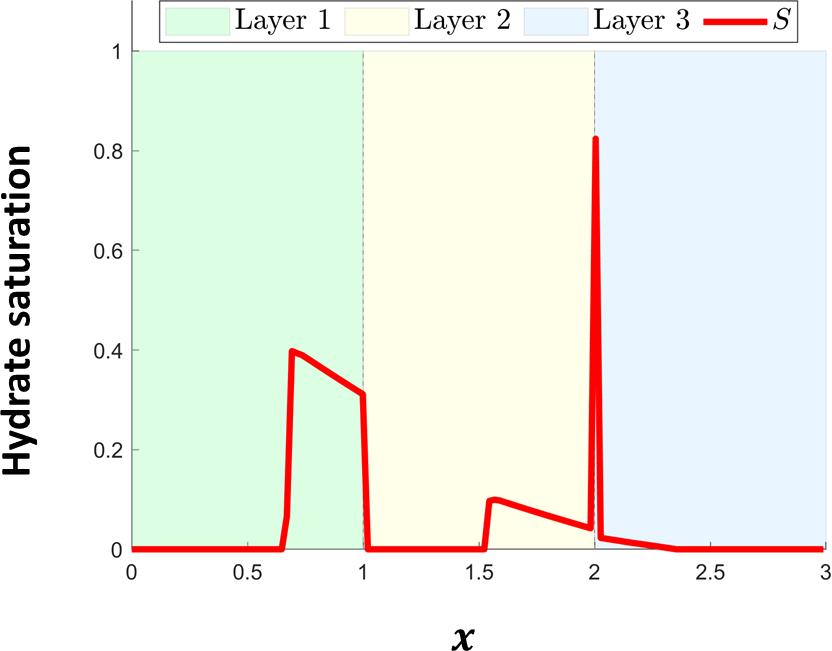}
\caption{Transport in  heterogeneous domain from Ex.~\ref{ex:adv_diff_heterogeneous_model} at $T=2.4$.  
Top: layers of heterogeneous sediment with different maximum solubility curves $\chi^*$ at $x = 1$ and $x= 2$. 
Bottom: numerical solution.
Of interest is behavior at the interfaces caused by the jumps of $\chi^*(x)$. \label{fig:multilayer_adv}}
\end{figure}

\begin{example}[Model problem in heterogeneous domain motivated by \cite{DD11}]\label{ex:adv_diff_heterogeneous_model}
\label{ex:hetero}
\mps Consider advection and diffusion of methane gas through 3 layers of sediments. Let $\Omega = \cup_{i=1}^3 \Omega_i$ where $i$ indicates each layer, each   with different methane solubility curves: $\chi_1^*(x) = -0.3x+1$, $\chi_2^*(x) = e^{-0.2(x-1)}-0.2$ and $\chi^*_3(x) = -0.1x +0.75$ shown in Fig.~\ref{fig:multilayer_adv}. 
  We use $R = 2$, $q = 1$, ${d_m} = 0$ and $u_{init}(x) = 0.8H(-x)$. 
\end{example}

The domain and the solutions are  illustrated in Fig.~\ref{fig:multilayer_adv} where the shaded blocks correspond to different layers.  We focus on the behavior near the interfaces at $x=1$ and $x=2$. 
As the front of methane enters Layer 2 at $x=1$, we expect to see methane hydrate {dissociation} since $\chi^*(1^+) > \chi^*(1^-)$ allows more methane gas to dissolve in the water.  In contrast, at $x = 2$, there is a reduction in maximum solubility; $\chi^*(2^+)<\chi^*(2^-)$: this cause a sudden formation of hydrate at the interface as in \cite{DD11,Rempel2011,VR2018}.
The simulation captures the hydrate {dissociation} at $x = 1$ and the formation at $x= 2$. The sharp spike at $x=2$ makes sense, since the weak derivative  $\partial_x f(x,u)$ at the discontinuity at $x= 2$ is a Dirac source $\delta(x-2)$.

\subsection{Examples for kinetic model}
Next we study convergence of the kinetic model, and compare the equilibrium model and kinetic models. Clearly kinetic rate should be fixed from experimental data; however,  we can investigate the case as $k_3$ increases to see how realistic it is to use (KIN3) in place of equilibrium models. As $k_3$ increases, we see the solutions to \eqref{eq:kinetic} become closer to and eventually indistinguishable from those for the equilibrium model.  

\begin{example}[{Model case: equilibrium and kinetic models}]\label{ex:kinetic}
Let $\Omega = (0,2), q= 1, R= 2$. We use $x_L = 0.35$, $\chi^*(x) = e^{-0.5x}$, and the initial condition $u_{init}(x) = \chi_LH(-x)$ with $\chi_L = 0.8395$. We simulate the problem using both the equilibrium model and scheme \eqref{eq:upwind}, and with the kinetic model and scheme \eqref{eq:schemek} when $k_3 = 10$ and $k_3=100$. Here $M = 100$ and {$\nu =q\tau/h=0.9$}. We compare with the equilibrium solution  at $T = 1$. 
\end{example}

\begin{figure}
\centering
\includegraphics[width=.4\textwidth]{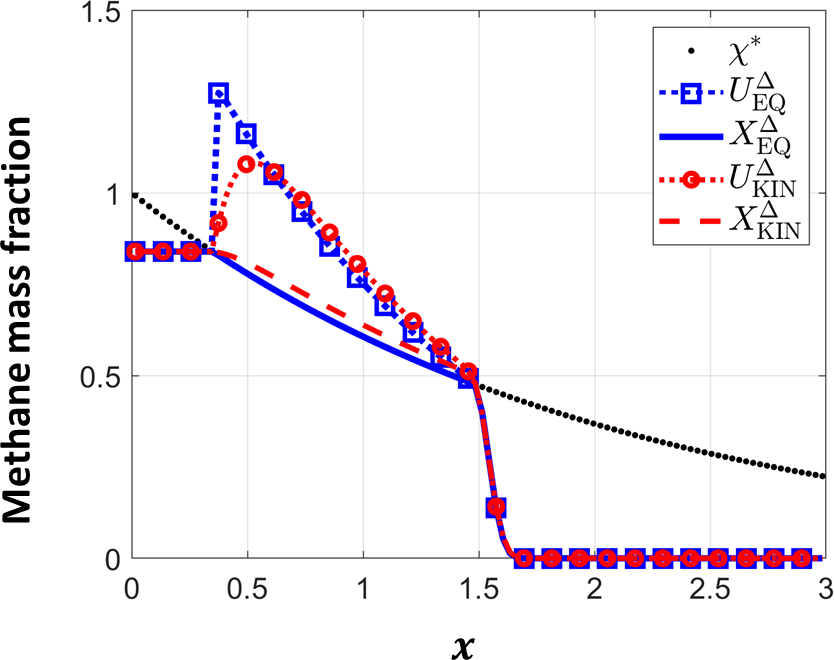}
\includegraphics[width=.4\textwidth]{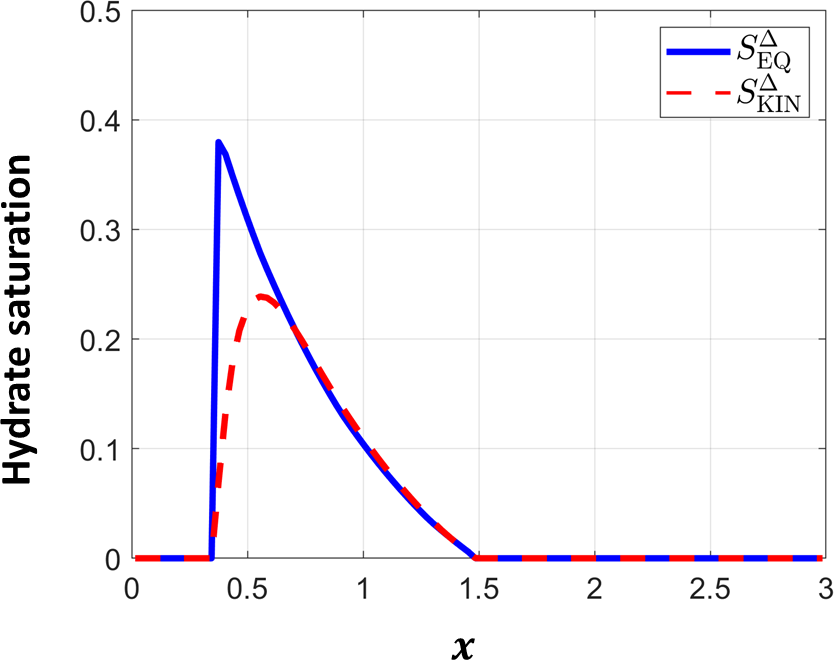}
\\
\includegraphics[width=.4\textwidth]{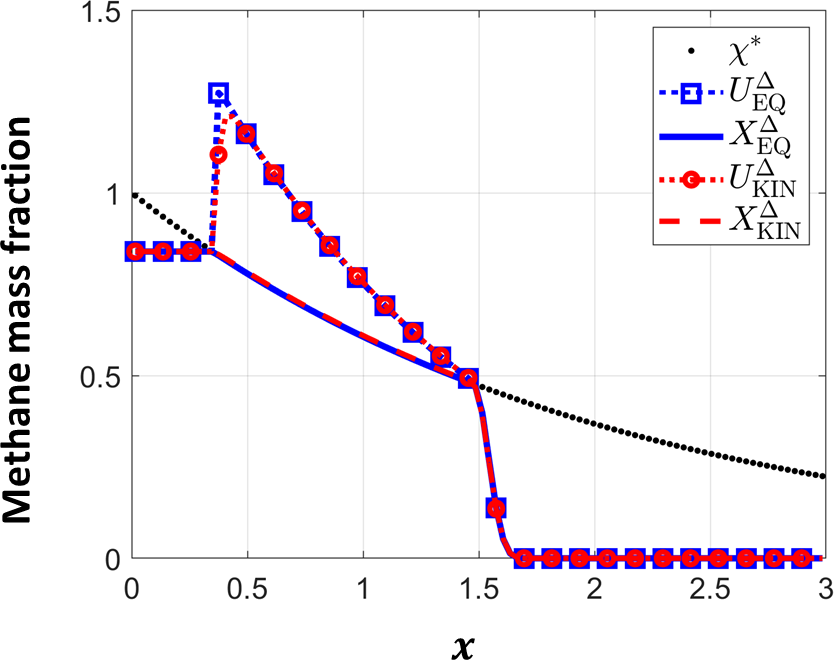}
\includegraphics[width=.4\textwidth]{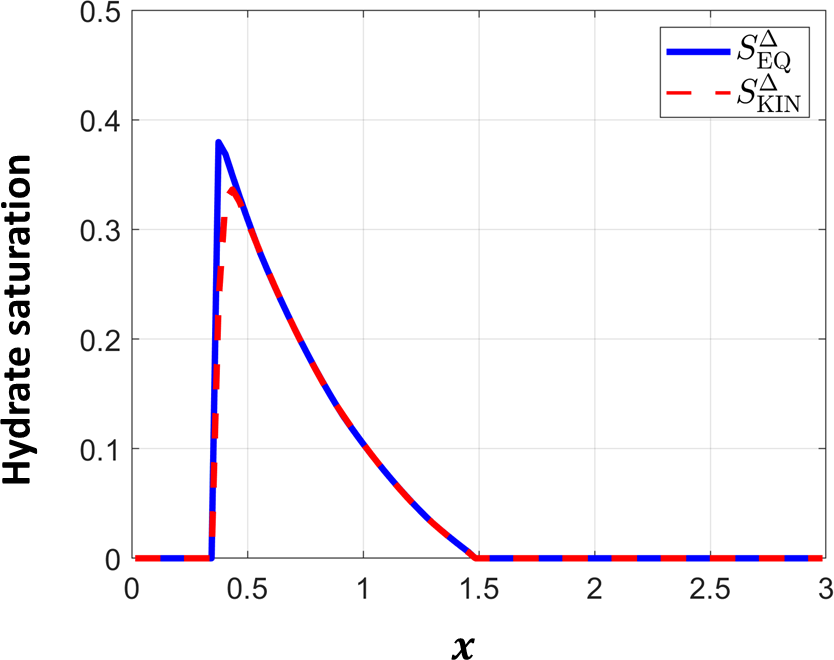}
\caption{Comparison of the kinetic and equilibrium model solutions at $T = 1$ with $M=100$, and $\nu = 0.9$ for Ex.~\ref{ex:kinetic}. Top: rate $k_3=10$. Bottom: rate $k_3=100$.}
\label{fig:kinetic_adv10}
\end{figure}

\begin{figure}
    \centering
    \includegraphics[width=.5\textwidth]{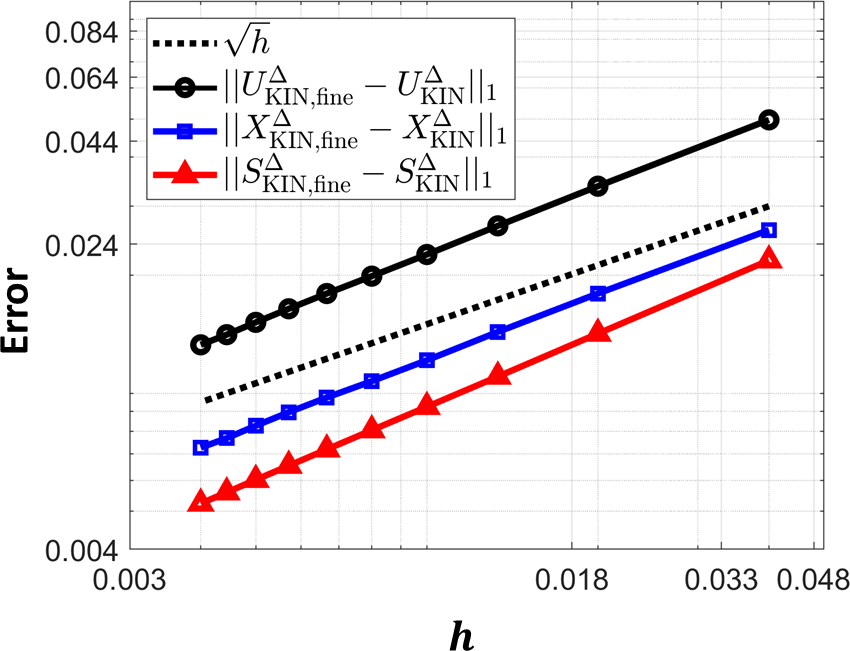}
    \caption{Convergence of the numerical solutions of Ex.~\ref{ex:kinetic} for $M = \{100,200,\dots,1000\}$ to the fine grid solution $M = 50000$ with $k_3 = 100$ at $T = 1$.
    \label{fig:KE_adv_convegence}}
\end{figure}

Fig.~\ref{fig:kinetic_adv10} illustrates the results.  We confirm that, as expected, the kinetic solution ``lives'' in the vicinity of the equilibrium solution. This closeness is more pronounced with larger $k_3$. 
In turn, Fig.~\ref{fig:KE_adv_convegence} shows that the numerical solutions $(U^{\Delta}_{\mathrm{KIN}},X^{\Delta}_{\mathrm{KIN}},S^{\Delta}_{\mathrm{KIN}})$ converges to the fine grid solutions $(U^{\Delta}_{\mathrm{KIN,fine}},X^{\Delta}_{\mathrm{KIN,fine}},S^{\Delta}_{\mathrm{KIN,fine}})$ at the order roughly of {$O(h^{0.5})$}. 
 \bas
 \|U^{\Delta}_{\mathrm{KIN,fine}}-U^{\Delta}_{\mathrm{KIN}}\|_1 &=& O(h^{0.57}),\;\;
 \\
 \|X^{\Delta}_{\mathrm{KIN,fine}}-X^{\Delta}_{\mathrm{KIN}}\|_1 &=& O(h^{0.56}),
 \;\;
\\
 \|S^{\Delta}_{\mathrm{KIN,fine}}-S^{\Delta}_{\mathrm{KIN}}\|_1 &=& O(h^{0.62}).
\eas


\subsection{Equilibrium and kinetic schemes under varying environmental and thermodynamic conditions and sensitivity to macro steps}
\label{sec:newex} 

{
Finally we illustrate the dependence of the solutions to the equilibrium model and kinetic models depending on the choice of macro steps $T^m=m \Delta T $ at which $\chi^*$ is recomputed. We allow  $P=P(x,t)$ and $T=T(x,t)$ to vary due to the changing environmental conditions and specifically due to the warming of ocean temperature and the sea level rise as predicted in \cite{MGT2016}, with the rate of sea level rise is $0.003\mpunit{m/y}$ and the rate of temperature rise at the seafloor of $0.01\mpunit{K/y}$. Then assume that the pressure (P), and temperature (T) at the seafloor vary linearly with respect to time $t\mpunit{y}$, with subscript $ref$ and $eq$ to denote the values at the seafloor and at the BHSZ, respectively. 
\bas
  P_{ref}(t) &=& \rho_l g D_{ref}(t),\; \;
  D_{ref}(t) = D_{ref}(0) + 0.003t,\\  
 T_{ref}(t) &=& T_{ref}(0) + 0.01t,
\eas
where $\rho_l\approx 1030\mpunit{kg/m^3}$ is the density of seawater, and $g = 9.8\mpunit{m/s^2}$. Over $150\mpunit{y}$ we see the sea level rise by $0.45 \mpunit{m}$ and $T_{ref}$ increase by $1.5\mpunit{K}$. Assume further that $(P(x,t),T(x,t))$ vary linearly in $\Omega$  
\bas
 P(x,t) &=& P_{ref}(t) + G_H(d_{sf}(x)-D_{ref}(t)),\\
 T(x,t) &=& T_{ref}(t) + G_T(d_{sf}(x)-D_{ref}(t)),
\eas
where $d_{sf}(x)$ is the depth below the sea level. We then recompute the equilibrium conditions at BHSZ using the parametric model from \cite{PHTK} responding to the increase in $T_{ref},P_{ref}$.
}

\begin{figure}
\centering
\includegraphics[width=.4\textwidth]{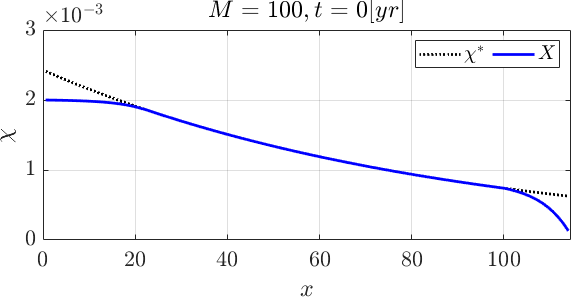}
\includegraphics[width=.4\textwidth]{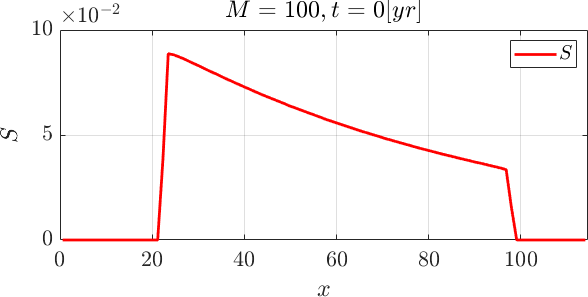}
\\
\includegraphics[width=.4\textwidth]{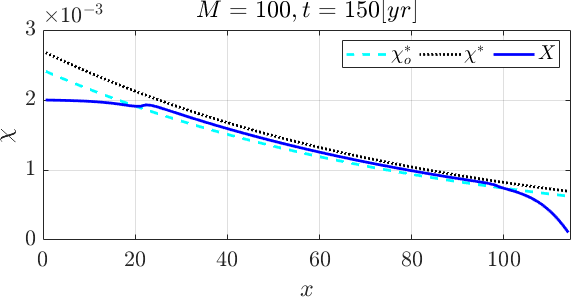}
\includegraphics[width=.4\textwidth]{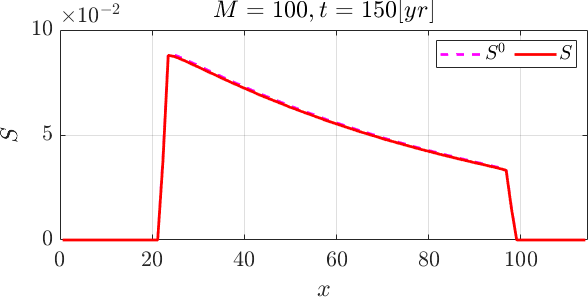}
\caption{Top: initial condition $\chi(x,0)$ and $S(x,0)$ from Sec.~\ref{sec:newex}. Bottom: results of Ex.~\ref{ex:warming} with  {$k_3 = 0.01$}. On both figures $\chi^*_0$ indicates the original $\chi^*\vert_{t=0}$. 
\label{fig:dissociation_IC}}
\end{figure}

{We apply these varying conditions to simulate the dissociation in a hydrate reservoir from the state obtained with simulation in Ex.~\ref{ex:regularization} run until $100\mpunit{ky}$. We consider this state to be the initial state for this simulation at $t=0$; see Fig.~\ref{fig:dissociation_IC} with $\chi_0^*$ as shown. The changes in $P$ and $T$ in time require we recompute $\chi^*=\chi^*(x,t)$ at the macro-time steps $T^m=m \Delta T$ as discussed in Sec.~\ref{sec:macro}. We adopt other parameters as in Ex.~\ref{ex:regularization} but use 
a fixed $\tau=1\mpunit{y}$ much smaller than that needed by CFL condition. We consider impact of $\Delta T=K\tau$, with $K=1,10,50,150$. }

\begin{example}[Hydrate dissociation due to warming waters]\label{ex:warming}
We start from the equilibrium state shown in Fig.~\ref{fig:dissociation_IC}. At every macro-time step $T^m=m \Delta T$, we recompute $\chi^*(x,T^m)$. Using the parameters as in Ex.~\ref{ex:regularization} {with $d_m = 3\times 10^{-2} \mpunit{m^2/y}$}, we simulate hydrate dissociation using the (EQ) model at $t\in [0,150]\,\mpunit{y}$, and plot the solutions at the final $t^N=150\mpunit{y}$.  We compare the results to the numerical solutions generated by the equilibrium model \eqref{eq:upwindqd} and the kinetic model \eqref{eq:kinetic}, both amended to include diffusion, depending on different rates {$k_3$} and the choice of $\Delta T$. 
\end{example}

At this time scale, dissociation proceeds slowly as shown by the decrease in the overall amount $\norm{S}{1}$ as well as the peak amount $\norm{S}{\infty}$; it is also interesting to test the magnitude $S_2$  of the last peak  before the decrease to seafloor. For all simulations the peak
$\norm{S}{\infty}$ is attained at $x=23.5 \mpunit{m}$, and the last saturation peak $S_2$ corresponds to $x=96.8\mpunit{m}$.

We find the difference between taking macro-time steps with $K=1$ up to $K=150$ very small for both equilibrium and kinetic models,  and the difference between kinetic model and equilibrium model indistinguishable when $k_3=100$. This rate is still about 100 times less than the rate used in \cite{Rempel2011,VR2018}.

\begin{table}
    \centering
    \begin{tabular}{|*{11}{c|}} \hline 
        $t\mpunit{y}$ & K
            & $P_{ref} \mpunit{MPa}$ & $T_{ref} \mpunit{K}$ & $D_{ref}\mpunit{mbsl}$ 
            & $P_{eq} \mpunit{MPa}$ & $T_{eq} \mpunit{K}$ & $D_{eq}\mpunit{mbsl}$ 
            & $100\norm{S}{\infty}$ & $100S_2$ & $\norm{S}{1}$ \\ \hline
        0 & -- & 21.6500 & 273.5500 & 2145.00 & 22.8849 & 294.6665 & 2268.49 & 8.8732 & 3.3524 & 4.3681\\ \hline 
        150 & 1 & 21.6545 & 275.0500 & 2145.45 & 22.7999 & 294.6364 & 2259.99 & 8.7212 & 3.2854 & 4.2848\\ \hline 
    \end{tabular}
    \caption{Simulation reference data and results generated by the equilibrium model (EQ) with $K=1$ for Ex.~\ref{ex:warming}.
    \label{tab:ex_eq}}
\end{table}

\begin{table}
    \centering
    \begin{tabular}{|*{11}{c|}} \hline 
    \multicolumn{2}{|c|}{ }    
            & \multicolumn{3}{c|}{{$k_3 = 0.01$}} 
            & \multicolumn{3}{c|}{$k_3 = 1$} 
            & \multicolumn{3}{c|}{$k_3 = 100$ and EQ} \\ \hline
        $t\mpunit{y}$ & $K$
            & $100\norm{S}{\infty}$ & $100 S_2$ & $\norm{S}{1}$
            & $100\norm{S}{\infty}$ & $100 S_2$ & $\norm{S}{1}$
            & $100\norm{S}{\infty}$ & $100 S_2$ & $\norm{S}{1}$ 
            \\ \hline
        150 & $1$ & 8.7946 & 3.3185 & 4.3295 & 8.7181 & 3.2844 & 4.2855 & 8.7212 & 3.2854 & 4.2848
        \\ \hline 
        150 & $10$ & 8.7905 & 3.3168 & 4.3276 & 8.7168 & 3.2839 & 4.2848 & 8.7212 & 3.2854 & 4.2847
        \\ \hline 
        150 & $50$ & 8.7736 & 3.3096 & 4.3196 & 8.7160 & 3.2836 & 4.2842 & 8.7212* & 3.2854 & 4.2841
        \\ \hline 
        150 & $150$ & 8.7381 & 3.2948 & 4.3044 & 8.7144 & 3.2831 & 4.2831 & 8.7213* & 3.2855 & 4.2830
        \\ \hline 
    \end{tabular}
\caption{Simulation results at $t=150 \mpunit{y}$ with kinetic model for Ex.\ref{ex:warming} with $\Delta t = K\tau$. 
(*) The results in the last macro-column differ are higher for the EQ model by one digit.  
    \label{tab:macrosteps150}}
\end{table}

%
\section{Conclusions and future work}
\label{sec:conclusions}

In this paper we considered equilibrium and kinetic phase behavior for hydrate in two-phase conditions typical in sediments above the bottom of Hydrate Stability Zone. Our objective was to study the stability of numerical models for {transport coupled with phase transitions} within  IMPES-like time stepping for pressure and temperature. We provided {rigorous} justification why the commonly used numerical scheme is stable and robust, and showed convergence with rate $O(\sqrt{h})$ consistent with that for monotone scheme and scalar conservation law in the presence of discontinuities. We also explained the presence of ``spikes'' of hydrate saturation similar to those observed in nature. 

In addition, we investigated robustness of a variety of kinetic models in the two phase liquid-hydrate conditions. Such models are needed, e.g., during sudden rearrangement of external controls on thermodynamic equilibria. Since the kinetic model popular in literature dubbed (KIN1) and its linear variant called (KIN2) work only in saturated conditions, we developed another model (KIN3) which is robust across the unsaturated and saturated conditions {and} is equivalent to (KIN2) in saturated conditions. We combined this model (KIN3) with the transport model and were able to show its numerical stability. 

The rigorous numerical analysis results we demonstrate are new; we are not aware of any other analysis of this type for hydrate models for either equilibrium and kinetic case. 

More work is needed. In particular, extensions of the analysis for both the equilibrium and extensions of (KIN3) model for three phase conditions are needed; this is subject to our current work.


%


\section{APPENDIX: Auxiliary results}
\label{sec:appendix}

\subsection{Proof of Proposition~\ref{th:stability} on the stability of upwind scheme \eqref{eq:upwind} for $u_t+f(x,t;u)_x=0$  when $f$ is smooth.}
\label{sec:app-EQ}
We adapt the proof \cite{RedLeveque}{[Chapter 12]} to the case when $f=f(x,t;u)$; we require the boundedness of $f_{xx}$ and $f_{xu}$ uniformly in time. 
The proof is broken to thee parts. First we bound the difference $\abs{\Delta U^n_j}$ between two adjacent values depending on $\abs{\Delta U^{n-1}_j}$ and $\abs{\Delta U^{n-1}_{j-1}}$. From this we conclude about $TV(U^n)$. Last we address $TV_T(U^{\Delta})$. 

\paragraph{\bf Local bounds on $\abs{\Delta U^n_j}$} 
We first subtract \eqref{eq:upwind} at $j-1$ from that at $j$ to get
\bas
\Delta U_j^n = \Delta U_j^{n-1} - \frac{\tau}{h}\big[\underbrace{F_j^{n-1}-F_{j-1}^{n-1}}_{(a)}\big] + \frac{\tau}{h}\big[\underbrace{F_{j-1}^{n-1} - F_{j-2}^{n-1}}_{(b)}\big],
\eas
where $\Delta U_j^{n*}$ denotes $U_j^{n*} - U_j^{n*-1}$ for $n* = n$ and $n-1$. Since $f$ is smooth, we can rewrite 
\bas
(a) &=& f(x_j, t^{n-1};U_j^{n-1}) - f(x_j,t^{n-1};U_{j-1}^{n-1})+f(x_j,t^{n-1};U_{j-1}^{n-1}) - f(x_{j-1},t^{n-1};U_{j-1}^{n-1})\\
&=& f_u(x_j,t^{n-1};\wt{U}_j^{n-1})\Delta U_j^n + f_x(\wt{x}_j,t^{n-1};U_{j-1}^{n-1})h
\eas
where $\wt{U}_j^{n-1} \in (U_{j-1}^{n-1},U_j^{n-1})$ and $\wt{x}_j\in (x_{j-1},x_j)$, and similarly
\bas
(b) = f_u(x_{j-1},t^{n-1};\wt{U}_{j-1}^{n-1})\Delta U_{j-1}^{n-1} + f_x(\wt{x}_{j-1},t^{n-1};U_{j-2}^{n-1}) h
\eas
where $\wt{U}_{j-1}^{n-1} \in (U_{j-2}^{n-1},U_{j-1}^{n-1})$ and $\wt{x}_{j-1}\in (x_{j-2},x_{j-1})$. After the substitution, we get
\begin{multline*}
\Delta U_j^n = \left(1-\frac{\tau}{h}f_u(x_j,t^{n-1};\wt{U}_j^{n-1})\right)\Delta U_j^{n-1}
+ \frac{\tau}{h}f_u(x_{j-1},t^{n-1};\wt{U}_{j-1}^{n-1})\Delta U_{j-1}^{n-1} \\
- \tau\big[\underbrace{f_x(\wt{x}_j,t^{n-1};U_{j-1}^{n-1})-f_x(\wt{x}_{j-1},t^{n-1};U_{j-2}^{n-1})}_{(c)}\big].
\end{multline*}
Applying mean value theorem to $f_x$ terms we rewrite (c) as
\bas
(c) &=& f_x(\wt{x}_j,t^{n-1};U_{j-1}^{n-1})-f_x(\wt{x}_j,t^{n-1};U_{j-2}^{n-1}) + f_x(\wt{x}_j,t^{n-1};U_{j-2}^{n-1})-f_x(\wt{x}_{j-1},t^{n-1};U_{J-2}^{n-1}),\\
&=& f_{xu}(\wt{x}_j,t^{n-1};\overline{U}_{j-1}^{n-1})\Delta U_{j-1}^{n-1} + f_{xx}(\overline{x}_j,t^{n-1};U_{j-2}^{n-1})(\wt{x}_j - \wt{x}_{j-1}),
\eas
where $\overline{U}_{j-1}^{n-1}\in (U_{j-2}^{n-1},U_{j-1}^{n-1})$ and $\overline{x}_j\in (\wt{x}_{j-1},\wt{x}_{j})\subseteq(x_{j-2},x_j)$. Next we substitute $(c)$ to get
\begin{multline*}
\Delta U_j^n = \left(1-\frac{\tau}{h}f_u(x_j,t^{n-1};\wt{U}_j^{n-1})\right)\Delta U_j^{n-1} + \frac{\tau}{h}f_u(x_{j-1},t^{n-1};\wt{U}_{j-1}^{n-1})\Delta U_{j-1}^{n-1} \\
- \tau\big[f_{xu}(\wt{x}_j,t^{n-1};\overline{U}_{j-1}^{n-1})\Delta U_{j-1}^{n-1} + f_{xx}(\overline{x}_j,t^{n-1};U_{j-2}^{n-1})(\wt{x}_j-\wt{x}_{j-1})\big].
\end{multline*}
Next we take the absolute value of both sides and apply the triangle inequality. Since the CFL condition \eqref{eq:tau} holds, we get
\begin{multline*}
    \abs{\Delta U_j^n} \leq \left(1-\frac{\tau}{h}f_u(x_j,t^{n-1};\wt{U}_j^{n-1})\right)\abs{\Delta U_j^{n-1}} 
    + \frac{\tau}{h}f_u(x_{j-1},t^{n-1};\wt{U}_{j-1}^{n-1})\abs{\Delta U_{j-1}^{n-1}}\\
    + \tau \abs{f_{xu}(\wt{x}_j,t^{n-1};\overline{U}_{j-1}^{n-1})\Delta U_{j-1}^{n-1}} + 2\tau h \abs{f_{xx}(\overline{x}_j,t^{n-1};U_{j-2}^{n-1})}
\end{multline*}
\paragraph{\bf Estimates on $TV(U^n)$.}
Now we take the sum over $j\in \mathbb{Z}$, keeping in mind the compact support of $U^{\Delta}$, which reduces any sums over $\mathbb{Z}$ to those over some finite set $\mathbb{Z}^0$. We obtain
\begin{multline}\label{eq:TVUn}
 TV(U^n) \leq TV(U^{n-1}) - \frac{\tau}{h}\sum_{j\in \mathbb{Z}^0} f_u(x_j,t^{n-1};\wt{U}_j^{n-1})\abs{\Delta U_{j-1}^{n-1}} + \frac{\tau}{h} \sum_{j\in \mathbb{Z}^0} f_u(x_{j-1},t^{n-1};\wt{U}_{j-1}^{n-1})\abs{\Delta U_{j-1}^{n-1}} \\
 + \tau \sum_{j\in \mathbb{Z}^0}\abs{f_{xu}(\wt{x}_j,t^{n-1};\overline{U}_{j-1}^{n-1})\Delta U_{j-1}^{n-1}} + 2\tau \sum_{j\in\mathbb{Z}^0}\abs{f_{xx}(\overline{x}_j,t^{n-1};U_{j-2}^{n-1})}h.
\end{multline}
Re-indexing the third term on the right-hand-side of \eqref{eq:TVUn}, the second and the third terms cancel each other. Using the definition of $L_1$, we have
\bas
TV(U^n) \leq TV(U^{n-1}) + \tau L_1 \sum_{j\in\mathbb{Z}^0} \abs{\Delta U_{j-1}^{n-1}}+ 2\tau L_1\sum_{j\in \mathbb{Z}^0} h.
\eas
Since $\abs{supp(f)}\leq \omega_S$, we have $\sum_{j\in \mathbb{Z}^0} h\leq \omega_S$. By re-indexing the second term, we get following:
\bas
TV(U^n) \leq TV(U^{n-1})(1+\tau L_1) + 2\tau L_1 \omega_S.
\eas
We repeat this inequality recursively to obtain
\bas
TV(U^n) \leq TV(U^0)(1+\tau L_1)^n + 2\tau L_1 \omega_S \sum_{k=0}^{n-1}(1+\tau L_1)^k.
\eas
From Bernoulli inequality, $1+\tau L_1\leq e^{\tau L_1}$, we get $(1+\tau L_1)^n\leq e^{n\tau L_1}\leq e^{TL_1}$ and we sum up the finite series to see 
that \eqref{th:TVstability_space} holds with $C_1(T) = TV(U^0)e^{TL_1} + 2\omega_S(e^{TL_1}-1)$. 

\paragraph{\bf Variation in time.} We rewrite \eqref{eq:upwind} as
\bas
U_j^n - U_j^{n-1} = -\frac{\tau}{h}\left[f_u(x_j,t^{n-1};\wt{U}_j^{n-1})\abs{\Delta U_j^{n-1}} + f_x(\wt{x}_j,t^{n-1}; U_{j-1}^{n-1}) h\right],
\eas
where $\wt{U}_j^{n-1}\in (U_{j-1}^{n-1},U_j^{n-1})$ and $\wt{x}_j\in (x_{j-1},x_j)$. Take the absolute values of both sides and apply the triangle inequality to get
\bas
\abs{U_j^n-U_j^{n-1}}\leq \frac{\tau}{h}L_2\left(\abs{\Delta U_j^{n-1}} + h\right).
\eas
Next, we multiply both sides by $h$ and sum over $j\in \mathbb{Z}^0$ to get
\bas
\norm{U^n-U^{n-1}}{1}\leq \tau L_2\left[TV(U^{n-1})+\omega_S\right].
\eas
Since $TV(U^n)\leq C_1(T)$ from \eqref{th:TVstability_space}, now \eqref{eq:gnorm} holds with $C_2(T) = L_2(C_1(T)+\omega_S)$.  Finally, to get \eqref{th:TVstability_time}, we combine \eqref{th:TVstability_space} and \eqref{eq:gnorm}, and obtain
\bas
TV_T(U^{\Delta}) \leq C_3(T)=\sum_{n=0}^{T/\tau} \tau (C_1(T)+C_2(T))=T(C_1(T)+C_2(T)).
\eas

\subsection{Properties for the kinetic model}
\label{sec:app-KIN}
In this Section we provide details of fully implicit schemes for models batch reactor models (KINj), $j=1,2,3$, respectively \eqref{eq:KE_model1}, \eqref{eq:KE_model2}, and \eqref{eq:KE_model3}. Our analysis motivates and supports the construction of the model (KIN3) which works across unsaturated and saturated conditions. Furthermore, our analysis helps to identify physically meaningful variables $(\chi,S)$ when working in non-isolated system, and to guide time stepping control. 
We define the discrete schemes in Sec.~\ref{sec:ckinetic}, and  analyze their solvability and properties of solutions in Sec.~\ref{sec:kin123-properties}. 
We illustrate the schemes and their properties in Sec.~\ref{sec:kin123-ex}.

Let each model (KINj) have its own rate $k_j>0$. We define $\okn{j}=\tau k_j$, and 
$\wt{k_j} = \frac{\okn{j}}{1+\okn{j}}$. 
We denote by $(X^{\infty},S^{\infty})$ the equilibrium values on graph $E$. 

We consider a uniform time step $\tau>0$, and ${t^n}=n\tau$, and we seek the approximations $X^n \approx \chi({t^n}),S^n \approx S({t^n})$  in one step $[{t^{n-1},t^n})$, using the initial conditions $X^{n-1},S^{n-1}$.  Other variables including $\Psi^n \approx \psi({t^n})$, and quantities such as $Q^n$, are denoted analogously. The total methane content $U^n = X^n + (R-X^n)S^n$. The solutions corresponding to model (KINj) are denoted with subscripts (KINj) e.g., we use notation $\kn{X}{j}{n}$. When more compact notation is desired, and there is no need to indicate the time step, we use simpler notation, e.g. $X_j$. When no distinction between models is needed, we drop subscript $j$, and denote the new time step value sought $X=\kn{X}{j}{n}$, while we set  the previous time step values equal $\chibar=\kn{X}{j}{n-1}$. With this notation, each scheme  advances $(\chibar,\sbar)$ to the new time step value $(X,S)$. 

\subsubsection{Discrete schemes for batch kinetic models and their properties}
\label{sec:ckinetic}
The schemes are fully implicit:  for (KIN2) and (KIN3) the solutions can be calculated with a closed formula, but (KIN1) requires an additional solvers. We prove various properties, and compare the models. 

\paragraph{\bf Discrete scheme for (KIN1).} 
Given $(\chibar,\sbar)=(\kn{X}{1}{n-1},\kn{S}{1}{n-1})$, find $(X,S)$. 
\bsub\label{eq:KIN1}
\ba
(1-S)X - (1-\sbar)\chibar &=& \okn{1}(\chi^*-X), \label{eq:KIN1_X0}\\
RS - R\sbar &=& \okn{1}(X-\chi^*) \label{eq:KIN1_S0}.
\ea
\esub

{\bf Solver for \eqref{eq:KIN1}:} the calculation of $(X,S)$ from \eqref{eq:KIN1} is coupled and not explicit. To get a useful formula,  we first calculate formally from \eqref{eq:KIN1_S0}
\ba\label{eq:KIN1_Sf}
S &=& \frac{\okn{1}(X-\chi^*)}{R} + \sbar.
\ea
Then we substitute \eqref{eq:KIN1_Sf} in \eqref{eq:KIN1_X0}, and rearrange to get a quadratic equation for $X$
\ba\label{eq:KIN1_Xf}
X &=& \frac{\okn{1}}{R}(R-X)(\chi^*-X) + \sbar X +(1-\sbar)\chibar.
\ea
The solvability of \eqref{eq:KIN1_Xf} is addressed in  Property (B) proven below; we also suggest a practical solver.  

\fbox{\parbox{.8\textwidth}{
{\bf (KIN1) summary:
Given $(\chibar,\sbar)=(\kn{X}{1}{n-1},\kn{S}{1}{n-1})$:
}  
\\
Solve \eqref{eq:KIN1_Xf} for $X$. 
\\
Calculate $S$ from \eqref{eq:KIN1_Sf}.
\\
Set $(\kn{X}{1}{n},\kn{S}{1}{n})=(X,S)$. 
}}

\paragraph{\bf Discrete scheme for (KIN2).}
Given $(\chibar,\sbar)=(\kn{X}{2}{n-1},\kn{S}{2}{n-1})$, calculate $\psibar=\sbar(R-\chibar)$, and find $(X,\Psi)$ for which
\bsub\label{eq:KIN2}
\ba
X - \chibar &=& \okn{2} (\chi^* - X), \label{eq:KIN2_X0}\\
\Psi-\psibar&=& \okn{2} (X-\chi^*) \label{eq:KIN2_S0}.
\ea
\esub

{\bf Solver for \eqref{eq:KIN2}:} since \eqref{eq:KIN2_X0} is linear, we rearrange to get
\ba
\label{eq:KIN2_Xf}
X &=& \wt{k_2}\chi^* +(1-\wt{k_2}) \chibar,
\ea
Substituting to \eqref{eq:KIN2_S0} we get
\ba
\label{eq:KIN2_Psif}
\Psi &=& \psibar + \wt{k_2}(\chibar-\chi^*).
\ea
After some algebra, we obtain also an explicit formula
\ba\label{eq:KIN2_Sf}
S &=& \frac{\okn{2}(\chibar-\chi^*) + (1+\okn{2})(R-\chibar)\sbar}{(R-\chibar)+\okn{2}(R-\chi^*)}.
\ea

\fbox{\parbox{.8\textwidth}{
{\bf (KIN2) summary: Given $(\chibar,\sbar)=(\kn{X}{2}{n-1},\kn{S}{2}{n-1})$:
} 
\\
Calculate $\psibar=\sbar(R-\chibar)$.
\\
Calculate $X$ from \eqref{eq:KIN2_Xf}, $\Psi$ from \eqref{eq:KIN2_Psif}, and $S$ from \eqref{eq:KIN2_Sf}.
\\
Set $(\kn{X}{2}{n},\kn{S}{2}{n})=(X,S)$. 
}}

\paragraph{\bf Discrete scheme for (KIN3).} 
Given $(\chibar,\sbar)=(\kn{X}{3}{n-1},\kn{S}{3}{n-1})$, calculate $\psibar=\sbar(R-\chibar)$, and find $(X,\Psi,W)$. 
\bsub\label{eq:KIN3}
\ba
X - \chibar &=& \okn{3}(W - X), \label{eq:KIN3_X0}
\\
\Psi - \psibar &=& \okn{3}(X - W),\label{eq:KIN3_Psi0}\\
W &\in& w_*(\Psi).
\ea
\esub

{\bf Solver for \eqref{eq:KIN3}:} At a first glance, the scheme is more complicated than \eqref{eq:KIN2}. However, we can exploit various properties of monotone graphs to simplify. First we calculate formally
\ba \label{eq:KIN3_X1}
X &=& \wt{k_3}W + (1-\wt{k_3})\chibar,
\ea
which we we substitute in \eqref{eq:KIN3_Psi0} and rearrange as 
\ba
\Psi+\okn{3}W = \psibar+\okn{3}(\wt{k_3}W+ (1-\wt{k_3})\chibar),\;  W\in w_*(\Psi).
\ea
After a few steps of algebra we get 
\ba
\label{eq:wdef}
\Psi+\wt{k_3}W
= \psibar+\wt{k_3}\chibar,\;  W \in w_*(\Psi).
\ea
Now we use the resolvent $\mathcal{R}^{w^*}_{\wt{k_3}}(\cdot)$ of $w^*$ as defined in \eqref{eq:re} to solve \eqref{eq:wdef} for $\Psi\in \mathrm{domain}(w_*)$
\bas
\Psi &=& \mathcal{R}^{w^*}_{\wt{k_3}}\left(\psibar + \wt{k_3}\chibar\right). 
\eas
Since this resolvent function has a simple form 
$
\mathcal{R}^{w^*}_{\wt{k_3}}(w) =(w-\wt{k_3}\chi^*)_+,
$
with a few more substitutions we get
\ba \label{eq:KIN3_Psif}
\Psi &=& (\psibar + \wt{k_3}(\chibar-\chi^*))_+, 
\ea
an explicit formula giving $\Psi$ in terms of $\chibar,\psibar$. Once $\Psi$ is known, we calculate the auxiliary variable $W$ from \eqref{eq:wdef} by back-substituting \eqref{eq:KIN3_Psif}, and we have $W=\frac{\psibar-\Psi}{\wt{k_3}}+\chibar$, thus $W=\chi^*$ if $\Psi\geq 0$, and $W=\chibar+\frac{\psibar}{\wt{k_3}}$ otherwise.  These calculations allow to calculate $X$ explicitly
\ba \label{eq:KIN3_Xf}
X =\wt{k_3}\chi^*+(1-\wt{k_3})\chibar,\;\; \mathrm{\ if\ } \Psi\geq 0, \mathrm{\ and\ }
X=\chibar+\psibar, \mathrm{\ otherwise}.
\ea 

\fbox{\parbox{.8\textwidth}{
{\bf (KIN3) summary: Given $(\chibar,\sbar)=(\kn{X}{3}{n-1},\kn{S}{3}{n-1})$: 
} 
\\
Calculate $\psibar=\sbar(R-\chibar)$.
\\
Calculate $\Psi$ from \eqref{eq:KIN3_Psif}.
\\
Given $\Psi$, calculate $X$ from \eqref{eq:KIN3_Xf}. 
\\
Calculate auxiliary variables $W=\frac{\psibar-\Psi}{\wt{k_3}}$, and $S = \frac{\Psi}{R-X}$.
\\
Set $(\kn{X}{3}{n-1},\kn{S}{3}{n-1})=(X,S).$}}

\subsubsection{Properties of schemes (KIN1), (KIN2), and (KIN3)}
\label{sec:kin123-properties}

Suppose that 
\ba
\label{eq:inita}
(\chibar,\sbar) \in D^0.
\ea
Also, denote $\psibar = (R-\chibar)\sbar$.  
Below we prove solvability of \eqref{eq:KIN1}, \eqref{eq:KIN2}, and \eqref{eq:KIN3} as well as analyze  qualitative properties of their solutions which we  arrange in a list (A-B-C-D-E). Since each of the schemes is a one-step scheme, it is sufficient to only consider the properties of one step solutions $(X,S)$ depending on $(\chibar,\sbar)$.

\bigskip
\textbf{Property (A): mass conservation.} 
\\
If the solutions to any scheme $j=1,2,3$ exist, they satisfy $U^n =u(X^n,S^n)=U^0$ where $u(X,S)$ is given by \eqref{eq:udef}. In other words, for each scheme, the solutions $(X^n,S^n)$ stay on the curve $u(X^n,S^n)=U^0$, and we have
\ba \label{eq:ALLS}
S_j =\frac{U^0-X_j}{R-X_j}=\frac{\chibar+(R-\chibar)\sbar - X_j}{R-X_j}.
\ea
The map $S_j=S_j(X_j)$ is smooth and invertible when $0\leq X_j<R$. 
\begin{proof}
The first part follows immediately by adding the two equations defining each scheme for one step, and following for $n>0$ inductively.  Analysis of \eqref{eq:ALLS} is straightforward. 
\end{proof}

\textbf{Property (B): solvability of schemes.} 
\\
Schemes (KIN2), (KIN3) are uniquely solvable, and (KIN1) is solvable depending on data and if $\tau$ is small enough. 
\begin{proof}
The solutions to schemes (KIN2) and (KIN3) can be calculated from explicit algebraic expressions depending on the data $(\chibar,\sbar)$, thus the conclusion is immediate. 

However, scheme (KIN1) \eqref{eq:KIN1} requires a solution to the quadratic equation \eqref{eq:KIN1_Xf} which we frame 
as $p(X)=0$ with 
\ba
\label{eq:fixedpt}
p(X)&=& \frac{\okn{1}}{R}(R-X)(\chi^*-X) + \overline{S}X + (1-\overline{S})\overline{X}-X.
\ea
We see that $p(\cdot)$ in \eqref{eq:fixedpt} is a quadratic polynomial, with $p(0)=\okn{1}\chi^*+(1-\sbar)\chibar$, and $p(R)=\sbar R -R + (1-\sbar)\chibar$. Also, $p'(X)=\frac{2\okn{1}}{R}(X-X_R)+\sbar-1$ where $X_R=\frac{R+\chi^*}{2}$. 
Now consider the root of $p(X)=0$.  From 
\eqref{eq:inita} we have that $p(0)>0$, and $p(R)=(1-\sbar)(\chibar-R)<0$. Since $p(\cdot)$ is continuous, we see that the root to $p(X)=0$ exists in $[0,R]$ and in fact is in $[0,R)$. Since, in addition, $p(\cdot)$ is convex, with $p''(X)=\frac{2\okn{1}}{R}>0$, we find that this root is unique in $[0,R)$, and is given from the quadratic formula 
\ba\label{eq:KIN1_Xsol}
X &=& \frac{R(1-\sbar)+\overline{k_1}(R+\chi^*) - \sqrt{(R(1-\sbar)+\overline{k_1}(R+\chi^*))^2 - 4\overline{k_1}R(\overline{k_1}\chi^*+\chibar(1-\sbar))}}{2\overline{k_1}}.
\ea
On the other hand, the second root given by a modification of \eqref{eq:KIN1_Xsol} always exists in $(R,\infty)$, but is unphysical. 
\end{proof}

\begin{lemma} \label{lemma:propC1} 
Consider (KIN1) scheme. Suppose \eqref{eq:inita} holds and consider the smaller root $X=X_1\in [0,R)$ of \eqref{eq:fixedpt}. Then (i) $S_1<1$. 
\\
(ii) If $\chibar<\chi^*$, then $\chibar<X<\chi^*$. If $\chi^*<\chibar$, then $\chi^*<X<\chibar$. If $\chibar=\chi^*$, then $X=\chibar$, and $S=\sbar$. 
\\
(iii)  If $(\chibar,\sbar) \in \Dp$
then $S_1\geq 0$ for any $\okn{1}$.
\\
(iv) On the other hand, suppose  $(\chibar,\sbar) \in \Dm$. If $\sbar=0$, then $S_1<0$. If $\sbar>0$, then for large $\okn{1}$ it is possible that $S_1 \leq 0$.
\end{lemma}
\begin{proof}
 To prove (i), we recall Property (A) and \eqref{eq:ALLS}. Since $u<R$, we have $S_1 =\frac{U^0-X_1}{R-X_1}<1$. 

 To prove (ii), assume $\chibar<\chi^*$.  First we collect the terms of \eqref{eq:fixedpt} with $X$ on the left hand side of the equation. Then subtract $(1-\sbar)\chi^*$ on the both sides of the equation to get
\bas
(\chi^*-X){\left[(1-\overline{S})+\frac{\okn{1}}{R}(R-X)\right]} &=& {(1-\overline{S})}(\chi^*-\overline{X}).
\eas
Since $\sbar<1$, and $X<R$, the second term on the left hand side and the first term on the right hand side are positive. Thus (i) the sign of $\chi^*-X$ is the same as that of $\chi^* -\chibar$. Further, rearrange $p(X)=0$ as in the proof of Property B to read
\bas
\frac{\okn{1}}{R}(R-X)(\chi^*-X) &=& (1-\overline{S})(X - \overline{X}).
\eas
Similarly as above we conclude that the sign of $X-\chibar$ is the same as that of $\chi^*-X$, which completes the proof of (ii). 

To prove (iii), take some $(\chibar,\sbar) \in \Dp$ so that $u(\chibar,\sbar)\geq \chi^*$. By property (ii), we can have $\chibar > X>\chi^*$, or $\chibar<X<\chi^*$. (We omit the trivial case $X=\chibar$). In the first case by \eqref{eq:KIN1_Sf} we have $S_1=\frac{\okn{1}(X-\chi^*)}{R} +\sbar\geq \sbar\geq 0$. In the second case by property A, $(X,S_1)$ is on the curve $u(X,S_1)=u(\chibar,\sbar)\geq \chi^*$ which is above the curve $u(X,S)=\chi^*$. Thus $S_1\geq 0$. 

To prove (iv), take $(\chibar,\sbar) \in \Dm$, so we must have $\chi<\chi^*$. 
With $\sbar=0$ we get 
from \eqref{eq:KIN1_Sf} that $S_1=\frac{\okn{1}(X-\chi^*)}{R} <0$. Even if $\sbar>0$, it is possible to find $\okn{1}$ large enough so that $\sbar\leq \frac{\okn{1}(\chi^*-\chi)}{R}$. For illustration, the curve $u(X,S_1)=u(\chibar,\sbar)\leq \chi^*$ is in this case below the curve $u(X,S)=\chi^*$, thus $S_1\leq 0$. 
\end{proof}

\textbf{Property C: the solutions to (KINj) stay physically meaningful: $(X_j,S_j) \in D^0$, under some conditions for $j=1,2$ and unconditionally for $j=3$.}
\\
(i) The solutions $(X_j,S_j)$ satisfy $0 \leq X_j<R$ and $S_j <1$ for all $j$. 
\\
(ii) In addition, $S_3 \geq 0$ unconditionally. 
\\
(iii) Let $j=1,2$. If $(\chibar,\sbar)\in \Dp$, then $S_j\geq 0$.  If $(\chibar,\sbar)\in \Dm$ and either  $\sbar>0$ with $\okn{j}$ small enough, then $S_j >0$.   However, if $(\chibar,\sbar)\in \Dm$ and either $\sbar=0$, or if $\okn{j}$ is large, then it is possible that $S_j<0$.  
\begin{proof}[Proof of property C]
(i) For (KIN1) Property B shows (i) for the correct root selected by the solver. 
For (KIN2), \eqref{eq:KIN2_Xf} shows that $X_2$ is a convex combination of $\chibar$ and $\chi^*$, thus $0\leq X_2<R$. 
For (KIN3), when $\sbar+\wt{k_3}(\chibar - \chi^*)\geq 0$, $X_3$ is a convex combination of $\chibar$ and $\chi^*$, and the same argument applies. Otherwise, $X_3=\chibar + \sbar(R-\chibar)<R$, and $0\leq X_3<R$.   
To prove $S_j<1$, we see that $X_j < R$ and by Property A  $u(X_j,S_j)=u(\chibar,\sbar)<R$, thus from \eqref{eq:ALLS} it follows that $S_j<1$. 

To show (ii) consider (KIN3) first. we have $\Psi_3\geq 0$ from \eqref{eq:KIN3_Psif}, thus $S_3\geq0$. For scheme (KIN1), we use Lemma~\ref{lemma:propC1}. For scheme (KIN2), we  recall \eqref{eq:KIN2_Psif}. To check if $0 \leq \Psi_2=\psibar+\wt{k_2}(\chibar-\chi^*)$, we first consider $\chibar\geq \chi^*$ (which implies  $(\chibar,\sbar)\in \Dp$). This yields $\Psi_2\geq \psibar \geq 0$, thus $S_2\geq0$. With  $\chibar<\chi^*$ however, we find that to guarantee $\psi\geq 0$, we must have $\wt{k_2} < \frac{\psibar}{(\chi^*-\chibar)}$. For these, we recall $\frac{\psibar}{(\chi^*-\chibar)}=\frac{u-\chibar}{(\chi^*-\chibar)}$, and this quantity $\frac{u-\chibar}{(\chi^*-\chibar)}\geq 1$ in $\Dp$, while we have that for any $\tau$, $0<\wt{k_2}<1$. We conclude that (KIN2) can produce unphysical $S_2\leq 0$ only for $(\chibar,\sbar) \in \Dm$. If $\psibar=0$, we always have $\Psi_2<0$ and $S_2<0$. 
\end{proof}
\textbf{Property D: stability of each scheme in $Q$}.
\\
We have $\abs{Q_j} <\abs{\qbar}$ for each scheme.  
\begin{proof}[Proof of property D]
We recall that $Q_j =k_j( X_j-\chi^*)$ for $j=1,2$, and $Q_3=k_3( X_3-W)$ for (KIN3). 
We consider the bounds for $j=1$ and (KIN1) first. We want to show $\abs{\chi^*-X_j} \leq \abs{\chi^* -\chibar}$. To this aim, we subtract $\chi^*$ from both sides of \eqref{eq:KIN1_Xf}, rearrange, and add $-\sbar \chi^*$ to both sides, and rearrange again to get
\bas 
\left(1-\sbar+\frac{\okn{1}}{R}(R-X)\right)(\chi^*-X_1) &=& (1-\sbar)(\chi^*-\chibar).
\eas
Next we take absolute value of both sides while we multiply them by $k_1$. Since $1-\sbar>0$ and  $\frac{\okn{1}}{R}(R-X) > 0$ from property C, we get, as desired 
\bas
\abs{Q_1} < \frac{1-\sbar}{1-\sbar+\frac{\okn{1}}{R}(R-X)}\abs{\qbar_1}<\abs{\qbar_1}.
\eas
For (KIN2),  we add $\chi^*-X_2$ to both sides of \eqref{eq:KIN2_X0} to get
\bas
\qbar=(1+\okn{2})Q_2.
\eas
By $1+\okn{2} >1$ it is easy to see $\abs{Q_2}<\abs{\qbar}$.

For (KIN3), the proof $\abs{\qbar_3}<\abs{Q_3}$ is a special homogeneous case of a more general proof.  
We first consider Yosida approximation $w_{\lambda} \approx w_*$, or some other regularization which maintains the monotonicity properties of the graph $w_*$. Given $(\chibar,\psibar)$ we seek the solution $(X_{\lambda},\Psi_{\lambda})$ to the regularization  of \eqref{eq:KIN3} with $Q_{\lambda} = X_{\lambda} - w_{\lambda}(\Psi_{\lambda})$, and $\qbar=\chibar-w_{\lambda}(\psibar)$

\bsub
\ba\label{eq:KIN3_Q}
X_{\lambda} - \chibar + \okn{3}Q_{\lambda} &=& 0, \label{eq:KIN3_Q_X0}\\
\Psi_{\lambda} - \psibar -\okn{3}Q_{\lambda} &=& 0. \label{eq:KIN3_Q_Psi0}
\ea
\esub
Next we multiply \eqref{eq:KIN3_Q_Psi0} by $w_{\lambda}'(\Psi')$ with some $\Psi'\in (\psibar,\Psi_{\lambda})$ to get
\ba \label{eq:KIN3_Q_Psi1}
w_{\lambda}(\Psi_{\lambda}) - w_{\lambda}(\psibar) - \okn{3}w_{\lambda}(\Psi')Q_{\lambda} &=& 0.
\ea
Subtract \eqref{eq:KIN3_Q_Psi1} from \eqref{eq:KIN3_Q_X0}, and take absolute value to get
\bas 
\left(1+\okn{3}(1+w_{\lambda}(\Psi')\right)\abs{Q_{\lambda}} = \abs{\qbar}.
\eas
Since $1+\okn{3}(1+w_{\lambda}(\overline{\Psi}))>1$, we get the inequality $\abs{Q_{\lambda}}<\frac{1}{1+\okn{3}}\abs{\qbar}$ as desired. Taking the limit as $\lambda \to 0$ gives the desired result. 
\end{proof}
\textbf{Property (E): (conditional) equivalence of schemes.}
\\
{\bf (i)} Let $k_2 =k_3$. The schemes (KIN2) and (KIN3) give the same numerical solutions $X_2=X_3$ iff $(\chibar,\sbar) \in \Dp$. 
{\bf (ii)} Moreover, if
\ba\label{eq:PropE_rate}
k_1 = \frac{R\okn{2}(1-\sbar)}{(R-\chibar) + \wt{k_2}(R-\chi^*)},
\ea
then the one-step solution to (KIN2) coincides with that for (KIN1).  

\begin{proof}
To prove (i), we want to check if $X_2 =X_3$ by setting the right hand sides of \eqref{eq:KIN2_Xf} and \eqref{eq:KIN3_X1} equal to each other. This identity holds if 
$
\chi^* = \min\left(\chibar+(1-\wt{k})\psibar,\chi^*\right),
$
which is equivalent to 
\ba
\label{eq:ktau}
\sbar \geq \frac{\wt{k}(\chi^*-\chibar)}{R-\chibar}. 
\ea
Now, if $(\chibar,\sbar) \in \Dp$, we have $U=\Ubar\geq \chi^*$, which means $\sbar\geq \frac{\chi^*-\chibar}{R-\chibar}\geq 
\frac{\wt{k}(\chi^*-\chibar)}{R-\chibar}$ for any $k$ and $\tau$. Conversely, for \eqref{eq:ktau} to hold, we must have $U\geq \chi^*$ since $\wt{k}$ can be made arbitrarily close to 1. 

To prove (ii), we want to calculate $k_2$ in terms of $k_1$ and previous time step data $(\chibar,\sbar)$.  Of course this is, in principle, always possible; the difficulty is to actually find this expression explicitly. We are able to do this and to obtain \eqref{eq:PropE_rate}. We explain the process below.

Recall from \eqref{eq:ALLS} that $S_j$ is a well defined invertible function of $X_j\in [0,R)$. In addition, for each scheme $j$, clearly each 
$(X_j,S_j)$ is a function of $(\chibar,\sbar)$ and of $k_j$.
If these were given explicitly, one could write, e.g., $S_1=S_2$ and attempt to solve for the dependence of $k_1$ on $k_2$ explicitly.  Alternatively, one could do the same starting with $X_1=X_2$ to get $k_1$ in terms of $k_2$. 
However, the solver for (KIN1) does not give either $X_1$ not $S_1$ explicitly depending on $k_1,\chibar,\sbar$, and these direct strategies fail. 

Instead, another possibility arises: we calculate $k_1=k_1(\chibar,\sbar;X_1)$, with $\frac{\partial k_1}{\partial X_1} \neq 0$, after some analysis. Next we assume $X_1=X_2$ and substitute $X_2=X_2(k_2,\chibar,\sbar)$
from \eqref{eq:KIN2_Xf}. With this, we get an expression with $k_1$ in terms of $(k_2,\chibar,\sbar)$ which is luckily explicit.

To get $k_1=k_1(\chibar,\sbar;X_1)$, we recall \eqref{eq:KIN1_Sf} which binds together the constants $k_1,\sbar$ and variables $X_1,S_1$. With \eqref{eq:ALLS} we eliminate $S_1$, and get a relationship between $k_1,\sbar,\chibar$ and $X_1$, and we solve for $\overline{k_1}$
\ba
\overline{k_1} = \frac{R(\chibar-X_1)(1-\sbar)}{(R-X_1)(X_1-\chi^*)}.\label{eq:E2_X1}
\ea
Now we assume $X_1=X_2$, recall \eqref{eq:KIN2_Xf} in which $X_2$ is given $X_2=X_2(\chibar,\sbar,k_2)$ explicitly and substitute this expression into \eqref{eq:E2_X1} to get
\bas
\overline{k_1}  = \frac{R(\chibar-\wt{k}_2\chi^*-(1-\wt{k}_2)\chibar)(1-\sbar)}{(R-\wt{k}_2\chi^* -(1-\wt{k}_2)\chibar)(\wt{k}_2\chi^*+(1-\wt{k}_2)\chibar-\chi^*)}
                = \frac{R\overline{k_2}(1-\sbar)}{(R-\chibar)+\wt{k}_2(\chibar-\chi^*)},
\eas
which, upon some algebra, is equivalent to \eqref{eq:PropE_rate}. 
\end{proof}

\subsubsection{Illustration of (KIN1), (KIN2), (KIN3) in batch setting.}
\label{sec:kin123-ex}

We illustrate now the three kinetic models with some numerical experiments. Our goal is to emphasize the similarities as well as the differences between these models. We employ the fully implicit schemes presented in Sec.~\ref{sec:ckinetic}.
In the examples we use data $R = 2$, $\chi^* = 1$, and $k_j=1$ for all $j=1,2,3$.

\begin{example}[Saturated case]
\label{ex:sat}
Suppose that $u^0 = 1.64>\chi^*$, thus the equilibrium state is  $(\chi^{\infty},S^{\infty}) =(1,0.64)$ on the $\Ep$ portion of the $E$ graph, and  this example falls in the saturated regime.  
We consider two  cases (I) $(X^0,S^0) = (0.2,0,8)$ and (II) $(X^0,S^0) = (1.4,0.4)$. Both are in saturated regime $(X^0,S^0) \in \Dp$. We use $\tau=1$. 
\end{example}
\begin{figure}
    \centering
    \includegraphics[width=.49\textwidth]{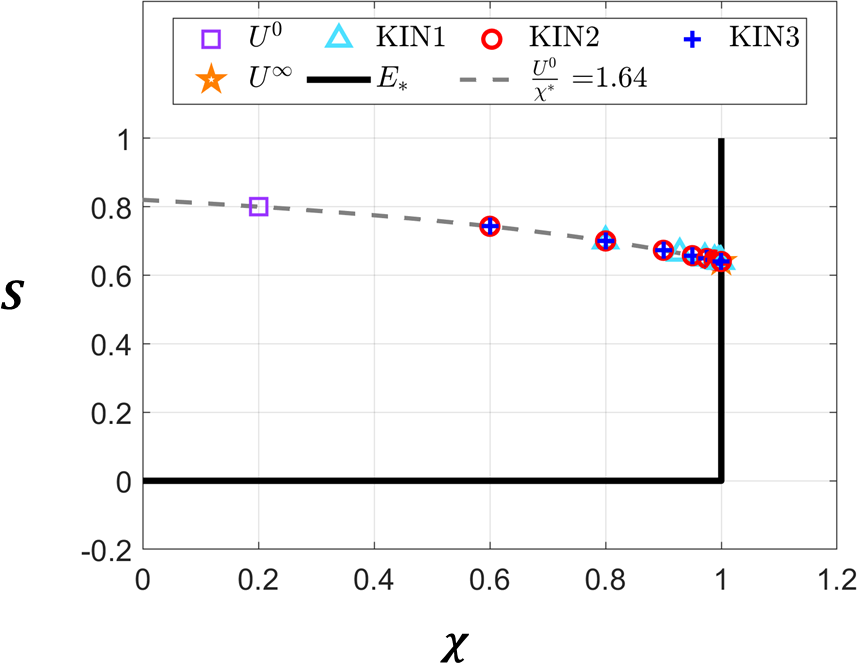}
\includegraphics[width=.49\textwidth]{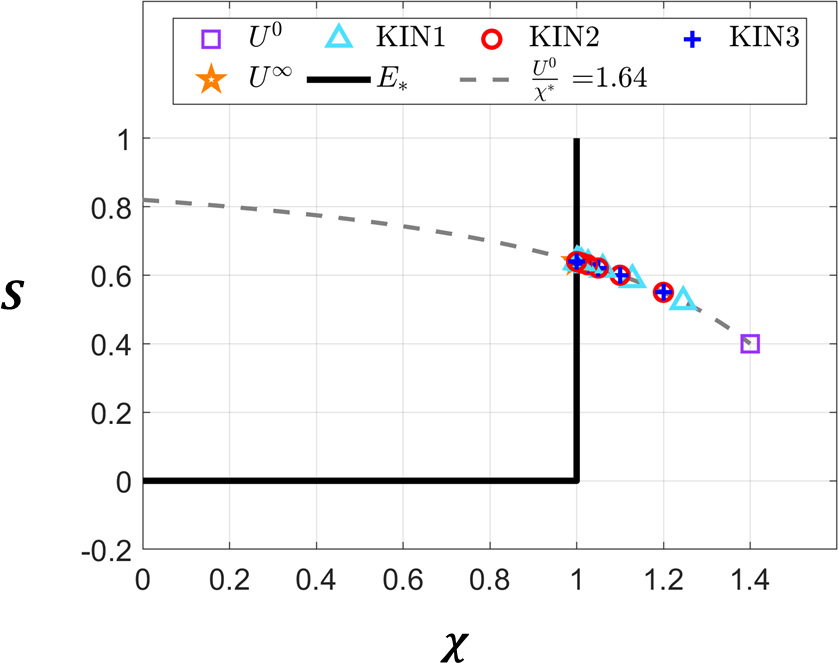}
    \caption{Simulation results for Ex.~\ref{ex:sat} illustrating schemes (KIN1), (KIN2), and (KIN3). Left: $(X^0,S^0) = (0.2,0.8)$. Right: $(X^0,S^0) = (1.4,0.4)$. The solutions to all schemes lie on the curve $u(X^n,S^n)=U^0$ and converge towards the equilibrium point $(X^{\infty},S^{\infty})$ on the portion $\Ep$ of the graph $E$. Solutions to (KIN2) and (KIN3) are indistinguishable.
    \label{fig:ex_saturated}}
\end{figure}

Fig.~\ref{fig:ex_saturated} illustrates the properties of the schemes from Sec.~\ref{sec:ckinetic}. We notice 
first that the property (A) holds: the numerical solutions $(X_j^n,S_j^n)$ given  by \eqref{eq:KIN1}, \eqref{eq:KIN2}, and \eqref{eq:KIN3} live {on the curve} $U^n=u(X^n,S^n)=U^0=u^0 = 1.64$, and as predicted by property (C), they stay in $D^0$ and are physical. 

Second, according to property (D) the solutions to every scheme converge towards the equilibrium point $(X^{\infty},S^{\infty})$ on the portion $\Ep$ of the graph $E_*$, i.e., their distance $Q$ from the equilibrium decreases.  
Third, as predicted by property (E), the solutions to (KIN2) and (KIN3) are indistinguishable, while the solutions to (KIN1) proceed at a rate different than that for (KIN2).

\begin{example}[Unsaturated case]
\label{ex:unsat}
Suppose that $u^0 = 0.6$, thus the equilibrium state $(\chi^{\infty},S^{\infty}) = (0.6,0) \in \Em$. Now we choose $(X^0,S^0) = (0.25,0.2)\in \Dm$. We use large $\tau=1$ or small $\tau=0.5$.
\end{example}
\begin{figure}
    \centering
    \includegraphics[width=.49\textwidth]{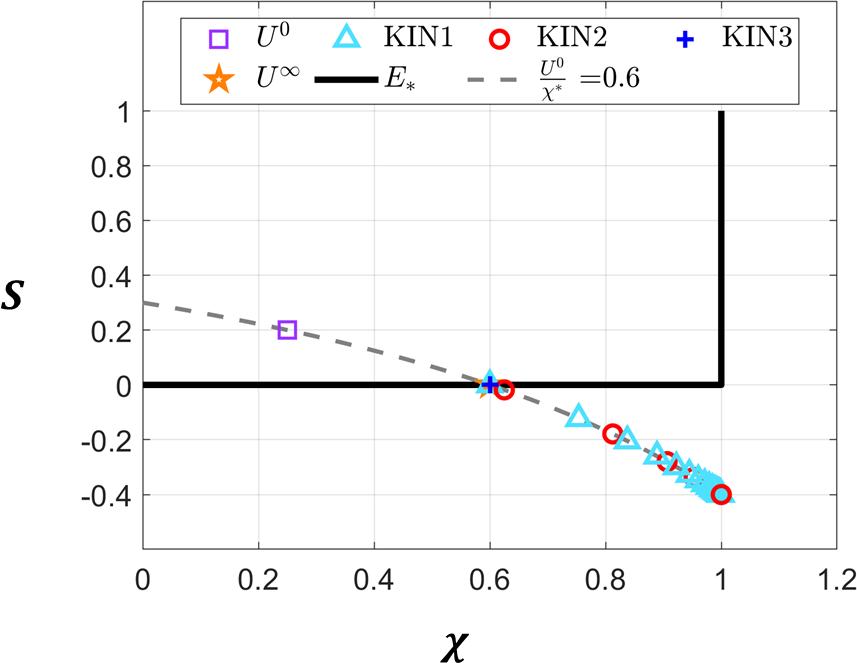}
      \includegraphics[width=.49\textwidth]{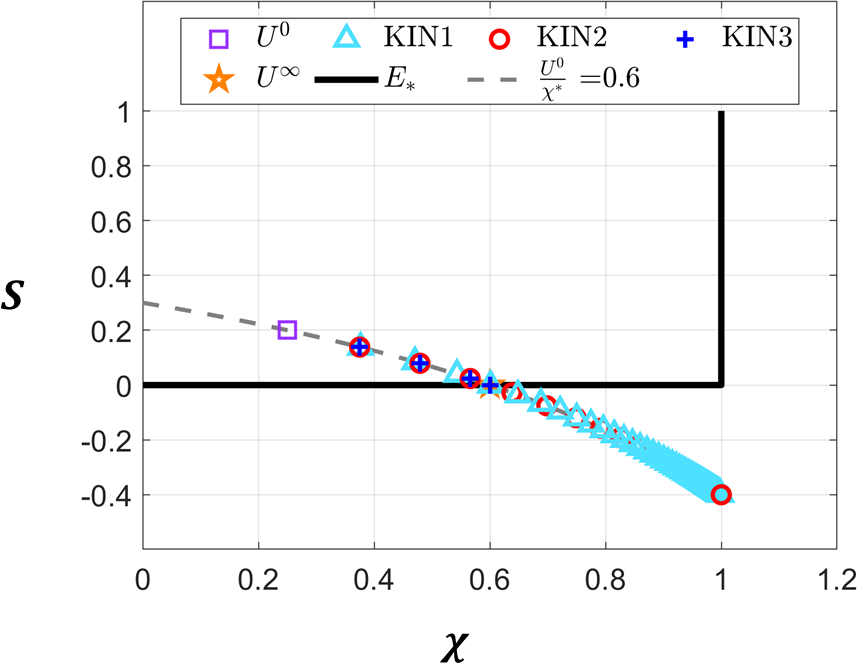}
    \caption{Simulation results of unsaturated case in Ex.~\ref{ex:unsat} illustrating the behavior of kinetic models (KIN1), (KIN2), and (KIN3) with $(X^0,S^0) = (0.25,0.2) \in \Em$. We use $\tau = 1$ (left) and $\tau=0.2$ (right). The solutions to (KIN1) and (KIN2) become unphysical after a few time steps when the curve $u(\chi,S)$ crosses the $\Em$ graph; the solutions to (KIN3) remain in $\Dm$. In addition, while in $D^0$, the solutions to (KIN2) and (KIN3) are indistinguishable.}
    \label{fig:ex_unsaturated}
\end{figure}
Fig.~\ref{fig:ex_unsaturated} demonstrates the results of the three models depending on the time step. First, we see that all solutions live on curve $U^n=u(X^n,S^n)=u^0=0.6$. Second, for smaller $\tau$ we see that (KIN2) and (KIN3) coincide while in $\Dp$. 

However, only  the solutions $(\kn{X}{3}{n},\kn{S}{3}{n})$ to model (KIN3) 
converge to the equilibrium state on $\Em$, and stay physical for all time steps. In contrast, the solutions to (KIN1) and (KIN2), $(\kn{X}{1,2}{n},\kn{S}{1,2}{n})$ give unphysical solutions with negative saturations  $S^n_j<0$, and appear to converge to $X_j^{\infty}=\chi^*$ with {$S_j^{\infty}=-0.4$} for which $Q_j=0$. In fact, $(\kn{X}{1,2}{n},\kn{S}{1,2}{n})$ cross the graph $\Em$, as predicted above. In particular,  for $\tau = 1$, we have $(\kn{X}{1}{n},\kn{S}{1}{n}) \approx (0.7536,-0.1232)$ for $n=2$ and $(\kn{X}{2}{n},\kn{S}{2}{n}) \approx (0.625,-0.0182)$ for $n=1$. For $\tau=0.5$ this happens for larger $n$. 

\begin{figure}
    \centering
    \includegraphics[width=.49\textwidth]{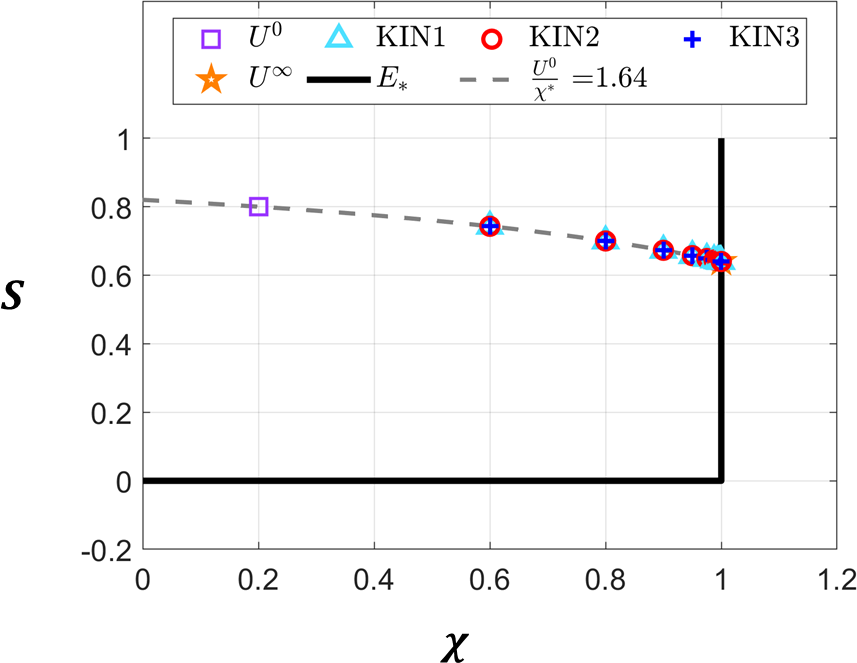}
    \includegraphics[width=.49\textwidth]{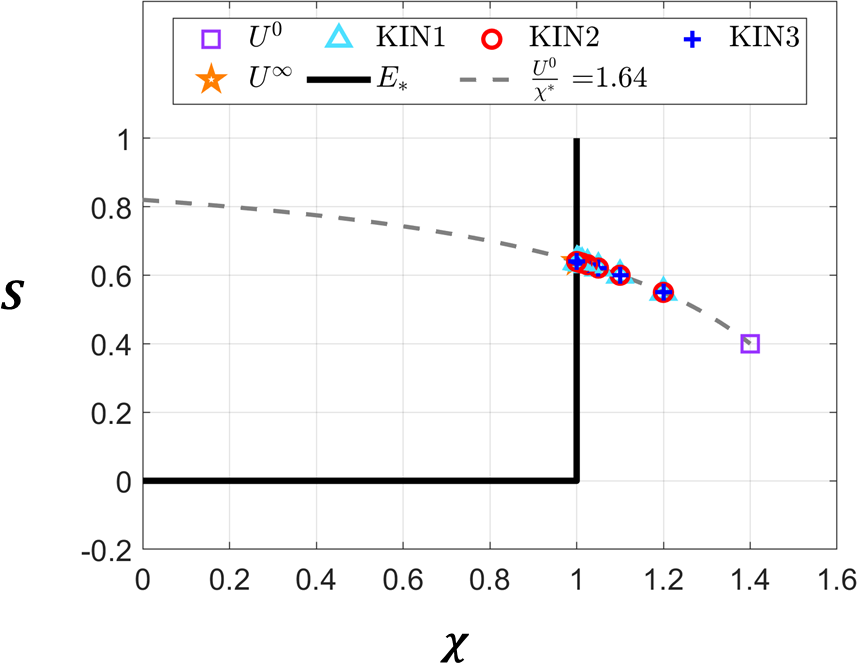}
    \includegraphics[width=.49\textwidth]{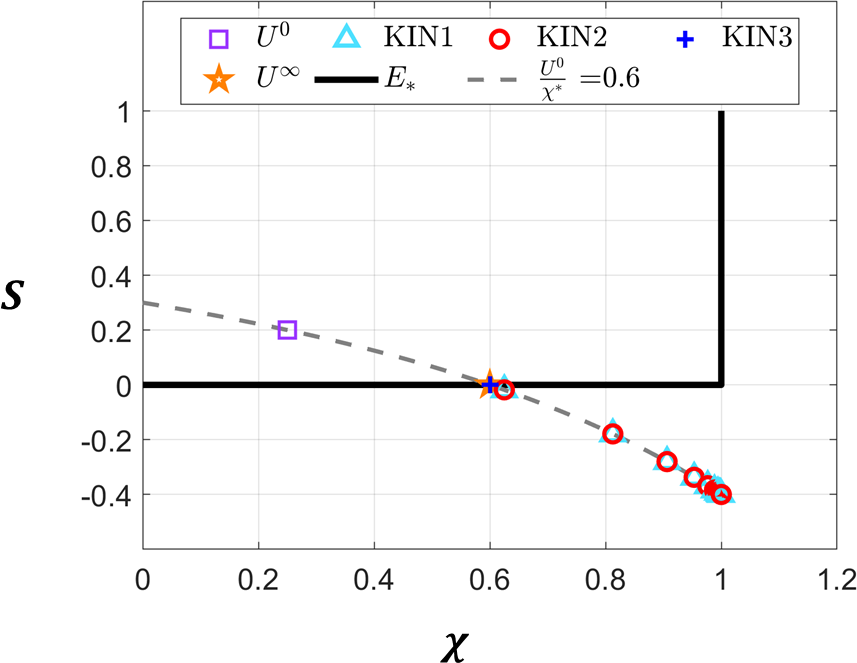}
    \includegraphics[width=.49\textwidth]{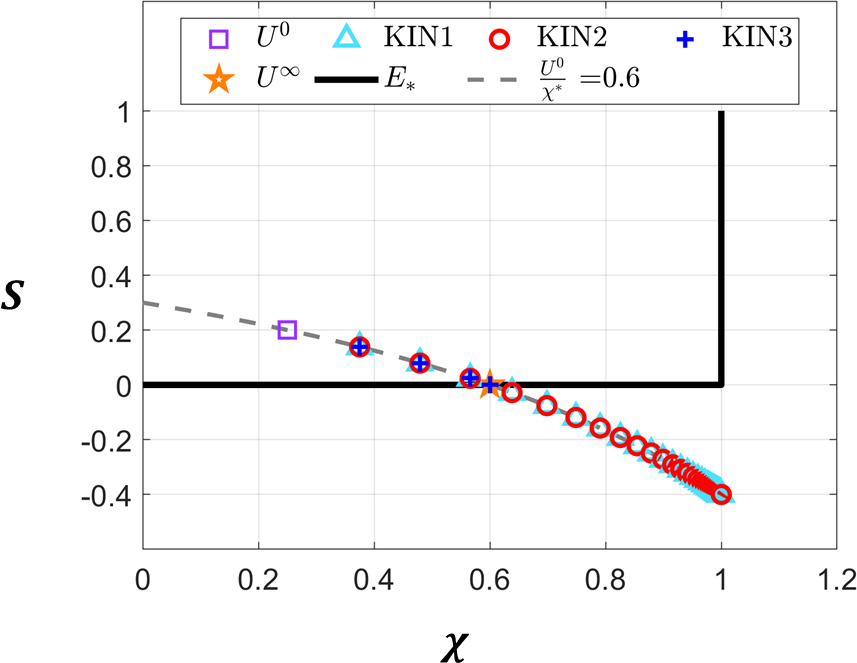}
    \caption{The numerical solutions generated by (KINj) for $j = 1,2,3$ are the same with \eqref{eq:PropE_rate} for $(X^0,S^0)\in D^0_+$. For unsaturated cases, $(\kn{X}{1}{n},\kn{S}{1}{n}) = (\kn{X}{2}{n},\kn{S}{2}{n})$.}
    \label{fig:KIN1_eq_KIN2}
\end{figure}

\begin{example}[Equivalence of (KIN1) and (KIN2)]
\label{ex:equivalent}
In this example we illustrate 
property E.ii. In each case we show that the solutions to (KIN1) are the same as those of (KIN2) when $k_1$ is appropriately chosen depending on $k_2$ and previous time step values.  In turn, (KIN3) solutions are identical to (KIN2) in $\Dp$. See Fig.~\ref{fig:KIN1_eq_KIN2}.
\end{example}

\paragraph{acknowledgements}
{The authors would like to thank the anonymous referees whose comments inspired additional results included in the paper as well as helped to improve the exposition}.
We also thank our colleague Ralph Showalter who made us aware of the paper \cite{HS1995}. We are grateful to our geoscience collaborators Marta Torres and Wei-Li Hong for motivating discussions. 
Choah Shin would like to thank {Larry Martin and Joyce O'Neill for the generous support with the endowed College of Science at Oregon State fellowship 2019-20. Malgorzata Peszynska would like to thank the NSF DMS IRD plan 2019-21 funding which partially made this research possible. {This material is based upon work supported by and while serving at the National Science Foundation. Any opinion, findings, and conclusions or recommendations expressed in this material are those of the authors and do not necessarily reflect the views of the National Science Foundation.}
This research was also partially supported by NSF DMS-1522734 ``Phase transitions in porous media across multiple scales'' and DMS-1912938 ``Modeling with Constraints and Phase Transitions in Porous Media''.}

\bibliographystyle{plain}  
\bibliography{PSpaper,peszynska,mpesz-extra}
\end{document}